\documentclass[12pt,a4paper]{article}

\usepackage{latexsym}
\usepackage{amsfonts}
\usepackage{amssymb}
\usepackage{ifthen}
\usepackage{theorem}
\usepackage{enumerate}
\usepackage{graphics}
\usepackage{color}

\theoremstyle{plain}

\newtheorem{theorem}{Theorem}[section]
\newtheorem{theorem*}[blank]{Theorem}

\newtheorem{conjecture*}[blank]{Conjecture}

\newtheorem{corollary*}[blank]{Corollary}
\newtheorem{lemma}[theorem]{Lemma}
\newtheorem{lemma*}[blank]{Lemma}
\newtheorem{proposition}[theorem]{Proposition}
\newtheorem{proposition*}[blank]{Proposition}

{\theorembodyfont{\rmfamily} 
  \newtheorem{definition}[theorem]{Definition}
  \newtheorem{definition*}[blank]{Definition}
  
  \newtheorem{category*}[blank]{Category}
  
  \newtheorem{functor*}[blank]{Functor}
  \newtheorem{algorithm}[theorem]{Algorithm}
  \newtheorem{algorithm*}[blank]{Algorithm}
  
  \newtheorem{remark*}[blank]{Remark}
  
  \newtheorem{notation*}[blank]{Notation}
  
  \newtheorem{problem*}[blank]{Problem}
  
  \newtheorem{question*}[blank]{Question}
  
  \newtheorem{terminology*}[blank]{Terminology}
  
  \theoremstyle{break}
    \newtheorem{example}[theorem]{Example}
    \newtheorem{example*}[blank]{Example}
}

\newcommand{\isz}{\setlength{\itemsep}{0pt}}
\newcommand{\ifnull}[3]{\def\nullstring{}\def\teststring{#1}\ifx\nullstring\teststring#2\else#3\fi}
\newenvironment{proof}[1][]{\begin{trivlist}\item\ifnull{#1}{\noindent{\sc Proof:\ }}{\noindent{\sc Proof of #1:\ }}}{\hspace*{\fill}$\Box$\vspace*{\baselineskip}\end{trivlist}}
\newcounter{enum}

\newcommand{\steplabel}[1]{\label{step:#1}}

\newcommand{\stepref}[1]{step \ref{step:#1}}

\newcommand{\figref}[1]{Figure~\ref{fig:#1}}

\newcommand{\figrefpart}[2]{Figure~\ref{fig:#1}(#2)}

\newcommand{\figscale}{1}

\newcommand{\fig}[2]{
  \begin{figure}[htb]\centering
    \ifnull{#1}{}{\scalebox{\figscale}{\includegraphics{#1.eps}}}
    \ifnull{#2}{}{\caption{#2}}
    \label{fig:#1}
  \end{figure}
}

\newcommand{\Z}{\mathbb{Z}}
\newcommand{\N}{{\mathbb N}}

\renewcommand{\th}{\ensuremath{^\mathrm{th}}}

\renewcommand{\emptyset}{\varnothing}
\renewcommand{\geq}{\geqslant}
\renewcommand{\leq}{\leqslant}

\newcommand{\card}[1]{\#(#1)}
\newcommand{\dirlim}[1]{\ifthenelse{\equal{#1}{}}{\underrightarrow{\lim}}{\underrightarrow{#1}}}
\newcommand{\embed}{\hookrightarrow}
\newcommand{\fto}{\longrightarrow}
\newcommand{\exto}{\fto}

\newcommand{\homotopic}{\sim}

\newcommand{\hniel}{h_\mathit{niel}}
\newcommand{\htop}{h_\mathit{top}}

\newcommand{\len}{{l}}

\newcommand{\p}{^\prime}
\newcommand{\pp}{^{\prime\prime}}
\newcommand{\restrict}[1]{|_{#1}}

\newcommand{\tendsto}{\rightarrow}

\newcommand{\norm}[2][]{|{#2}|_{#1}}

\newcommand{\id}{id}

\newcommand{\mod}{\mathrm{mod}}

\newcommand{\cl}[1]{\mathrm{cl}(#1)}
\newcommand{\interior}[1]{\mathrm{int}(#1)}
\newcommand{\inr}[1]{\interior{#1}}

\newcommand{\Alpha}{\aleph}

\newcommand{\link}[1]{\mathrm{Lk}(#1)}

\newcommand{\pre}[1]{\textrm{Pre-}#1}

\newcommand{\tang}{\makebox[0pt][l]{$\frown$}\raisebox{0.5ex}{$\smile$}}
\newcommand{\trans}{\makebox{\Large $\times$}}
\newcommand{\sig}{\mathcal{B}}

\newcommand{\us}{{U/S}}

\newcommand{\edge}{E}
\newcommand{\control}{C}
\newcommand{\vertex}{V}
\newcommand{\init}[1]{\imath(#1)}

\newcommand{\final}[1]{\imath(\bar #1)}

\renewcommand{\card}[1]{|#1|}
\newcommand{\cut}[2][]{\mathcal{C}^{#1}{#2}}

\newcommand{\normalmargins}{
  \setlength{\oddsidemargin}{-0.4mm}
  \setlength{\evensidemargin}{-0.4mm}
  \setlength{\textwidth}{160mm}

  \setlength{\topmargin}{-0.4mm}
  \setlength{\headheight}{4mm}
  \setlength{\headsep}{7mm}
  \setlength{\textheight}{225mm}
  \setlength{\footskip}{11mm}

  \setlength{\parindent}{0mm}
  \setlength{\parskip}{0.5\baselineskip}
  \addtolength{\topsep}{-\parskip}
  \addtolength{\partopsep}{-\parskip}

  \setlength{\headheight}{0mm}
  \setlength{\headsep}{0mm}
  \setlength{\textheight}{236mm}
}

\renewcommand{\figscale}{1}
\normalmargins

\begin{document}
\sloppy


\title{Dynamics of surface diffeomorphisms relative to homoclinic and heteroclinic orbits}
\author{Pieter Collins
  \thanks{This work was partially funded by Leverhulme Special Research Fellowship SRF/4/9900172}
  \\Department of Mathematical Sciences\\University of Liverpool\\Liverpool L69 7ZL, U.K.\\pieter.collins@12move.nl
 }
\date{\today}
\maketitle

\begin{abstract}
In this paper, we show how to obtain information about the dynamics of a two-dimensional discrete-time system from its homoclinic and heteroclinic orbits.
The results obtained are based on the theory of \emph{trellises}, which comprise finite-length subsets of the stable and unstable manifolds of a collection of saddle periodic orbits.
For any collection of homoclinic or heteroclinic orbits, we show how to associate a canonical \emph{trellis type} which describes the orbits.
Given a trellis type, we show how to compute a \emph{graph representative} which gives a combinatorial invariant of the trellis type.
The orbits of the graph representative represent the dynamics forced by the homoclinic/heteroclinic orbits in the sense that every orbit of the graph representative is \emph{globally shadowed} by some orbit of the system, and periodic, homoclinic/heteroclinic orbits of the graph representative are shadowed by similar orbits.
By constructing suitable surface diffeomorphisms, we prove that these results are optimal in the sense that the topological entropy of the graph representative is the infemum of the topological entropies in the isotopy class relative to the homoclinic/heteroclinic orbits.
\end{abstract}


\newpage

{\setlength{\parskip}{0pt} \tableofcontents}


\section{Introduction}
\label{sec:Introduction}

The importance of homoclinic orbits in dynamical systems theory was first realised by Poincar\'e \cite{Poincare1892}, who showed that the presence of a homoclinic tangle in the three-body problem was enough to show that the system was non-integrable.
More work on homoclinic tangles was later undertaken by Birkhoff \cite{Birkhoff49}, but modern interest was inspired by Smale, who showed that any system with a transverse homoclinic point must have a horseshoe in some iterate \cite{Smale63}.
In this paper, we extend the local analysis of Smale to obtain a global understanding of the dynamics for diffeomorphisms in two dimensions.
The analysis depends on constructing a \emph{graph representative} of the system, from which we obtain information about the dynamics in terms of \emph{symbolic dynamics} relative to a Markov partition and \emph{global shadowing}.

The use of Markov partitions to give symbolic dynamics is a fundamental tool in dynamical systems theory. 
For uniformly-hyperbolic systems, the dynamics is described completely by the symbolics, but for non-uniformly hyperbolic systems, the Markov partition typically only gives a \emph{lower bound} for the dynamics.
From this, we can often deduce information on the \emph{forcing relation}, which describes when the presence of a particular geometric structure implies the existence of another structure.
Sharkovskii's theorem for interval maps is a classic example of a forcing result: a total ordering $\succ$ of the natural numbers is given such that the presence of a periodic orbit of period $p$ forces periodic orbits of all periods $q$ with $p\succ q$.
Sharkovskii's theorem can be strengthened to consider the ordering of points on the orbit, and this gives a complete (though less elegant) description of the forcing relation for periodic orbits of interval maps.
Interval maps can be studied by elementary methods, but obtaining symbolic dynamics for more complicated classes of system typically requires the use of some topological index theory, such as the fixed point index \cite{Fried83,Jiang93} or the Conley index \cite{Conley78}.

The theory of isotopy classes of surface diffeomorphisms relative to periodic orbits is known as \emph{Nielsen-Thurston theory} \cite{CassonBleiler88}.
Here, the period is no longer a useful characterisation of the orbit; instead, we need to consider the \emph{braid type} which is the conjugacy class of isotopy classes relative to the orbit in question.
Just as for periodic orbits of interval maps, orbits of different periods may have different braid types.
For each braid type there is a canonical representative such that a braid type $BT_1$ forces $BT_2$ if and only if the \emph{Thurston minimal representative} of $BT_1$ contains an orbit of type $BT_2$.
A Thurston minimal representative has the minimum number of periodic points of any period, and minimal topological entropy in the isotopy class \cite{Handel85}.
The most interesting case is when the canonical representative is pseudo-Anosov.
A pseudo-Anosov map preserves a pair of transverse foliations, the leaves of the stable foliation being uniformly contracted and the unstable foliation being uniformly expanded.
Pseudo-Anosov maps therefore exhibit a type of uniformly hyperbolic behaviour, and Thurston's classification theorem can be seen as a persistence result for this hyperbolicity.

Unfortunately, Thurston's original proof was non-constructive, and thus of little use in actually determining the forcing relation.
It was not until an algorithm of Los \cite{Los93} that computations became possible.
Los' algorithm was soon followed by an algorithm which was discovered independently by Franks and Misiurewicz \cite{FranksMisiurewicz93} and Bestvina and Handel \cite{BestvinaHandel95}.
A description of the moves of this algorithm with an emphasis on implementation details was given by Keil \cite{Keil97}.
Even so, the global properties of the forcing relation are still only partially understood and this is an active area of research \cite{deCarvalhoHallPP4}.

A more challenging problem is to consider the dynamics of surface diffeomorphisms relative to homoclinic and heteroclinic orbits of tangles.
The standard results of Poincar\'e and Smale for transverse homoclinic points have since been extended to topologically transverse intersections which are not geometrically transverse \cite{BurnsWeiss95COMMP,Rayskin98PhD}, and a great deal is known about bifurcations near a homoclinic tangency \cite{Newhouse80,MoraViana93ACTAM}.
The global dynamics and the forcing relation have been previously considered by McRobie and Thompson \cite{McRobieThompson94DSTSYS}, who carried out a detailed study of the \emph{spiral $3$-shoe}, by Rom-Kedar \cite{RomKedar94NONLIN}, who obtained results on entropy bounds and homoclinic forcing using a \emph{topological approximation method}, and by Easton \cite{Easton00}.
Algorithms for computing dynamics forced by homoclinic orbits were recently given independently by Handel \cite{Handel99TOPOL} and Hulme \cite{Hulme00PHD}.
Homoclinic tangles have been used by Christiansen and Politi \cite{ChristiansenPoliti96NONLIN} to compute generating partitions for symbolic dynamics.
Orbits of the H\'enon map have also been obtained using continuation methods by Sterling, Dullin and Meiss \cite{SterlingDullinMeiss99PHYSD}.

The aim of this paper is to give a complete and general theoretical basis for describing and computing the dynamics forced by a homoclinic or heteroclinic tangle.
The results are a completion of the theory begun in \cite{Collins99AMS}, and are related to the pruning theory of de Carvalho \cite{deCarvalho99}.
Instead of considering an entire tangle, which consists of immersed curves of infinite length, we only consider a subset consisting of curves of finite length.
We call such a figure a \emph{trellis}, after Birkhoff \cite{Birkhoff49} and Easton \cite{Easton86TRAAM}.
Aside from being easier to work with, trellises have the advantage that there exist a number of algorithms to compute them to a very high degree of accuracy, such as those of Simo \cite{Simo89} and Krauskopf and Osinga \cite{KrauskopfOsinga98JCOMP}.
We therefore obtain results which can be directly applied to physically-relevant systems.
In Section~\ref{sec:tangletrellis} we give a formal description of trellises and discuss the features of their geometry which we shall need later.

In Section~\ref{sec:biasymptoticminimal}, we we show how to relate a collection of homoclinic and/or heteroclinic orbits to a canonical trellis type forcing the same dynamics, thus unifying the viewpoints of homoclinic/heteroclinic orbits and tangles.
This result can be applied to classify homoclinic orbits by braid type.
In particular, is is straightforward to obtain a classification of \emph{horseshoe homoclinic orbits}, which are homoclinic orbits of the Smale horseshoe map.

For the rest of the paper, we try to mimic the results of Nielsen-Thurston theory as closely as possible.
While the fact that we are dealing with homoclinic orbits often complicates matters, it sometimes can yield simplifications.
In Section~\ref{sec:graph}, we show how to associate a graph representative to a given trellis map, in much the same way that a pseudo-Anosov braid type gives rise to a train track.
The graph we obtain is unique (unlike a train track) and gives a convenient way of specifying the trellis.
From this, we also obtain a canonical graph representative for a homoclinic orbit.

We consider the problem of computing the graph representative in Section~\ref{sec:algorithm}.
We give and algorithm modelled on the Bestvina-Handel algorithm, but again the invariant manifolds simplify matters by enabling us to avoid the notorious valence-$2$ homotopies of the graph map.
We follow an idea of Los \cite{Los96} by replacing the valence-$2$ homotopy by a \emph{valence-$3$ homotopy} which does not require a global knowledge of the map.
Unfortunately, Los' analysis of the effect of the move on entropy was erroneous; the entropy may increase during this move, and this meant that the algorithm was not proved to terminate.
Our analysis goes along different lines and avoids entropy computations altogether.
We also unify the process of absorbing into the peripheral subgraph with the main algorithm, which also simplifies matters somewhat.
As with the Bestvina-Handel algorithm, reductions are obtained and dealt with accordingly.
Finally, we give a simplified algorithm for braid types based on our analysis of the algorithm for trellises.

In Section~\ref{sec:shadowing}, we discuss how the dynamics of the graph representative model the dynamics of the surface diffeomorphism.
These results are based on the \emph{relative periodic point theory} developed in \cite{Collins01TOPOA}.
We first deduce the existence of periodic orbits; the existence of other orbits follows by taking limits.
In particular, these results show that the topological entropy of the graph representative is a lower bound for the topological entropy of the surface diffeomorphism.
This lower bound is called the \emph{Nielsen entropy}, since it is related to the growth rate of the number of essential Nielsen classes of period-$n$ orbits.
Orbits are related by \emph{global shadowing}, a relation introduced by Katok, though here we give a new definition more in keeping with the spirit of our topological analysis.
One important fact is that periodic orbits which globally shadow each other have the same braid type, so we obtain information about the braid types of the orbits of the system.
Finally, we show be a careful construction that homoclinic and heteroclinic orbits of the graph map are shadowed by homoclinic and heteroclinic orbits of the surface diffeomorphism.

The final section, Section~\ref{sec:entropy}, is devoted to showing that the dynamics computed by our algorithm is indeed all the dynamics forced by the trellis.
In particular, the lower entropy bound computed by the algorithm is sharp.
Unfortunately, unlike Nielsen-Thurston theory, there may not be a canonical diffeomorphism whose topological entropy equals the entropy bound given by the Nielsen entropy.
We give necessary and sufficient conditions for the existence of a diffeomorphism whose entropy equals the Nielsen entropy.
We also show that by carefully extending the original trellis, we can always find a diffeomorphism whose entropy approximates the Nielsen entropy arbitrarily closely.
Finally, we complete the link with Nielsen-Thurston theory by finding pseudo-Anosov maps which approximate the Nielsen entropy when such maps exist for the given trellis type.

A unifying theme throughout this paper is to consider isotopy classes of curves embedded in the surface.
We are particularly concerned with curves that have the minimum possible number of intersections with the stable curves of the trellis.
Stable curves can be thought of as cohomology classes, though we do not use this interpretation here, as the abelian homology theory is insufficient for many of our constructions, which require the non-abelian homotopy groups.
Many of the results we need are intuitively obvious ``folk theorems,'' but we give complete proofs in Appendix~\ref{sec:isotopytheory}.


\section{Tangles and Trellises}
\label{sec:tangletrellis}

In this section we discuss the geometry of homoclinic and heteroclinic tangles, and introduce trellises, trellis mapping classes and trellis types, which are our main topic of study.
We give a number of definitions which allow us to describe some of the properties of trellises.
We discuss reducibility of trellis mapping classes, and explain how this differs from reducibility of mapping classes of surface diffeomorphisms.
Finally, we discuss the the procedure of cutting along the unstable set, which gives us a way of relating trellis maps to maps of topological pairs, 
 and hence a way of studying the dynamics.


\subsection{Tangles}
\label{sec:tangle}

A tangle is the figure formed by the stable and unstable manifolds of a collection of periodic saddle orbits.

\begin{definition}[Tangle]
Let $f$ be a diffeomorphism of a surface $M$, and $P$ a finite invariant set of saddle points of $f$.
Let $W^S(f;P)$ be the stable set of $P$, and $W^U(f;P)$ the unstable set of $P$.
The pair $W=(W^U(f;P),W^S(f;P))$ is a \emph{tangle} for $f$.
If $P$ consists of a single point, then $W$ is a \emph{homoclinic} tangle, otherwise $W$ is a \emph{heteroclinic} tangle.

If $p$ is a point of $P$, then the stable and unstable curves passing through $p$ are denoted $W^S(f;p)$ and $W^U(f;p)$ respectively.
\end{definition}

We use the notation $W^{\us}$ in statements which hold for both stable and unstable manifolds.
If the diffeomorphism $f$ under consideration is clear from the context, we abbreviate $W^{\us}(f;P)$ to $W^{\us}(P)$,
 and if the periodic point set $P$ is also clear, we simply write $W^{\us}$.

We use the obvious notation for intervals of $W^{\us}$:
The closed interval in $W^{\us}$ with endpoints $a$ and $b$ is denoted $W^{\us}[a,b]$ and the open interval is denoted $W^{\us}(a,b)$.
Half-open intervals are denoted $W^{\us}[a,b)$ and $W^{\us}(a,b]$.

Each stable or unstable curve $W^{\us}(p)$ is naturally separated into two pieces by the periodic point $p$.
\begin{definition}[Branch]
A \emph{branch} of a tangle is a component of $W^{\us}(P)\setminus P$.
We orient each branch of $W^S(p)$ so that the positive orientation is in the direction of $p$, and each branch of $W^U(p)$ so the the positive orientation is away from $p$.
A branch is \emph{trivial} if it does not intersect any other branch.
\end{definition}
We denote a branch of $W^\us(p)$ by integers $W^\us(p,b)$.
Note that the branches of $W^\us$ are permuted by $f$.

\begin{definition}[Fundamental domain]
A \emph{fundamental domain} is a half-open interval in $W^\us$ bounded by points $q$ and $f^n(q)$, where $q\not\in W^P$,
 and $n$ is the least integer such that $q$ and $f^n(q)$ lie in the same branch.
\end{definition}
A branch of a tangle can be found by taking the union of the $n$-th iterates of a fundamental domain.

\begin{definition}[Biasymptotic points]
An \emph{intersection point} of a tangle $W(P)$ is a a point in $W^U\cap W^S$.
An intersection point not in $P$ is a \emph{biasymptotic} point; an \emph{asymptotic} point is a point in $(W^U\cup W^S)\setminus P$.
A biasymptotic point $q\in W^U(p_1)\cap W^S(p_2)$ is a \emph{homoclinic} point if $p_1$ and $p_2$ are part of the same periodic orbit,
 otherwise $q$ is a \emph{heteroclinic} point.
If $W^U$ and $W^S$ cross transversely at a biasymptotic point $q$, then $q$ is a \emph{transverse} homoclinic or heteroclinic point.
\end{definition}
There are six different types of isolated intersection point, two transverse intersections $\trans$, with positive ($+$) and negative ($-$) crossings respectively, 
 and are four tangential intersections $\tang$.
It is well known \cite{Smale63} that if a diffeomorphism $f$ has a transverse homoclinic point, then $f$ must be chaotic and have positive topological entropy.
If $f$ is a uniformly hyperbolic diffeomorphism, then any biasymptotic point is a transverse intersection.
\begin{definition}[Primary intersection point]
An intersection point $q\in W^U(p_1)\cap W^U(p_2)$ is a \emph{primary intersection point} or \emph{pip} if 
 $W^U(p_1,q)$ and $W^S(q,p_2)$ are disjoint.
\end{definition}


\subsection{Trellises}
\label{sec:trellis}

There are a number of difficulties in working with tangles directly.
These are due to the fact that if any branch intersects another, then it will be an immersed, rather then embedded, curve.
The geometry of such a curve is extremely complicated, as is the global topology of the tangle.
Further, since branches may have infinite length, tangles are impossible to compute directly.
Instead, we work with subsets of tangles known as \emph{trellises}, for which the stable and unstable curves have finite lengths.

\begin{definition}[Trellis]
Let $f$ be a diffeomorphism of a surface $M$ with a finite invariant set $P$ of hyperbolic saddle points.
A \emph{trellis} for $f$ is a pair $T=(T^U,T^S)$, where $T^U$ and $T^S$ be subsets of $W^U(f;P)$ and $W^S(f;P)$ respectively such that 
{\setlength{\parskip}{0pt}
\begin{enumerate}\isz
\item $T^U$ and $T^S$ both consist of finitely many compact intervals with non-empty interiors, and
\item $f(T^U)\supset T^U$ and $f(T^S)\subset T^S$.
\end{enumerate}
We denote the set of periodic points of $T$ by $T^P$, and the set of intersections of $T^U$ and $T^S$ by $T^V$.
}
\end{definition}
\begin{definition}[Trellis map]
If $f$ is a diffeomorphism and $T$ is a trellis for $f$, then the pair $(f;T)$ is called a \emph{trellis map}.
\end{definition}
There are two important relations between trellis maps, namely \emph{isotopy} and \emph{conjugacy}.
\begin{definition}[Isotopy and conjugacy]
Trellis maps $(f_0;T)$ and $(f_1;T)$ are \emph{isotopic} if there is an isotopy of surface diffeomorphisms $(f_t)$ such that for each $t$, $(f_t;T)$ is a trellis map.
Trellis maps $(f_0;T_0)$ and $(f_1;T_1)$ are \emph{conjugate} if there is a homeomorphism $h$ such that $h\circ f_0=f_1\circ h$ and $h(T_0^{\us})=T_1^{\us}$.
\end{definition}
These these relations allow us to define equivalence classes of trellis map. 
It is these equivalence classes which will be our primary object of study.
\begin{definition}[Trellis mapping class and trellis type]
The \emph{trellis mapping class} $([f];T)$ is the set of all trellis maps which are isotopic to $f$.
The \emph{trellis type} $[f;T]$ is the set of all trellis maps which are conjugate to a map in $([f];T)$.
Trellis maps $(f_0;T_0)$ and $(f_1;T_1)$ are considered equivalent if $[f_0;T_0]=[f_1;T_1]$.
\end{definition}

\fig{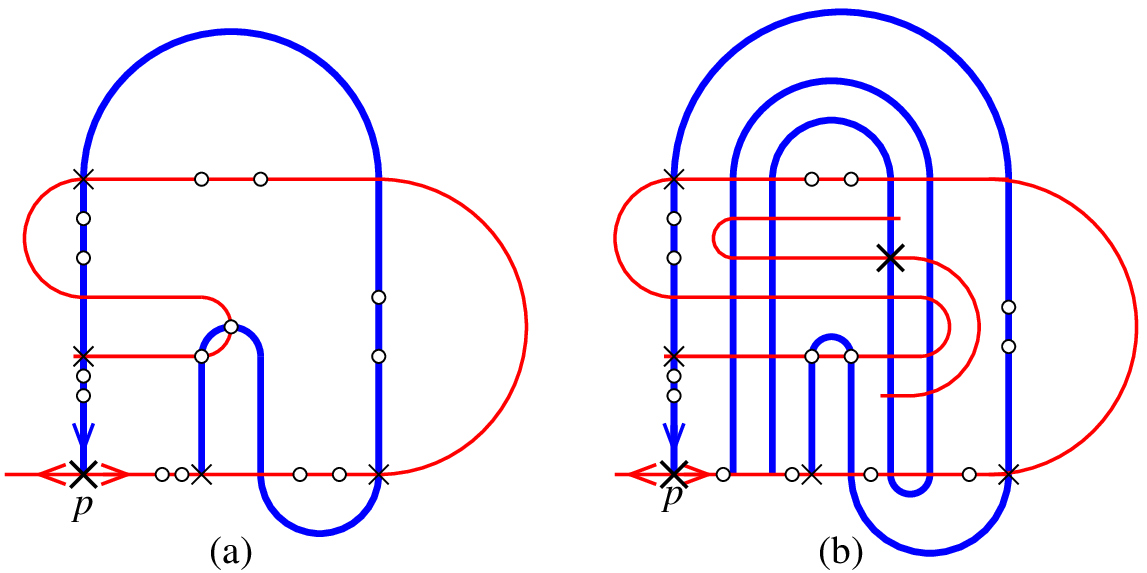}{Trellises in the H\'enon map. (a) is a homoclinic trellis, (b) a heteroclinic trellis.}
In \figref{henontrellises} we show two trellises, each of which is homeomorphic to a trellis for a diffeomorphism in the H\'enon family.
In each case the small crosses mark primary intersection points on a single orbit.
It turns out that this is enough information to specify a trellis type for these trellises.

Just as we can define branches of a tangle, we can also define branches of a trellis.
\begin{definition}[Branch]
If $W^{\us}(p,b_{u/s})$ is branch of a tangle for $f$, then the corresponding branch of a trellis $T\subset W$ is given by 
\[ T^{\us}(p,b_{u/s})=\{p\}\{p\} \cup W^{\us}(p,b_{u/s}) \cap T^{\us} . \] 
\end{definition}
Note that a branch of a trellis is a union of embedded closed intervals, whereas a branch of a tangle is an immersed open interval.
In particular, a branch of a trellis contains the periodic saddle point.

Most numerical algorithms compute trellises by growing out the branches from $P$, and the resulting trellises have connected branches.
While this is the most natural case to visualise, it is often important to consider trellises with branches that are not connected.
Disconnected branches can be useful when considering trellises which are subsets of tangles with unbounded branches.

The most important subsets of $T^\us$ are intervals with endpoints in $T^V$, especially those which have no intersections in their interior.
\begin{definition}[Arc, segment and end]
An \emph{arc} of $T^\us$ is a closed subinterval of $T^\us$ with endpoints in $T^V$.
A \emph{segment} of $T^\us$ is an arc of $T^\us$ with no topologically transverse intersection points in its interior.
An \emph{end} of a trellis is a subintervals of $T^\us$ which does not lie in any segment.
An \emph{end intersection} of a trellis is an intersection which does not lie in two segments of $T^\us$ (so lies in at least one end of $T^\us$).
\end{definition}

Similarly, the most important subsets of the surface $M$ are those bounded by $T^U$ and $T^S$.
\begin{definition}[Domain and region]
An \emph{open domain} of $T$ is an open subset of $M$ with boundary in $T^U\cup T^S$.
A \emph{closed domain} is the closure of an open domain.
An \emph{open region} of $T$ is a component of $M\setminus(T^U\cup T^S)$. 
A \emph{closed region} is the closure of an open region, and hence includes the stable and unstable boundary segments.
We will denote domains by the letter $D$, and regions by the letter $R$.
\end{definition}
There are two special types of region which play an important role later.
\begin{definition}[Rectangle]
A \emph{rectangle} is a region which is a topological disc bounded by two stable and two unstable segments which make acute or obtuse angles with each other
 (i.e. do not form reflex angles).
\end{definition}
\begin{definition}[Bigon]
A \emph{bigon} is a region which is a topological disc bounded by one stable and one unstable segment which make acute or obtuse angles with each other.
\end{definition}
Bigons play a similar role in the trellis theory as do critical points in the kneading theory, and may be called \emph{critical} regions.

In general, we do not have much control over the geometry of a trellis.
The exception is near a point of $T^P$, where the dynamics are conjugate to a linear map.
It is important to consider how the nontrivial branches $T$ divide the surface in a neighbourhood of $T^P$.
\begin{definition}[Quadrant, secant and coquadrant]
Let $T$ be a trellis, and $\overline{T}$ be the nontrivial branches of $T$, and $p$ be a point of $T^P$.
Then a \emph{sector} of $T$ at $p$ is a local component of a region of $\overline{T}$ in a neighbourhood of $p$.
A sector is a \emph{quadrant} $Q$ if it intersects no trivial branches, in which case the boundary includes a single unstable branch $T^U(Q)$ and a single stable branch $T^S(Q)$.
A sector is an \emph{secant} if it intersects a single trivial branch; if this is a branch of $T^U$ the secant is \emph{attracting}, and if this is a branch of $T^S$, the secant is repelling.
A sector is a \emph{coquadrant} if it intersects two trivial branches, in which case the only other sector at $p$ is a quadrant.
\end{definition}

\fig{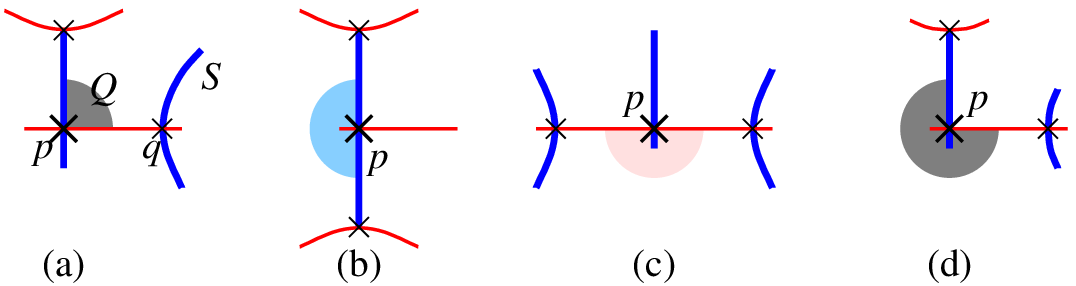}{Sectors at $p$. (a) is a quadrant, (b) an attracting secant, (c) a repelling secant and (d) a coquadrant.}

Quadrants, secants and coquadrants are shown in \figref{quadrantsecant}.
The region containing a quadrant $Q$ is denoted $R(Q)$.
The \emph{image} of a quadrant $Q$ is the quadrant containing $f(Q)$; note that $T^\us(f(Q))$ is the same branch as $f(T^\us(Q))$.
A stable segment $S$ with endpoint $q\in T^U(Q)$ intersects lies on the \emph{$Q$-side} of $T^U$ if locally $S$ lies on the same side of $T^U(Q)$ at $q$ as $T^S(Q)$ does at $p$.
Similarly, an unstable segment $U$ with endpoint $q\in T^S(Q)$ intersects lies on the \emph{$Q$-side} of $T^S$ if locally $U$ lies on the same side of $T^S(Q)$ at $q$ as $T^U(Q)$ does at $p$.
The segment $S$ in \figrefpart{quadrantsecant}{a} lies on the $Q$-side of $q$.

We now illustrate the concepts introduced so far with one of the most important trellis types, the \emph{Smale horseshoe trellis}.
\begin{example*}[The Smale horseshoe trellis]

\fig{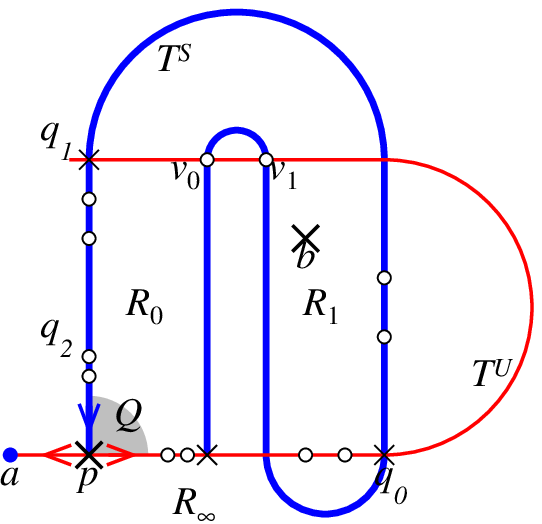}{The Smale horseshoe trellis.}

The Smale horseshoe trellis is formed by the stable and unstable manifolds of the direct saddle fixed point (i.e. the saddle point with positive eigenvalues)
 of the Smale horseshoe map, and is depicted in \figref{smale}.
The stable and unstable sets are subsets of the stable and unstable manifolds of the saddle fixed point $p$.
The branches of this trellis are connected and the all intersection points are transverse.
The points $q_0$, $q_1$ and $q_2$ are primary intersection points on a single homoclinic orbit.
The orbits of the biasymptotic points $v_0$ and $v_1$ are shown in white dots.
These orbits are called \emph{forcing orbits} for the Smale horseshoe trellis; 
 it is essentially impossible to remove any intersection points of the Smale horseshoe trellis by an isotopy
 without first removing $v_0$ and $v_1$ in a homoclinic bifurcation.

Throughout this paper we follow the convention of showing forcing orbits as which dots, and primary intersection points as crosses.

There are eight regions, an unbounded region $R_\infty$, three bigons and four rectangles.
The quadrant $Q$ is contained in the region $R_0$.
Under the Smale horseshoe map, there is a Cantor set of nonwandering points contained in the (closed) rectangular regions $R_0$ and $R_1$, including a fixed point $b$ in $R_1$.
All other points are wandering except for an attracting fixed point at the end of one unstable branch of $T^U(p)$
 in $R_\infty$.
The topological entropy of the Smale horseshoe map is $\log2$.
We shall see later that any diffeomorphism with this trellis type must have topological entropy $\htop\geq\log2$.

Notice that one of the unstable branches of the Smale horseshoe trellis ends in an attracting fixed point,
 and one of the stable branches is non-existent.
Both of these are therefore trivial branches.

\fig{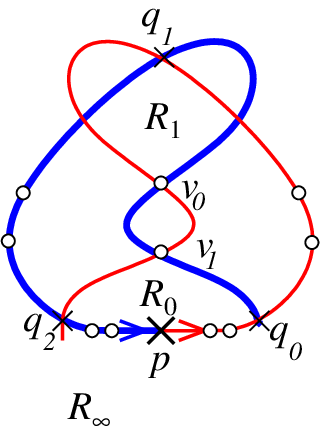}{The singular Smale horseshoe trellis.}

We could isotope the Smale horseshoe map in a neighbourhood of the saddle point $p$ to obtain a homeomorphism with a one-prong singularity at $p$.
The new trellis is shown in \figref{singularsmale}.
It is topologically equivalent to the Smale horseshoe trellis (minus the unstable branch ending at $a$) but has a different geometric type,
 and is known as the \emph{singular Smale horseshoe trellis} 
Many of the properties of the singular Smale horseshoe trellis are the same as that of regular Smale horseshoe trellis,
 except that now there is no homeomorphism in the trellis mapping class with a one-prong singularity at $p$ which has topological entropy $\log2$.
Instead, for any $\epsilon>0$ there is a pseudo-Anosov homeomorphism $f$ with topological entropy $\log2<\htop(f)<\log2+\epsilon$.
\end{example*}

The most important dynamical feature of a trellis type is its \emph{entropy}.
\begin{definition}[Entropy]
The \emph{entropy} of a trellis type $[f;T]$, denoted $\htop[f;T]$ is the infemum of the topological entropies of diffeomorphisms in $[f;T]$; that is
\[ \htop[f;T]=\inf\{\htop(\widehat{f}):\widehat{f}\in[f;T]\} \;. \]
If this infemum is not a minimum, we sometimes write $\htop[f;T]=\inf\{\htop(\widehat{f})\}+\epsilon$.
\end{definition}

We finally give two subclasses of trellis maps which are important for technical reasons.
In both these subclasses, we treat the stable and unstable sets differently; this anticipates the loss of time-reversal in the cutting procedure of Section~\ref{sec:cutting}.
We consider \emph{proper trellises} in order that the stable set be well-behaved after cutting.
\begin{definition}[Proper trellis]
We say a trellis $T=(T^U,T^S)$ is \emph{proper} if $\partial T^S\subset T^U$ and $\partial T^U\cap T^S=\emptyset$.
In other words, a trellis is proper if the endpoints of intervals in $T^S$ lie in $T^U$, but the endpoints of intervals in $T^U$ do not lie in $T^S$.
A trellis map $(f;T)$ is proper if $f^n(U)$ does not contain a point of $T^V$ for any open unstable end $U$ (so the ends of $T^U$ do not force any dynamics).
\end{definition}
Note that considering proper trellises is no great loss of generality; an unstable curve with an endpoint in $T^S$ can be extended by an arbitrarily short piece of curve, 
 and a stable curve with an endpoint not in $T^U$ can be shrunk by removing its end.
Further, since the ends of a trellis give no information about the dynamics, both these procedures are harmless and do not affect dynamical computations in any way.
Henceforth, all trellises will be taken to be proper trellises unless otherwise stated.

We consider \emph{well-formed} trellises in order that the isotopy classes of curves with endpoints in $T^S$ are well behaved under iteration.
This will be important when considering the graph representative defined in Section~\ref{sec:compatiblegraph}, as it gives a necessary condition for the topological entropy of the graph representative to be an optimal entropy bound for the trellis mapping class.
\begin{definition}[Well-formed trellis]
We say a proper trellis $T=(T^U,T^S)$ for a diffeomorphism $f$ is \emph{well-formed} if every component of $T^U\cup f(T^S)$ contains a point of $T^P$.
\end{definition}
Again, taking well-formed trellises is no great loss of generality, as we will usually take trellises for which the unstable branches are connected, and it remains to ensure that for every stable segment interval $S$, its image $f(S)$ intersects $T^U$, which we can do by removing the stable intervals whose iterates do not intersect $T^U$.


\subsection{Attractors}
\label{sec:attractor}

As well as saddle periodic orbits, generic diffeomorphisms also have attracting and repelling periodic orbits.
Clearly, a stable segment cannot be in the basin of attraction of an attracting orbit, and neither can an unstable segment be in the basin of a repelling orbit.
However, it is possible for the interior of an entire region to lie in the basin of an attracting or repelling periodic orbit.
Such regions are called \emph{stable} or \emph{unstable} regions.
\begin{definition}[Stable and unstable regions]
Let $([f];T)$ be a trellis mapping class. A domain $D$ of $T$ is \emph{stable} if there is a diffeomorphism $\widetilde{f}$ in $([f];T)$ and a periodic orbit $P$ of $f$ such that
 $\widetilde{f}^n(x)\tendsto P$ as $n\tendsto\infty$ for all $x$ in the interior of $D$.
If $D$ contains a point of $P$, then $D$ is \emph{attracting}.

Similarly, a domain $D$ is \emph{unstable} if there is a diffeomorphism $\widetilde{f}$ in $([f];T)$ and a periodic orbit $P$ of $f$
 such that $\widetilde{f}^{-n}(x)\tendsto P$ as $n\tendsto\infty$ for all $x$ in the interior of $D$, and if $D$ contains a point of $P$, then $D$ is \emph{repelling}.
\end{definition}

\fig{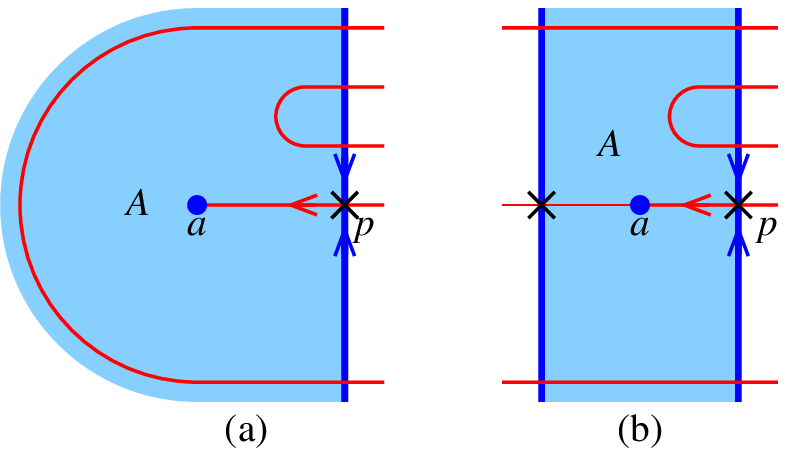}{Attracting domains.}
Two examples of attracting domains are shown in \figref{attractor}.
All points in the shaded domains are in the basin of attraction of the stable periodic point.
Note that $T^S$ is disjoint from the interior of $A$.
Repelling domains are similar.
Stable and unstable domains are wandering sets, unless they contain the attracting or repelling periodic orbit.

\fig{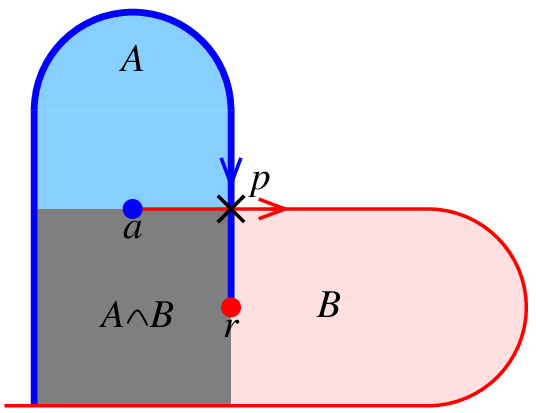}{An unbounded region.}
A region which contains an attractor and a repellor is called \emph{unbounded} region, by analogy with the unbounded region of the planar Smale horseshoe, for which the repellor is the point at infinity.
\figref{infinity} shows an unbounded region $R$ near a fixed point $p$. 
All points in $A$ lie in the basin of attraction of an attracting fixed point $a$,
 and all points in $B$ lie in the basin of a repelling fixed point $b$.
Since $A$ and $B$ cover the interior of $R$, all points in the interior of $R$ are wandering, apart from $a$ and $b$.

An open stable segment which is in the basin of a repelling periodic orbit is also a wandering set,
 as is an open unstable segment which is in the basic of an attracting periodic orbit.
There are some nonwandering segments such that any finite iterate need not contain an intersection point.
We call such segments \emph{almost wandering}.
\begin{definition}[Almost wandering segment]
\label{defn:almostwandering}
Let $([f];T)$ be a trellis mapping class.
An open segment is \emph{almost wandering} if it is nonwandering, but for any integer $n$ there exists a diffeomorphism $\tilde{f}\in([f];T)$ the $n\th$ iterate of the segment by $\tilde{f}$ contains no intersection points.
\end{definition}

Therefore, an open stable segment $S$ is almost wandering if for any positive integer $n$, there exists $\widetilde{f}\in([f];T)$ such that $\widetilde{f}^{-n}(S)\cap T^U=\emptyset$, but for any diffeomorphism $\widehat{f}\in([f];T)$, there exists $n$ such that $\widehat{f}^{-n}(S)\cap T^U\neq\emptyset$.
A similar statement holds for unstable segments.

\begin{definition}[Chaotic region]
Let $([f];T)$ be a trellis mapping class.
A region $R$ is \emph{chaotic} if for every diffeomorphism $\widehat{f}$ in $([f];T)$, there exists an integer $n$ such that $f^n(T^U)\cap f^{-n}(T^S)$ contains a point in the interior of $R$.
\end{definition}
A chaotic region supports entropy, in the sense that for any diffeomorphism $\widehat{f}\in([f];T)$, there is an ergodic measure $\mu$ with positive entropy such that $\mu(R)>0$.


\subsection{Extensions}
\label{sec:extension}

A surface diffeomorphism $f$ with a periodic saddle orbit has infinitely many trellises, which are partially ordered by inclusion.
Taking a smaller trellis gives a \emph{subtrellis}, and a larger trellis a \emph{supertrellis}.
\begin{definition}[Subtrellis and supertrellis]
Let $(f;T)$ and $(f;\widehat{T})$ be trellis maps.
If $T^\us\subset\widehat{T}^\us$, then $T$ is a \emph{subtrellis} of $\widehat{T}$, and $\widehat{T}$ is a supertrellis of $T$.
Similarly, if $([f];T)$ is a trellis mapping class, a trellis map $(\widehat{f};\widehat{T})$ is a supertrellis for $([f];T)$ if $\widehat{f}\in([f];T)$ and $T^\us\subset \widehat{T}^\us$.
\end{definition}
Since taking a subtrellis loses information about the trellis map, we are more interested in supertrellises.
Of particular importance are those trellises which can be obtained by iterating segments or branches.
\begin{definition}[Iterate]
Let $(f;T)$ be a trellis map.
Then a trellis $\widehat{T}$ is an \emph{$f$-iterate} of $f$ if there exist positive integers $n_u$ and $n_s$ such that $\widehat{T}^U=f^{n_u}(T^U)$ and $\widehat{T}^S=f^{-n_u}(T^S)$.
An iterate is a \emph{stable iterate} if $\widehat{T}^U=T^U$ and an \emph{unstable iterate} if $\widehat{T}^S=T^S$.
A trellis map $(\widehat{f};\widehat{T})$ is an iterate of the trellis mapping class $([f];T)$ if $\widehat{f}\in([f];T)$ and $\widehat{T}$ is an $\widehat{f}$-iterate of $T$.
\end{definition}
The concept of an iterate is rather restrictive in practice, as we may wish to iterate only a few segments of $T$.
Further, trellises are often computed by growing branches out from $T^P$, so we may obtain a supertrellis for which $\widehat{T}^P=T^P$ but which is not an iterate.
A class of supertrellis which is more restrictive than that of an iterate is an \emph{extension}
\begin{definition}[Extension]
Let $(f;T)$ be a trellis map. 
A trellis $\widehat{T}$ is an \emph{$f$-extension} of $T$ if there exists $n$ such that
  \[ T^U\subset \widehat{T}^U\subset f^n(T^U) \quad \textrm{and} \quad T^S\subset \widehat{T}^S\subset f^{-n}(T^S) . \]
A trellis map $(\widehat{f};\widehat{T})$ is an \emph{extension} of a trellis mapping class $([f];T)$ if
 $\widehat{f}\in([f];T)$ and $\widehat{T}$ is a $\widehat{f}$-extension of $T$.
\end{definition}
In other words, a trellis $\widehat{T}$ is an extension of $T$ if it is a subtrellis of some iterate of $T$.
Just as for iterates, a \emph{stable extension} has $\widehat{T}^U=T^U$ and an \emph{unstable extension} has $\widehat{T}^S=T^S$.

If $S$ is a segment of $T^S$ and $f^{-1}(S)$ is not a subset of $T^S$, we can take a stable $f$-extension $\widehat{T}$ with $\widehat{T}^S=T^S\cup f^{-1}(S)$.
We say this extension is formed by taking a backward iterate of $S$.
In a similar way, we can form extensions by taking backward iterates of arcs of $T^S$, or forward iterates of arcs of $T^U$.
For these extensions, the endpoints of any curve of $\widehat{T}$ lie in the set $X=\bigcup_{n=-\infty}^{\infty}f^n(T^V)$.

By taking infinitely many iterates of $T^U$ and $T^S$, we obtain a tangle known as an \emph{infinite extension}.
\begin{definition}[Infinite extension]
Let $(f;T)$ be a trellis map.
Then the \emph{infinite $f$-extension} of $T$ is the tangle $W(f;T^P)$.
\end{definition}
As long as no point of $T^P$ is the end of a curve of $T$, the tangle $W(f;T^P)$ can be obtained by iterating all branches of $T$ infinitely often under $f$.

The main difference between extensions and supertrellises is that $\widehat{T}^P=T^P$ for an extension, but $T^P$ may be a strict subset of $\widehat{T}^P$ for a supertrellis.
This difference makes the analysis of supertrellises slightly more complicated than that of extensions.
In particular, a curve in $\widehat{T}^U$ may contain more than one periodic point, in which case all periodic points lie in the same region of $T$, and we can isotope so the curve contains an attractor between any two periodic points.
This property makes supertrellises particularly useful, since we can first construct some invariant curves, and then decide how many points of $\widehat{T}^P$ they should contain.

\fig{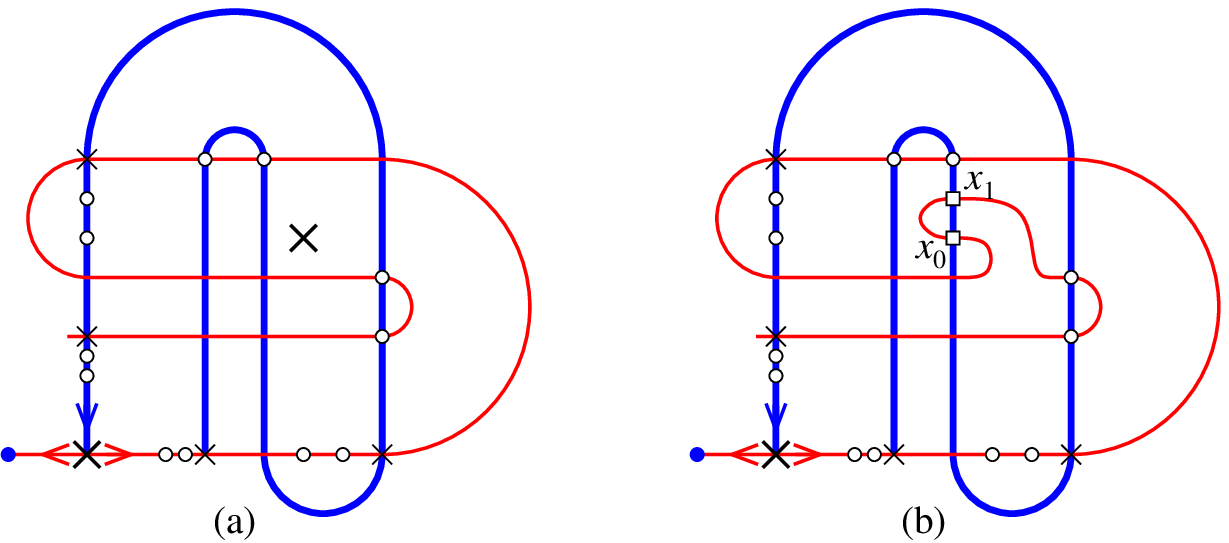}{Extensions of the Smale horseshoe trellis (a) is a minimal extension, (b) is not.}
Of particular interest are \emph{minimal} iterates, extensions and supertrellises, which introduce as few extra intersections as possible given the endpoints of the branches.
The trellis shown in \figrefpart{smaleexts}{a} is a minimal extension of the Smale horseshoe trellis, but the trellis shown in
\figrefpart{smaleexts}{b} is not, since the intersection points $x_0$ and $x_1$ can be removed.


\subsection{Reducibility}
\label{sec:reducibility}

The concept of reducibility for trellis mapping classes is similar to that for mapping classes of surface diffeomorphisms,
 but is complicated by the fact that it may not be possible to decompose a system into a number of invariant subsets
 of negative Euler characteristic, but may be possible to decompose the system into an attractor-repeller pair.
This gives two slightly different notions of a reductions, namely an \emph{invariant curve reduction} and an \emph{attractor-repeller reduction}
\begin{definition}[Invariant curve reduction]
A trellis mapping class $([f];T)$ has an \emph{invariant curve reduction} if there is a diffeomorphism $\widetilde{f}\in([f];T)$ and a one-manifold $C$ which is disjoint from $T$,
 invariant under $\widetilde{f}$, and each component of the complement of $C$ either contains a subtrellis of $T$, or has negative Euler characteristic.
If the complement of $C$ is disconnected, we have a \emph{separating reduction}, otherwise we have a \emph{non-separating reduction}.
\end{definition}
\begin{definition}[Attractor-repeller reduction]
A trellis mapping class $([f];T)$ has an \emph{attractor-repeller reduction} if there is a diffeomorphism $\widetilde{f}\in([f];T)$ and a one-manifold $C$
 such that $C$ divides $M$ into subsets $A$ and $B$ such that $\cl{\widetilde{f}(A)}\subset\inr{A}$, and each component of the complement of $C$ is either contains a point of $T^P$,
 or has negative Euler characteristic.
\end{definition}
In both cases, the set $C$ is called a set of \emph{reducing curves}.
It is clear that reducibility is a property of the trellis type, and not just the trellis mapping class.

The essence of both these definitions is that the manifold simplifies into two subsets which we can treat independently.
However, the invariant curve reduction allows for curves do not divide the manifold, but instead decrease its genus.
In both cases, if we find a reduction, we cut along the reducing curves,
 to obtain a simpler problem, or a number of simpler problems, which we deal with independently. 
After cutting along the reducing curves of an invariant curve reduction, we immediately obtain a diffeomorphism.
After cutting along the reducing curves of an attractor-repellor curve reduction, we need to glue on an annulus to each curve.
  After glueing, we extend the diffeomorphism $\widetilde{f}$ to the annulus in such a way that the remaining boundary is an attractor/repeller.
For most of this paper, we only consider irreducible trellis mapping classes.
We consider the problem of finding reductions in Section~\ref{sec:algorithm}.

\fig{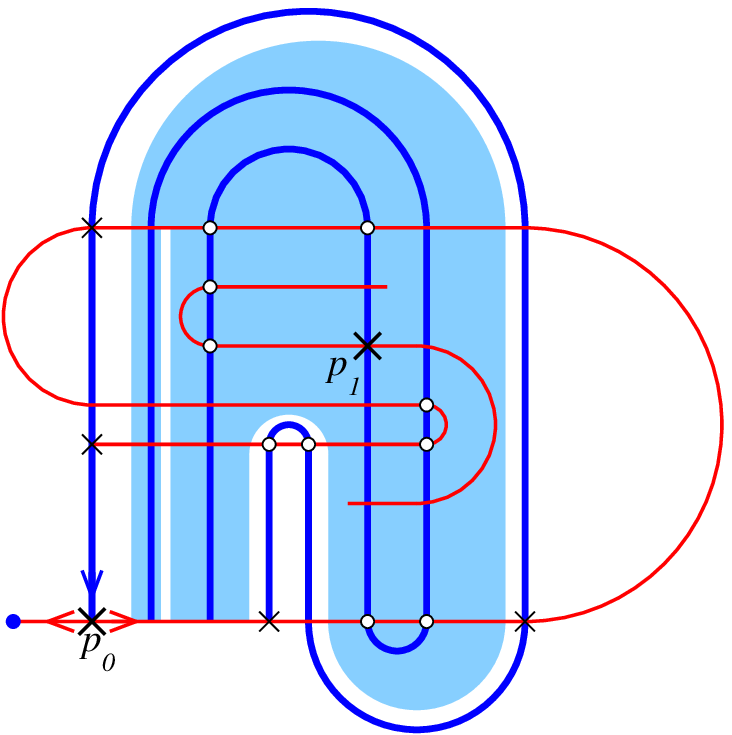}{An trellis with an attractor-repellor reduction.}
In \figref{arreduction} we show a trellis type with an attractor-repellor reduction.
A sub-basin $A$ of an attractor is shaded.
Notice that the unstable set $T^U(p_1)$ lies in $A$, which means it is possible to find a diffeomorphism $\widetilde{f}$ for which $W^U(\widetilde{f};p_1)\cap W^S(\widetilde{f};p_0)=\emptyset$.
The closure of the unstable manifold of $p_1$ is then a non-trivial chaotic attractor, though it may have non-trivial sub-attractors.

We now examine the hypothesis of irreducibility in more detail.
Recall that the main effect of reducibility is that the dynamics of a reducible trellis mapping class can be separated into simpler pieces.
We now show that the converse is also true; the dynamics forced by an irreducible trellis essentially lies in one minimal component.
To prove this, we need to carefully consider the intersections of branches of the trellis.
We write $W^U(p_u,b_u)\trans W^S(p_s,b_s)$ if the branches $W^U(p_u,b_u)$ and $W^S(p_s,b_s)$ intersect transversely.

For the case of trellises, we consider a relationship based on \emph{heteroclinic chains}.
\begin{definition}[Heteroclinic chains]
Let $([f];T)$ be an transverse trellis mapping class, and $p_u,p_s\in T^P$.
We say there is a \emph{heteroclinic chain} from $p_u$ to $p_s$ of length $n$ if there are points $p_u=p_0,p_1,\ldots p_n=p_s$ in $T^P$
 and an integer $m$ such that for every $\widehat{f}\in([f];T)$ we have $T^U(p_i) \cap \widehat{f}^{-m}(T^S)(p_{i+1})\neq\emptyset$ for $i=0,1,\ldots n-1$.
We write $p_u\succ p_s$ if there is a heteroclinic chain from $p_u$ to $p_s$.
A subset $P$ of $T^P$ is \emph{chain transitive} if $p_1\succ p_2$ for all $p_1,p_2\in P$.
\end{definition}
Clearly, $\succ$ is a transitive relation, though it need not be reflexive or symmetric.

Note that although $p_u\succ p_s$ implies $W^U(f;p_u)$ intersects $W^S(f;p_s)$ for any $f$ in the trellis mapping class, it is not necessarily true
 that there exists an $n$ such that $T^U(p_u)$ intersects $f^{-n}(T^S(p_s))$ for any $f$ in the trellis mapping class.
In other words, the intersections may lie in arbitrarily high iterates.
For the next results, we therefore consider tangles rather than trellises.

\begin{lemma}
Let $([f];T)$ be a well-formed trellis mapping class, and $p_u=p_0,p_1,\ldots p_n=p_s$ be a heteroclinic chain.
Let $W$ be the infinite $f$-extension of $T$.
Then $W^U(p_u)\trans W^S(p_s)$.
\end{lemma}
\begin{proof}
By the Lambda lemma, if $W^U(p_{i-1},b_{i-1})\trans W^S(p_i)$ and $W^U(p_i)\trans W^S(p_{i+1},b_{i+1})$,
 then $W^U(p_{i-1},b_{i-1})$ transversely intersects $W^S(p_{i+1},b_{i+1})$.
Induction on the length of the heteroclinic chain completes the proof.
\end{proof}

We now use the above observation to show that the chain transitive components of $T$ are permuted.
\begin{lemma}
Let $([f];T)$ be a trellis mapping class.
Suppose $P$ is a chain transitive subset of $T^P$.
Then there exists least $N\in\N$ such that $f^N(P)=P$.
Further, if $0<i<N$, then $p\not\succ f^i(p)\not\succ p$ for any $i$.
\end{lemma}
\begin{proof}
Suppose $p_u\succ p_s$.
Let $\widetilde{f}$ be a diffeomorphism such that $\widetilde{f}^{n}(T^U(p_u))\trans T^S(p_s)$ and $(\widetilde{f}^{n+1}(T^U),T^S)$ is a minimal extension of $T$.
Then $\widetilde{f}^{n+1}(T^U)(p_u)\trans \widetilde{f}(T^S(p_s))$, so $\widetilde{f}^{n}(T^U(f(p_u)))\trans T^S(f(p_s))$, so $f(p_u)\succ f(p_s)$.
Similarly, we can see that $f^{-1})(p_u)\succ f^{-1}(p_s)$, and hence  $f^n(p_u)\succ f^n(p_s)$ for any $n\in\Z$.

Since the points of $P$ are periodic and reflexive under $\succ$, there exists least $N\in\N$ such that $p\succ f^{N}(p)$ for some $p\in P$.
Then $f^{iN}(p)\succ f^{(i+1)N}(p)$, for any $i$, and since $p$ is periodic, we deduce $p\succ f^{iN}(p)\succ p$ for any $i$.

Now suppose $q\in P$, so $p\succ q\succ p$.
Then for any $i$, we have $p\succ f^{iN}(p)\succ f^{ni}(q)\succ f^{iN}(p)\succ p$, so $f^{iN}(q)\in P$.
Conversely, we can show that if $f^{iN}(q)\in P$, then $q\in P$.
Therefore, $P=f^{N}(P)$.
\end{proof}
Note that the construction implies that $N$ must divide the period of any point of $P$.

The relation $\succ$ induces a relation on the chain transitive components which is antisymmetric and transitive, but is not necessarily reflexive.
If $P$ is a reflexive component, then $P$ contains a homoclinic orbit, so induces some chaotic behaviour.
A chain transitive component $P$ which is not reflexive must consist of a single periodic point $p$ with $W^U(p)\cap W^S(p)=\emptyset$, so there are no homoclinic orbits to $p$.

If $P_R$ and $P_A$ are chain transitive components and $P_R\succ P_A\not\succ P_R$, then we hope to find an attractor-repellor reduction with $P_A$ in the attractor and $P_R$ in the repellor.
If $P=f(P)$, then the tangle of the orbits of $P$ is connected, but if $P=f^N(P)$ for some least $N>1$, then the tangle of the orbits of $P$ may have several components,
 which may form the basis for some invariant curve reduction.
If $P$ is not reflexive, then $W^U(P)$ may not intersect $W^S(P)$, and $P$ is a non-chaotic saddle.
The following theorem makes this interpretation precise.

\begin{theorem}[Intersections of irreducible tangles]
\label{thm:branchintersection}
Let $([f];T)$ be a well-formed irreducible transverse trellis mapping class, and let $W$ be the infinite $f$-extension of $T$.
Then if $T^U(p_u,b_u)$ and $T^S(p_s,b_s)$ are closed branches of $T$, the extensions $W^U(p_u,b_u)$ and $W^S(p_s,b_s)$ intersect transversely.
\end{theorem}

\begin{proof}
The relationship $\succ$ clearly defines a partial order on the chain transitive sets.
Choose a set $P$ which is minimal in this partial order, and let $P_A=\bigcup_{n=1}^\infty f^n(P)$.

Let $([f_0];T_0)=([f];T)$, and for $i=1,\ldots\infty$, let $([f_i];T_i)$ be a minimal extension of $([f_{i-1}];T_{i-1})$ with $T_i=(f_i(T_{i-1}^U),T_{i-1}^S)$.
Let $K_i$ be a set consisting of $T_i^U(P_A)$, and all regions of $T_i$ which are topological discs or annuli and all of whose stable boundary segments lie in $T^S(P_A)$.
Since $T_i^U(P_A)\cap T_i^S(P_R)=\emptyset$, $K_i\subset K_{i+1}$.
Let $A_i$ be a closed neighbourhood of $K_i$ which deformation retracts onto $K_i$, and such that $A_i\subset A_{i+1}$.
Now, $\chi(K_i\setminus T^U)\geq\chi(K_{i+1}\setminus T^U)\geq\chi(M)$, so the Euler characteristics are a decreasing bounded sequence of integers, so have a limit,
 which occurs for some $A_n$.
Then $A_{n+1}$ deformation retracts onto $K_n$.
Further, $f(K_n)\subset K_{n+1}\subset A_{n+1}$, so by taking $A_n$ sufficiently small, we can ensure $f(A_n)\subset A_{n+1}$.
Choose an isotopy $h_t$ such that $h_1(A_{n+1})\subset A_n$ and $h_t$ is fixed on $K$, and let $\widetilde{f}_t=h_t\circ f_n$.
Then $\widetilde{f}_1(A_n)=h_1(f_n(A_n))\subset h_1(A_{n+1})\subset A_n$, so $\widetilde{f}_1\in([f];T)$ and maps $A_n$ into itself.

If $P_A\neq P$, then this gives an attractor-repellor reduction, a contradiction.
If $P_A=P$ and $\chi(A_n)\neq\chi(M)$, we can isotope $f$ so that the boundary of $A_N$ maps to itself, giving an invariant curve reduction, also a contradiction.
Therefore $P_A=P$ and there is no reduction.
Since $T^U(p_u,b_u)$ and $T^S(p_s,b_s)$ are closed branches, there are periodic points $\widetilde{p}_u$ and $\widetilde{p}_s$ such that $T^U(p_u,b_u)\trans T^S(\widetilde{p}_u)$ and $T^U(\widetilde{p}_s)\trans T^S(p_s,b_s)$.
Then $T^U(\widetilde{p}_u)\trans T^S(\widetilde{p}_s)$, and so $T^U(p_u,b_u)\trans T^S(p_s,b_s)$ as required.
\end{proof}

In other words, if $([f];T)$ is irreducible, then for there exists $n$ such that the $n$th image of every unstable branch under $f$ intersects every stable branch.
This is, essentially, the converse to the trivial observation that if $([f];T)$ has an attractor-repellor reduction, then there is a diffeomorphism $\widetilde{f}$ in $([f];T)$ and branches of $W(\widetilde{f};T^P)$ which do not intersect.


\subsection{Cutting}
\label{sec:cutting}

Given a trellis, the first step in obtaining rigorous information on the dynamics is to \emph{cut} along the unstable manifold to obtain a new space $\cut{T}$, a topological pair which we denote $(\cut[U]{T},T^S)$.
A formal description of cutting is given in \cite{Collins99AMS}, so we shall not give one here as it is intuitively obvious.

\fig{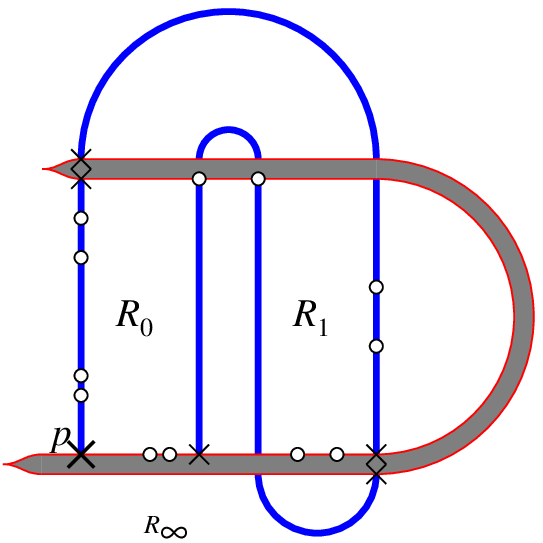}{Cutting along the unstable curve.}
Cutting along the unstable curve of the Smale horseshoe trellis gives the set shown in \figref{cutsmale}.
The ends of the unstable curves gives \emph{cusps} in the new surface.
If $T$ is a proper trellis, the stable segments lift to cross-cuts in the new surface $\cut[U]{T}$.
End intersections of $T^S$ also give isolated points in $\cut[U]{T}$ which are discarded.
The topology of the lift of $T^S$ to $\cut[U]{T}$ is sufficiently similar to that of $T^S$ that we also denote this lift by $T^S$.

If $f$ is a trellis map, $f(T^U)$ covers $T^U$, so we have a well-defined function $\cut{f}$ on $\cut[U]{T}$.
This function is not a diffeomorphism since a point $q$ of $T^U$ which maps to a point $f(q)\not\in T^U$
 will lift to two points with the same image.
However, $\cut{f}$ is surjective and at most two-to-one; in particular, $f$ has the same topological entropy.
Since $\cut{f}$ maps $T^S$ into itself, so is a map of the topological pair.

Notice that the pair $\cut{T}=(\cut[U]{T},T^S)$ contains the pair $(M\setminus T^U,T^S\setminus T^U)$ as an open subset which is invariant under $\cut{f}$.
Indeed, $\cut[U]{T}$ can be regarded as a natural compactification of $M\setminus T^U$.
This compactness is important in the rigorous application of the Nielsen theory in Section~\ref{sec:shadowing}.
However, the homotopy properties of $(\cut[U]{T},T^S)$ and $(M\setminus T^U,T^S\setminus T^U)$ are essentially the same.


\subsection{Curves}
\label{sec:curves}

Our main tool for studying the geometry, topology and dynamics associated with trellis maps will be to consider curves embedded in the cut surface $(\cut[U]{T},T^S)$.
As this surface is a topological pair, our curves will be maps in this category, so a curve $\alpha$ in $T$ is a mapping $\alpha:(I,J)\fto(\cut[U]{T},T^S)$,
 where $I$ is the unit interval $[0,1]$, and $J$ is a closed subset of $I$.
The \emph{path} of such a curve $\alpha$ is the set $\alpha(I)$.

Since the topological pair $\cut{T}$ is obvious from the trellis, we will usually draw curves in $\cut[U]{T}$ as curves embedded in the original surface $M$, and,
 wherever possible, ensure these curves are disjoint from $T^U$.
However, as long as a curve in $M$ does not cross $T^U$, it lifts to a curve in $\cut[U]{T}$.
This simplifies many of the diagrams, and hence clarifies the geometry of the situation.

For the most part, we are only interested in curves up to homotopy or isotopy, and we always take homotopies and isotopies of curves through maps of pairs.
For the Nielsen theory and shadowing results of Section~\ref{sec:shadowing}, we will always keep the endpoints fixed during the homotopy.
For most other purposes, we only consider curves in $\cut[U]{T}$ whose endpoints lie in $T^S$ and for which the set $J$ consists of finitely many points including $0$ and $1$.
In this case, we only consider homotopies of curves which restrict to isotopies on $J$, or in other words, we ensure that $\alpha_t$ is one-to-one on $J$.
However, we do allow the endpoints to move, though they are of course restricted to a segment of $T^S$.

If $\alpha:(I,J)\fto(X,Y)$ is a curve and $J$ contains $\{0,1\}$ we say $\alpha$ has endpoints in $Y$; if $J$ equals $\{0,1\}$ we say $\alpha$ has endpoints only in $Y$.

Reparameterising a curve does not change its path, but may change the set $J$ which maps into $T^S$.
This means that different parameterisations of the same path may not even be comparable under homotopy.
However, we consider different parameterisations of the same curve as equivalent.
\begin{definition}[Equivalence of curves]
Curves $\alpha_1:(I_1,J_1)\fto(X,Y)$ and $\alpha_2:(I_2,J_2)\fto(X,Y)$ are \emph{equivalent} if there is a homeomorphism $h:(I_1,J_1)\fto(I_2,J_2)$ with $h(J_1)=J_2$
 such that $\alpha_1\homotopic\alpha_2\circ h$ as curves $(I_1,J_1)\fto(X,Y)$.
(The homotopy may be taken relative to endpoints, as appropriate.)
\end{definition}

To study a trellis mapping class $([f];T)$, we consider the iterates of an exact curve $\alpha:(I,J)\exto(\cut[U]{T},T^S)$ with endpoints in $T^S$ under models of $([f];T)$.
We are especially interested in iterates with minimal number of intersections with $T^S$.
\begin{definition}[Minimal iterate]
Let $([f];T)$ be a trellis mapping class, and $\alpha:(I,J)\exto\cut{T}$ be a simple curve with endpoints in $T^S$ and transverse intersections with $T^S$.
Then a \emph{minimal iterate} of $\alpha$ is a curve $\beta$ which is homotopic to $\cut{f}\circ\alpha$ relative to $J$ and which minimises the number of intersections with $T^S$.
We let $J\p=\beta^{-1}(T^S)$, and consider $\beta$ as a curve $(I,J\p)\fto\cut{T}$.
\end{definition}

\fig{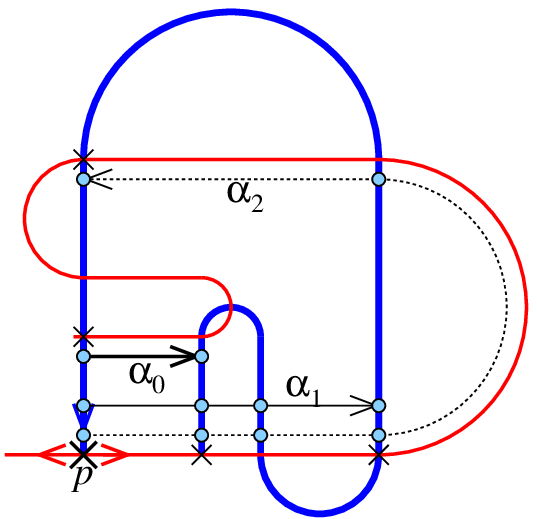}{Curve $\alpha_0$ and minimal iterates $\alpha_1\in f_{\min}[\alpha_0]$ and $\alpha_2\in f_{\min}[\alpha_1]$.}

Theorem~\ref{thm:minimalcurve} shows that the curve $\beta$ is only well-definite up to equivalence.
We therefore obtain a well-defined map $f_{\min}$ on equivalence classes of curves given by $f_{\min}[\alpha]=[\beta]$.
An example of the first two minimal iterates of a curve $\alpha_0$ is given in \figref{iteratecurves}.
Note that by $f_{\min}^n[\alpha]$ we mean ${(f_{\min})}^n[\alpha]$, and not ${(f^n)}_{\min}[\alpha]$, which may be different.
By $f_{min}(\alpha)$, we mean some curve in the homotopy class $f_{\min}[\alpha]$.
Minimal iterates are closely related to minimal extensions of trellises.

\fig{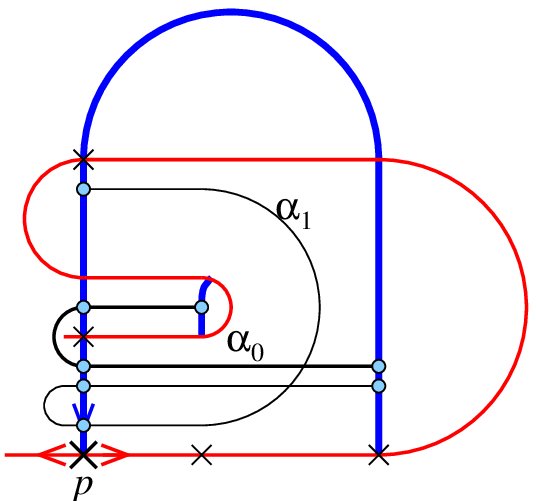}{Curve $\alpha_0$ and minimal iterate $\alpha_1$ which self-intersects.}
In \figref{badminimaliterate} we show a trellis type which is not well-formed, and a curve $\alpha_0:(I,J)\exto(\cut[U]{T},T^S)$ with endpoints in $T^S$.
The minimal iterate of $\alpha_0$ has self-intersections; note that by our notion of isotopy for such curves, points in $\alpha_t(J)$ cannot cross each other.
However, this is impossible for a well-formed trellis mapping class.

\begin{lemma}
\label{lem:wellformed}
Suppose $([f];T)$ is a well-formed transverse trellis mapping class, and $\alpha:(I,J)\exto(\cut[U]{T},T^S)$ is a simple exact curve in $\cut[U]{T}$ with endpoints in $T^S$ and transverse intersections.
Then there is a minimal iterate $f_{\min}[\alpha]$ which is a simple curve.
\end{lemma}

\begin{proof}
Let $K$ be any disc in $\cut[U]{T}$ such that the boundary of $K$ consists of an arc of $\beta$ with endpoints $x_0$ and $x_1$, and a subinterval of a segment of $T^S$.
If neither $x_0$ or $x_1$ lie in $f(S)$, then there is a curve homotopic to $\beta$ with fewer intersections with $T^S$ than $\beta$ obtained by removing the intersections $x_0$ and $x_1$.

Suppose there is a point $q\in\alpha(J)$ such that $f(q)\in T^S(x_0,x_1)$, as shown in \figref{iterateisotope}, and let $S$ be the component of $T^S$ containing $q$.
\fig{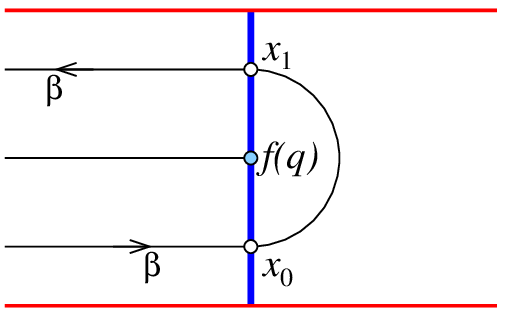}{If neither $x_0$ or $x_1$ lie in $f(T^S)$, then $T$ is not well-formed.}
Then $q$ lies in a component $S$ of $T^S$.
Then \[f^{-1}(T^U)\cap S=f^{-1}(T^U\cap f(S)\subset f^{-1}(T^U\cap T^S(x_0,x_1))=\emptyset,\] which contradicts the fact that $T$ is well-formed.
Therefore, if there are any points of $\beta$ in $T^S(x_0,x_1)$, these also do not lie in $\beta(J)$, so can be removed without introducing any self-intersections.
\end{proof}
This result is a \emph{homotopy} result, unlike the isotopy results of Section~\ref{sec:entropy}.

Sometimes we want to relate a curve $\widehat{\alpha}:(\widehat{I},\widehat{J})\fto(\cut[U]{T},T^S)$ to a curve $\alpha:(I,J)\fto(\cut[U]{T},T^S)$.
If $(\widehat{I},\widehat{J})\neq(I,J)$, the curves cannot be homotopic, but if $(I,J)\subset(\widehat{I},\widehat{J})$ then the curve $\widehat{\alpha}$ may be in some sense more complicated than $\alpha$, as defined below.
\begin{definition}[Tightening curves]
\label{defn:curvetightening}
Let $\widehat{\alpha}:(\widehat{I},\widehat{J})\exto(X,Y)$ and $\alpha:(I,J)\exto(X,Y)$ be exact curves such that $(I,J)\subset(\widehat{I},\widehat{J})$.
We say $\widehat{\alpha}$ \emph{tightens to} $\alpha$ if the curve $\alpha_0:(I,J)\fto(X,Y)$ given by $\alpha_0(t)=\widehat{\alpha}(t)$ is homotopic to $\alpha_1=\alpha$.
\end{definition}
It is clear that the tightening relation is a partial order on homotopy classes of curves.



\section{Biasymptotic Orbits and Minimal Trellises}
\label{sec:biasymptoticminimal}

We use the term \emph{biasymptotic orbit} to refer to an orbit which is either homoclinic or heteroclinic to the set of periodic points.
In this section, we consider how to extend a set of biasymptotic orbits to a trellis, and how to extend a trellis to a larger trellis.
We shall usually consider extensions with a given \emph{Birkhoff signature}.
As previously remarked, this does not guarantee a unique trellis type.
Instead, we define what it means for a trellis to be \emph{minimal} given a biasymptotic mapping class.
By carefully isotoping a given trellis map, we can construct a minimal trellis, and we then show that the type of a minimal trellis is almost unique given its end points, the only ambiguities arising from the orientation (transverse or tangent) of the intersections.

We also consider minimal iterates and minimal extensions of a given trellis mapping class.
A minimal extension is forced by the biasymptotic orbits of the intersections of the original trellis.
Unlike a minimal trellis for a biasymptotic mapping class, the orientations of the intersections are uniquely determined.
We also define minimal supertrellises, for which we need to consider changes to to periodic point set.


\subsection{Biasymptotic orbits}
\label{sec:biasymptotic}


Just as we can consider isotopy classes of surface diffeomorphisms relative to periodic orbits and trellises,
 we can also consider isotopy classes relative to homoclinic and heteroclinic orbits to periodic saddle orbits. 

\begin{definition}[Biasymptotic mapping class]
Let $f$ be a surface diffeomorphism, and $X$ a closed, invariant set consisting of a set of periodic saddle orbits $X^P$ and a set $X^V$ of biasymptotic orbits to $X^P$.
The isotopy class of $f$ relative to $X$ is called the \emph{biasymptotic mapping class} $([f];X)$.
The conjugacy class of a biasymptotic mapping class is a \emph{biasymptotic type}.
\end{definition}

Methods for computing the dynamics forced by a biasymptotic mapping class have been given by  Handel \cite{Handel99TOPOL} and Hulme \cite{Hulme00PHD}.
Here we show how to relate a biasymptotic type to a trellis type class, and hence derive a new method for computing the dynamics.

Since $X^P$ is a finite set of saddle orbits, we can compute the stable and unstable manifolds $W^\us(f;X^P)$.
Since the points of $X^V$ are biasymptotic to $X^P$, they are all intersection points of $W^U$ and $W^S$.
However, $W^U$ and $W^S$ may have extra intersections which are not in $X$, and indeed, typically have infinitely many such intersections.
Different representatives of a biasymptotic mapping class will have different tangles.

\begin{definition}[Compatible tangle/trellis]
Let $([f];X)$ be a biasymptotic mapping class. 
Then if $\widehat{f}$ is any diffeomorphism in $([f];X)$, the tangle $\widehat{W}=W(\widehat{f};X^P)$ is a \emph{compatible tangle} for $([f];X)$.
A trellis map $(\widehat{f};\widehat{T})$ is \emph{compatible} with $([f];X)$
 if $\widehat{f}\in([f];X)$ and $\widehat{T}$ is a subtrellis of the compatible tangle $\widehat{W}$ with $X^P=T^P(=W^P)$.
\end{definition}

There are many different trellis types compatible with a biasymptotic type $([f];X)$.
To restrict the possible trellises, we consider only trellises with connected branches, and specify the endpoints of each branch as a point of $X$.
This gives rise to the notion of a \emph{Birkhoff signature}.
\begin{definition}[Birkhoff signature]
A \emph{Birkhoff signature} for a biasymptotic mapping class $([f];X)$ is a pair $\sig=(\sig^U,\sig^S)$ of subsets of $X$
 such that $\sig^U$ contains one point of $X$ from each branch of $W^U(f;X^P)$ and $\sig^S$ contain one point of $X$ from each branch of $W^S(f;X^P)$.
If $\sig_1$ and $\sig_2$ are two signatures, we say $\sig_1\leq \sig_2$ if the trellis $T_1$ for $f$ compatible with $(f;X)$ and with endpoints in $\sig_1$ is a subtrellis of the compatible trellis $T_2$  with endpoints in $\sig_2$.
\end{definition}
We can further restrict the compatible trellises by specifying whether the intersection of $W^U$ and $W^S$ is transverse $\trans$ or tangential $\tang$ at each point of $X$,
Note that all points on the same orbit must have the same intersection type.

\fig{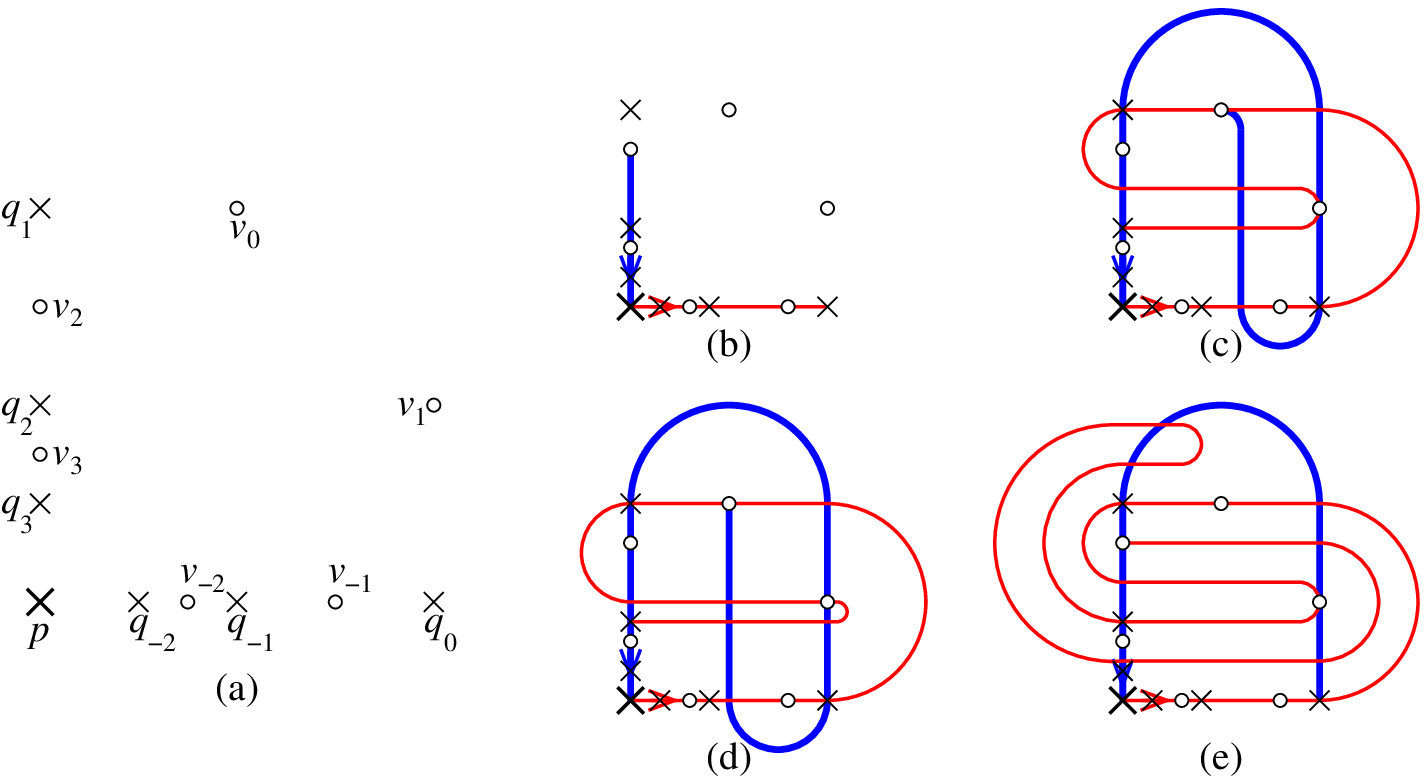}{(a) Orbits $\{q_i\}$ and $\{v_i\}$ of the Smale horseshoe map. (b-e) Trellises compatible with $\{q_i\}$ and $\{v_i\}$.}
\figref{horseshoebraid} shows some trellises compatible with the biasymptotic mapping class $([f];X)$, where \[X=\{p\}\cup\{q_i\}\cup\{v_i\}\]
The orbits $q_i$ and $v_i$ shown in \figrefpart{horseshoebraid}{a} are orbits of the Smale horseshoe map homoclinic to the saddle point $p$.
Both these orbits lie in the same branch of $W^U(p)$ and $W^S(p)$, so a signature consists of a pair of points from either orbit.
The trellis in \figrefpart{horseshoebraid}{b} has signature $\left(\{q_0\},\{v_2\}\right)$ and no intersections other than $p$.
The trellises in \figrefpart{horseshoebraid}{c} and \figrefpart{horseshoebraid}{d} have signature $\left(\{q_2\},\{v_0\}\right)$.
However, the trellis in \figrefpart{horseshoebraid}{c} has a tangency at $v_1$, whereas the trellis in \figrefpart{horseshoebraid}{d} has a transverse intersection
 with negative orientation.
The trellis in \figrefpart{horseshoebraid}{e} has signature $\left(\{v_2\},\{q_0\}\right)$, but has a pair of transverse intersections which is not forced by $X$.

\fig{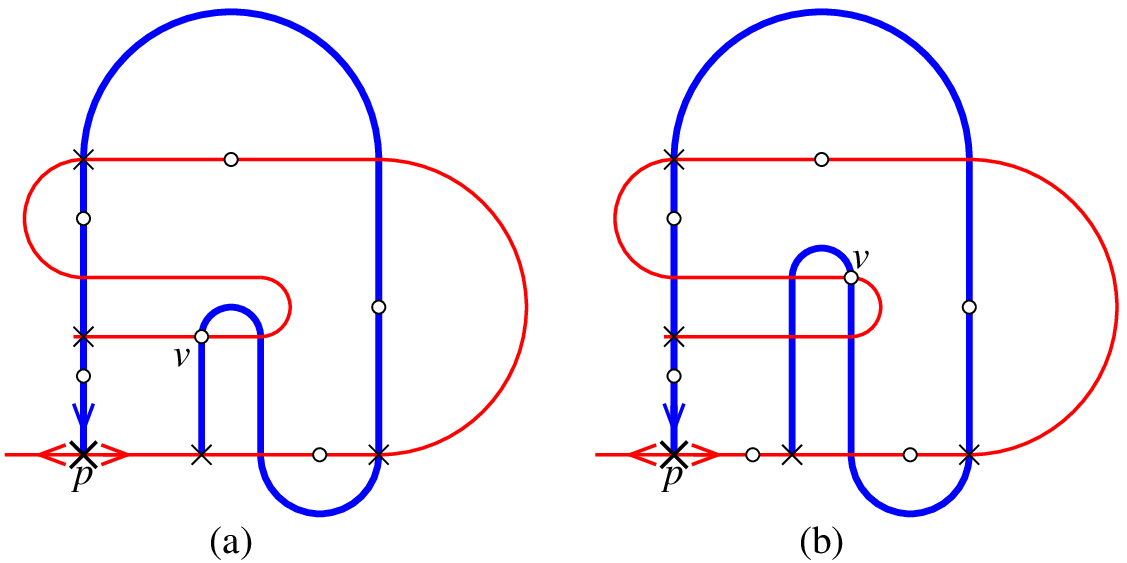}{The trellis mapping types shown in (a) and (b) are not isotopic through maps preserving the intersection at $v$.}
If we specify the orientation of a crossing, it is not necessarily true that every compatible trellis map lies in the same isotopy class relative to the oriented intersections.
An example is shown in \figref{orientationspecification}.
Both trellis types shown have a set $X$ consisting of the fixed saddle point $p$ and the orbit of a homoclinic point $v$, and both yield the same biasymptotic mapping class $(f;X)$.
Further, the orientation of the crossing at $v$ is the same in each case.
However, there is no isotopy between the maps in the different isotopy classes relative to $X$ such that the orientation of the intersection at $v$ does not change.

Clearly, even given a Birkhoff signature, there are still many, indeed infinitely many, trellis types compatible with a biasymptotic type.
However, some trellises have extra intersections which can be removed by an isotopy in the biasymptotic mapping class.
A trellis with no extra intersections is called a \emph{minimal compatible trellis}.
However, it is easiest to give a \emph{local} definition which characterises the properties of the trellis more precisely,
 and from which minimality of the number of intersections follows.
\begin{definition}[Minimal trellis]
Let $([f];X)$ be a biasymptotic mapping class.
Then a trellis map $(\widetilde{f};\widetilde{T})$ is a \emph{minimal} compatible trellis map if $\widetilde{f}\in([f];X)$ and
\begin{enumerate}\isz
\item The endpoints of $\widetilde{T}^U$ and $\widetilde{T}^S$ are points of $X$.
\item Every bigon of $\widetilde{T}$ either contains a point of $X$ in the interior of one of its boundary segments,
 or both vertices are points of $X$.
\item Every tangency of $\widetilde{T}^U$ and $\widetilde{T}^S$ is a point of $X$.
\end{enumerate}
The trellis $\widetilde{T}$ is called a \emph{minimal trellis} for $([f];X)$.
\end{definition}

We can also extract the biasymptotic orbits from a trellis.
Clearly, if $(f;T)$ is a trellis map, the set 
\[ 
  X=\bigcup_{n=-\infty}^{\infty} f^n(T^V)
\]
is a set of biasymptotic orbits, and $T$ is a minimal compatible trellis with $([f];X)$.
However, there may be smaller sets of biasymptotic orbits for which $T$ is a minimal compatible trellis.
Such a set is called a set of \emph{forcing orbits}.
\begin{definition}[Forcing]
Let $(f;T)$ be a proper trellis map. 
A set $X$ is a set of \emph{forcing orbits} for $T$ if $T$ is a minimal compatible trellis for $([f];X)$.
We say that $T$ is \emph{forced} by $X$.
\end{definition}

The rest of this section is mostly devoted to proving the existence and (essential) uniqueness of minimal trellises.


\subsection{Existence of minimal trellises}
\label{sec:minimaltrellisexistence}

To find a minimal trellis $\widetilde{T}$ for a biasymptotic mapping class $([f];X)$ with a given Birkhoff signature $\cal{B}$,
 we start with any compatible trellis map $(f_0;T_0)$ and try to remove intersection points by isotopy of maps in $([f];X)$.
Such an isotopy is called a \emph{pruning isotopy}.
\begin{definition}[Pruning isotopy]
Let $([f];X)$ be a biasymptotic mapping class, and $\sig$ be a Birkhoff signature for $([f];X)$.
Let $f_t$ be an isotopy of maps in $([f];X)$, and $T_t$ be the trellis for $f_t$ with signature $\sig$.
Then $(f_t;T_t)$ is a \emph{pruning isotopy} if for every open set $U$ of the surface $M$, and every $\tau\in[0,1]$, there is a neighbourhood $V$ of $\tau$ such that the number of intersections of $T_t$ in $U$ decreases as $t$ increases in $V$,
\end{definition}
In other words, no extra local intersections are created in $T_t$ as $t$ increases through a pruning isotopy.
An example of a pruning isotopy is shown in \figref{pruningisotopy}.
\fig{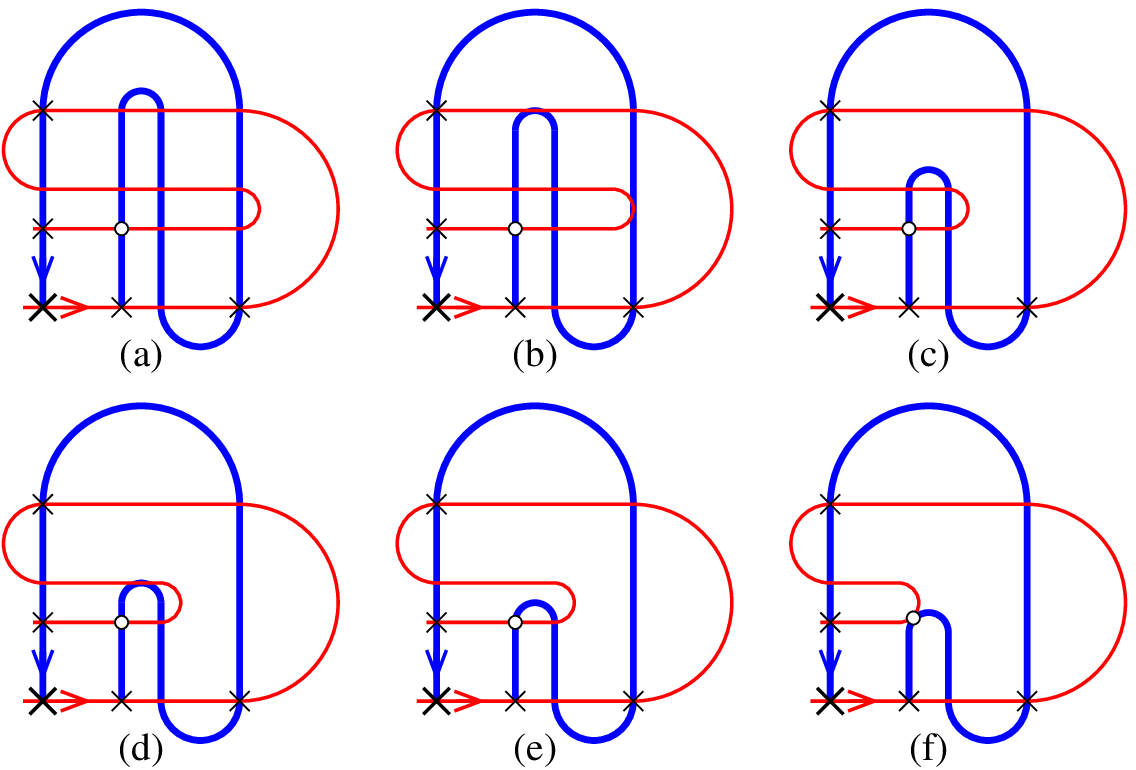}{A pruning isotopy starting with the Smale horseshoe trellis (a) and ending with a minimal trellis relative to the marked point (f).}

Assuming we always have finitely many intersections, the trellis type only changes at a tangency of $T_t^U$ and $T_t^S$.
Passing through such a tangency allows us to remove a pair of intersection points.
These intersection points must be the vertices of a special kind of bigon called an \emph{inner bigon}.
\begin{definition}[Inner bigon]
\label{defn:innerbigon}
A bigon $B$ of a trellis $T$ with connected branches is \emph{inner} if there exist positive integers $n_u$ and $n_s$
 such that both $f^{n_u}(B)$ and $f^{-n_s}(B)$ are bigons,
 but $f^{n_u+1}(B^U)\cap T^U=\emptyset$ (or is an end of $T^U$) and $f^{-(n_s+1)}(B^S)\cap T^S=\emptyset$,
 where $B^U$ and $B^S$ are respectively the open unstable and stable boundary segments of $B$.
A \emph{first inner bigon} is an inner bigon for which $n_s=0$, and a \emph{last inner bigon} is an inner bigon for which $n_u=0$, 
\end{definition}
Essentially, an inner bigon has no bigons inside it or any iterate. 
In particular, the sets $f^n(B)$ are disjoint for $-n_s\leq n\leq n_u$.
Clearly, every inner bigon is the image of some first inner bigon.
Note that the conditions preclude an inner bigon with a vertex in $T^P$.

\fig{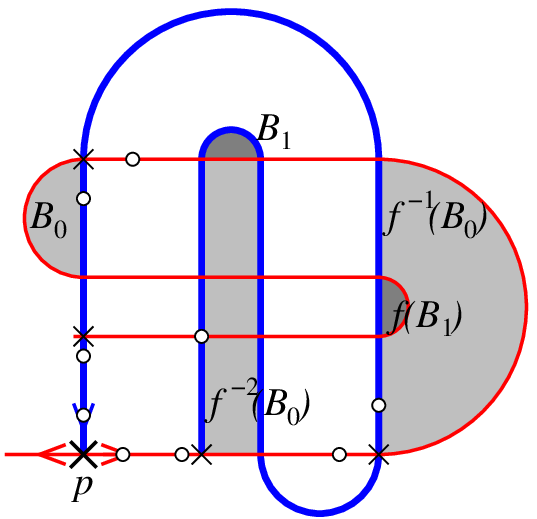}{The Smale horseshoe trellis with bigons $B_0$ and $B_1$.}
In \figref{innerbigon} we show a trellis with two labelled bigons, $B_0$ and $B_1$.
The bigon $B_1$ is a first inner bigon, but $B_0$ is not an inner bigon since $f^{-1}(B_0)$ is a domain but not a region.

The following lemma shows that if we have a trellis map with a bigon for which neither vertex is a periodic point, there must be an inner bigon.
\begin{lemma}
\label{lem:bigonexistence}
Suppose $(f;T)$ is a trellis map, where $T$ is a proper trellis with connected branches, and let $B$ be a bigon which does not contain a point of $T^P$.
Then $T$ has an inner bigon.
\end{lemma}

\begin{proof}
Let $B$ be any bigon, and let $n_s$ be the greatest integer such that $f^{-n_s}(B_0^S)$ is bounded by an interval in $T^U$ and an interval in $T^S$.
Then $f^{-n_s}(B)$ contains no stable curves in its interior, but may contain unstable curves.
Since $f^{-n_s}(B)$ is a topological disc, it must either be a bigon, or contain a bigon $B$ as a subset.
Then there exists a least integer $n_u$ such that $f^{n_u+1}(B^U)\cap T^U=\emptyset$.
Since $B$ is a bigon, $n_u<n_s$, and therefore $f^{n_u}(B)\cap T^S\subset f^{n_u-n_s}(B)\cap T^S\subset f^{n_u-n_s}(B^S)\subset f^{n_u}(B^S)$, so $f^{n_u}(B)$ is also a bigon.
Hence $B$ is an inner bigon.
\end{proof}

Henceforth we will only consider pruning isotopies which preserve the unstable set $T^U$, but allow us to change the stable set $T^S$.
Unfortunately, changes made to one piece of stable manifold may affect other pieces.
Therefore, we only consider of the form $f_t=f\circ h_t^{-1}$, where $h_t$ is supported on a set $K$ whose iterates are controlled as follows.
\begin{definition}[Simple pruning disc]
Let $([f];T)$ be a trellis mapping class.
An open disc $K$ is a \emph{simple pruning disc} if
\begin{enumerate}\setlength{\itemsep}{0pt}
\item $T^U\cap f(K)= T^S\cap f^{-N}(K)=\emptyset$, 
\item $T^U\cap f^{-N}(K) = f^{-N}(T^U\cap K)$, and 
\item the sets $K_n=f^{-n}(K)$ are disjoint for $0\leq n<N$.
\end{enumerate}
\end{definition}

The following lemma shows how precomposing with the inverse of a diffeomorphism $h$ supported in a simple pruning disc changes the stable and unstable manifolds.
\begin{lemma}
\label{lem:precompose}
Let $([f];T)$ be a trellis mapping class with signature $\sig$, and $K$ be an open subset of $M$ such that
Let $h$ be a diffeomorphism supported on $K$ such that $h$ is the identity on orbits of $\sig$.
Let $\widetilde{f}=f\circ h^{-1}$, and $\widetilde{T}$ be the trellis of $\widetilde{f}$ with signature $\sig$.
Then $\widetilde{T}^U=T^U$, $\widetilde{T}^S=T^S$ outside of $\bigcup_{n=0}^{N-1}f^{-n}(K)$, and $\widetilde{T}^S=h(T^S)$ in $K$.
\end{lemma}

\begin{proof}
Since $f^{-1}(T^U)\cap K=\emptyset$, if $x\in T^U$, then $f_t^{-1}(x) = h_t^{-1}(f^{-1}(x)) = f^{-1}(x)$, so $x\in \widetilde{T}^U$.
Hence $\widetilde{T}^U=T^U$.

$f^N(T^S)\cap K = f^N(T^S\cap f^{-N}(K)) = \emptyset$.
Therefore, if $x\in T^S$ and $f^n(x)\not\in K$ for $0\leq n<N$, then $x$ never enters $K$, so $f^m(x)=f^m(x)$ for all $m$.
Hence $x\in \widetilde{T}^S$.

If $1\leq n<N$, then 
\[ f^n(T^S\cap K)\cap K \subset f^n(K)\cap K = \emptyset , \]
 and if $N\leq n$, then 
\[ f^n(T^S\cap K)\cap K \subset f^n(T^S)\cap K \subset f^N(T^S)\cap K =\emptyset . \]
Therefore, if $x\in h(T^S)\cap K$, we have $\widetilde{f}(x)=f(h^{-1}(x)\in f(T^S\cap K)$ since $h(T^S)\cap K=h(T^S)\cap K$.
So $\widetilde{f}^n(x)=f^n(h^{-1}(x))$ for any $n>0$ since $f^n(h^{-1}(x))\not\in K$.
Hence $x\in \widetilde{T}^S$.
\end{proof}
Notice that since $K$ is a disc, the diffeomorphism $h$ is isotopic to the identity, so $\widetilde{f}$ is isotopic to $f$.

 Neighbourhoods of inner bigons and points of $T^V$ give suitable simple pruning discs.
\begin{lemma}
\label{lem:simplepruningdisc}
Let $([f];T)$ be a trellis map compatible with a biasymptotic mapping class $([f];X)$.
If $B$ is a last inner bigon such that no point of $B$ lies in $X$, then there is a simple pruning disc $K\supset B$.
Similarly, if $x$ is a point of $T^V$ such that $f(x)\not\in T^U$, there is a simple pruning disc $K\ni x$.
\end{lemma}
The proof is straightforward from the definitions of inner bigon and simple pruning disc.

We can now give the main result of this section
\begin{theorem}[Existence of minimal trellis]
\label{thm:minimalexistence}
Suppose $([f];X)$ is a biasymptotic mapping class and $\sig=(\sig^U,\sig^S)$ is a signature in $X$.
Then there is a minimal trellis map $(\widetilde{f};\widetilde{T})$ compatible with $([f];X)$ with signature $\sig$.
\end{theorem}

\begin{proof}
Suppose $(f;T)$ is not minimal.
Then either $T$ has a tangency $v$ which is not in $X$, or by Lemma~\ref{lem:bigonexistence} we can find a last inner bigon $B$
 such that at least one vertex of $B$ is not a point of $X$.

\fig{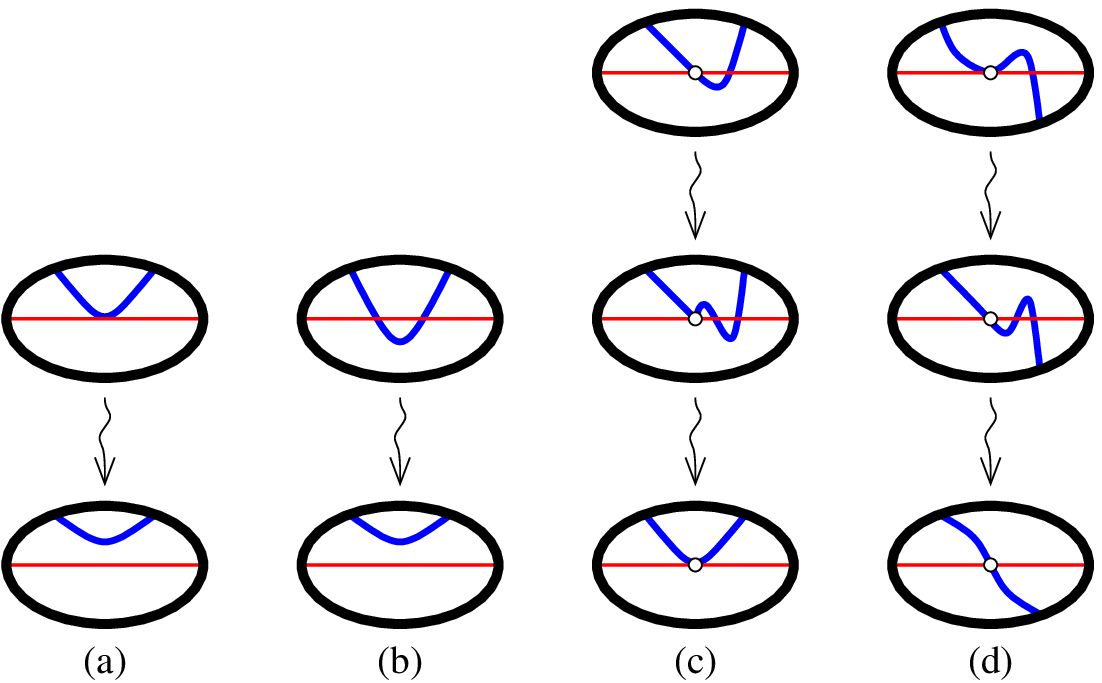}{Isotopies removing intersections in $K$: (a) removing intersections from a bigon, (b) removing a tangency, 
 (c) making a tangency at a point of $X$ and (d) moving a transverse intersection to a point of $X$.}

If $v$ is a tangency point which is not in $X$, then there exists a greatest integer $n_u$ such the $f^{n_u}(v)\in T^U$.
By Lemma~\ref{lem:simplepruningdisc} there is a simple pruning disc $K$ containing $f^{n_u}(v)$ and no other intersections of $T$.
Take $h$ such that $h(T^S\cap K)$ does not intersect $T^U$ and $\widetilde{f}=f\circ h^{-1}$.
Then by Lemma~\ref{lem:precompose}, $\widetilde{T}$ has fewer intersection points than $T$.
Then $f^{n_u}(B)\cap T^U=f^{-(n_s+1)}(B)\cap T^S=\emptyset$ since the end vertices of $T$ are points of $X$.
This case is shown in Figure~\figrefpart{bigonisotopy}{a}.

Similarly, suppose $T$ has a last inner bigon $B$ such that neither vertex of $B$ is a point of $X$.
Then by Lemma~\ref{lem:simplepruningdisc}, there is a simple pruning disc $K$ containing $B$, and no points of $T^V$ other than the vertices of $B$, 
so we can isotopy $f$ to $\widetilde{f}$ with a trellis $\widetilde{T}$ which has fewer intersection points than $T$.
Then $f^{n_u}(B)\cap T^U=f^{-(n_s+1)}(B)\cap T^S=\emptyset$ since the end vertices of $T$ are points of $X$.
Since the sets $f^{n}(B)$ are disjoint for $-n_s\leq n\leq n_u$, we can find a contractible neighbourhood $K$ of $f^{n_u}B$
 satisfying the conditions of Lemma~\ref{lem:precompose}.
Then, as before, we can find isotope $f$ to obtain a trellis with fewer intersections.
This case is shown in Figure~\figrefpart{bigonisotopy}{b}.

It remains to consider the cases where there is a last inner bigon $B$ for which at least one vertex $x$ is a point of $X$.
Now, it may be the case that $f(x)\in T^U$ or $f^{-(n_s+1)}(x)\in T^S$, so Lemma~\ref{lem:precompose} does not apply directly to $B$.
In this case, we can first precompose with an isotopy supported in a neighbourhood of $f^(x)$ to obtain an inner bigon whose vertices are not points of $X$,
 as shown in \figrefpart{bigonisotopy}{c} and \figrefpart{bigonisotopy}{d}, and then remove this inner bigon as previously discussed.
This again reduces the number of intersection points.

The result follows by induction on the number of vertices, since we can always obtain a trellis with finitely many vertices by perturbing $f$.
Note that this construction does not change the unstable set $T^U$.
\end{proof}


\subsection{Uniqueness of minimal trellis}
\label{sec:orbit}

We now show that a minimal trellis compatible with a given biasymptotic mapping class with a given signature is essentially unique.
Unfortunately, there is one problem; if there are two points $x_0$ and $x_1$ of $X$ which form the intersections of an inner bigon, the orientation of these intersections is unspecified by the minimality conditions.
Note that if there is a bigon with both vertices in $X$, we may have different local behaviour as shown in \figref{bigonnonunique}.
\fig{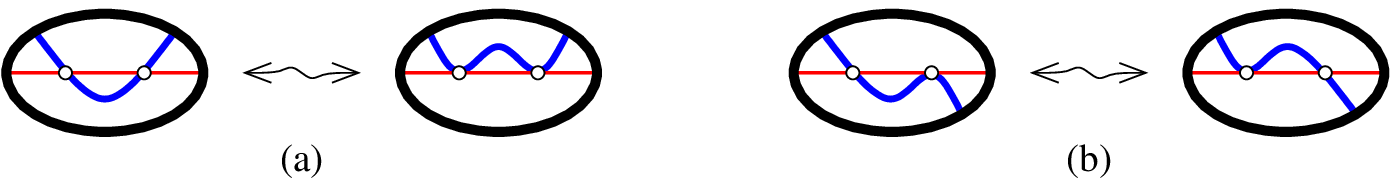}{Possibilities for bigons with both vertices in $X$.}

As shown in \figref{smalenonunique}, there are minimal compatible trellises with transverse and tangential crossings at $x_0$ and $x_1$.
\fig{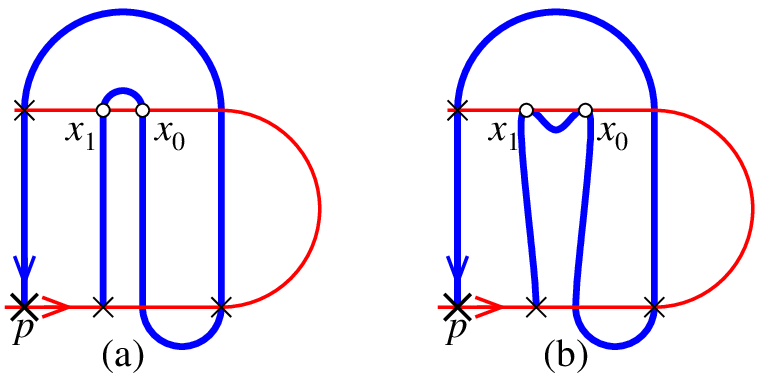}{Different minimal compatible trellises.}

We use the following preliminary theorem, which is derived from Theorem~\ref{thm:twoisotopy} and Theorem~\ref{thm:minimalcurve}.
\begin{proposition}
\label{prop:orbitspec}
Suppose $(f_0;T)$ and $(f_1;T)$ are minimal trellis maps compatible with a biasymptotic mapping class $([f];X)$.
Then $(f_0;T)$ and $(f_1;T)$ are in the same trellis mapping class.
\end{proposition}

\begin{proof}
Let $h=f_0^{-1}\circ f_1$.
Then $h$ is isotopic to $\id$ relative to $X$.
Since $T^U$ consists of cross-cuts to $X$, there is an isotopy $\widetilde{h}_t$ such that
 $\widetilde{h}_0=\id$, $\widetilde{h}_1=h$ and $\widetilde{h}_t(T^U)=T^U$ for all $t$.
If we now cut along $T^U$, we find that $h$ is isotopic to the identity by the isotopy $\widetilde{h}$, and $T^S$ is a set of cross-cuts.
Therefore, where is an isotopy $h_t$ such that $h_0=\id$, $h_1=h$, and for all values of $t$, $h_t(T^U)=T^U$ and $h_t(T^S)=T^S$.
To show that $f_t=f_0\circ h_t$ gives the required isotopy from $f_0$ to $f_1$, we need only check that $f_t$ preserves $T^U$ and $T^S$.
We have $f_t^{-1}(T^U)=h_t^{-1}(f_0^{-1}(T^U))=f_0^{-1}(T^U)$ since $f_0^{-1}(T^U)\subset T^U$, and $f_t(T^S)=f_0(h_t(T^S))=f_0(T^S)$, as required.
\end{proof}

This result is important in its own right, since it shows that as long as the topology of the surface $M$ is simple,
 we can specify a trellis mapping class by giving the trellis up to homeomorphism, and then only need to give one biasymptotic orbit on each branch.
This greatly simplifies the representation of a trellis mapping class required for computations,
 and justifies our custom of specifying a trellis mapping class by drawing the trellis and part of a biasymptotic orbit.

We can now show that a minimal trellis compatible with a set of biasymptotic orbits is essentially unique.

\begin{theorem}[Uniqueness of minimal trellis]
\label{thm:minimaluniqueness}
Let $([f];X)$ be a biasymptotic mapping class and $\sig$ a Birkhoff signature for $X$.
Let $(f_0;T_0)$ and $(f_1;T_1)$ be trellis maps forced by $([f];X)$, with end intersections at points of $\sig$.
Suppose that whenever there are points $x_1$ and $x_2$ of $X$ such that $T^U[x_1,x_2]\cup T^S[x_1,x_2]$ bounds a bigon $B$ in $T_0$,
 and there are no other points of $X$ on the boundary of $B$, that $T_0$ and $T_1$ have the same type of crossing (transverse or tangential) at $x_1$ and $x_2$.
Then $(f_0;T_0)$ and $(f_1;T_1)$ have the same trellis type.
\end{theorem}

\begin{proof}
Since the curves $T^U$ and $T^S$ are mutually homotopic, have the same types of intersection at $X$, 
 and have minimal intersections with respect to $X$, by Theorem~\ref{thm:minimalcurve}, 
 there is a homeomorphism $h$ such that $h(T_0)=T_1$ and $h$ is isotopic to the identity relative to $X$.
Then $(h^{-1}\circ f_1\circ h;h^{-1}\circ T_1)=(h^{-1}\circ f_1\circ h;T_0)$, and $f_0$ and $h^{-1}\circ f_1\circ h$ are isotopic
 relative to $X$.
Then by Proposition~\ref{prop:orbitspec}, $(f_0;T_0)$ and $(h^{-1}\circ f_1\circ h;T_0)$ are in the same trellis mapping class,
 so $(f_0;T_0)$ and $(f_1;T_1)$ have the same type.
\end{proof}


\subsection{Minimal extensions}
\label{sec:minimalextension}

We have considered how to construct a from a set of biasymptotic orbits a trellis with a given signature.
We can also extend trellises to longer trellises.
Of particular interest are those which do not introduce any more intersections than necessary.
Such trellises are called \emph{minimal iterates}  and \emph{minimal extensions}.
They are most simply characterised as being extensions which are minimal with respect to the orbits of the original trellis.
\begin{definition}[Minimal iterate and minimal extension]
Let $([f];T)$ be a trellis mapping class.
An iterate $(\widetilde{f};\widetilde{T})$ of $([f];T)$ is a \emph{minimal iterate} of $([f];T)$ if $\widetilde{T}$ is forced by $X$,
 where the set $X$ is given by $X=\bigcup_{n\in\Z}\widetilde{f}^n(T^V)$, the union of the orbits of the intersections of $T$.
An extension $(\widetilde{f};\widetilde{T})$ of $([f];T)$ is a \emph{minimal extension} of $([f];T)$ if it is a subtrellis of some minimal iterate.
\end{definition}

The following result follows almost directly from the definition of a minimal extension. 
\begin{theorem}[Existence and uniqueness of minimal extensions]
\label{thm:minimalextension}
Let $([f];T)$ be a trellis mapping class, $X=\bigcup_{n\in\Z}f^n(T^V)$, $\sig$ the signature class of $T$ in $X$,
 and $\widetilde{\sig}$ a signature class in $X$ with $\widetilde{\sig}\geq\sig$.
Then there is a unique trellis type $[\widetilde{f};\widetilde{T}]$ with signature $\widetilde{\sig}$
 which is a minimal extension of $([f];T)$.
Further, if $T$ is a transverse trellis, so is $\widetilde{T}$.
\end{theorem}
The proof is essentially the same as that of Theorem~\ref{thm:minimalexistence} and Theorem~\ref{thm:minimaluniqueness}.
Note that any inner bigon of $\widetilde{T}$ is an iterate of some inner bigon of $T$, so the intersection type is given by $T$.

The following trivial result gives a characterisation of certain minimal extensions.
\begin{lemma}
\label{lem:minimalextensiongoodend}
Let $(\widetilde{f};\widetilde{T})$ be a minimal extension of $([f];T)$.
Then every bigon of $\widetilde{T}$ contains a point of $\widetilde{X}=\bigcup_{n=-\infty}^{\infty}f^n(T^V)$.
Further, if $(\widetilde{f};\widetilde{T})$ is any extension for which the end intersections in $\widetilde{X}$, and every bigon of $\widetilde{T}$ contains a point of $\widetilde{X}$, then $(\widetilde{f};\widetilde{T})$ is a minimal extension.
\end{lemma}
Note that it is possible to find extensions of $([f];T)$ for which every bigon contains a point of $\widetilde{X}$, but which is not a minimal extension.
By Lemma~\ref{lem:minimalextensiongoodend}, this may only occur if there is an end intersection not in $\widetilde{X}$.

The definition of a minimal supertrellis is complicated since we need to ensure that the extra periodic points correspond to orbits of the original trellis mapping class.
\begin{definition}
Let $([f];T)$ be a trellis mapping class.
Then a supertrellis $(\widetilde{f};\widetilde{T})$ is a \emph{minimal supertrellis} if for every extension of $(\widehat{f};\widehat{T})$ of $([\widetilde{f}];\widetilde{T})$, every bigon of $\widehat{T}$ contains a point of $\widehat{X}=\bigcup_{n=-\infty}^{\infty}\widehat{f}^n(T^V)$.
\end{definition}
In particular, this means that every point of $\widetilde{T}^P\setminus T^P$ shadows an essential periodic orbit of $([f];T)$.

Since minimal extensions and minimal supertrellises have the same forcing orbits as the original trellis,
 we expect the dynamics to be the same as that forced by the original trellis, or indeed, the dynamics forced by the biasymptotic orbits themselves.
This is in fact the case, as we shall see later.
Minimal extensions provide a way of getting more information about the dynamics without increasing entropy.


\section{Graphs Representatives}
\label{sec:graph}

Any connected compact surface with nonempty boundary is homotopy-equivalent to a one-dimensional space.
Further, such spaces, and maps on them, are very easy to describe combinatorially.
We can use these properties to provide a framework for representing surface diffeomorphisms and computing their dynamical properties.
In classical Nielsen-Thurston theory, surface homeomorphisms are represented on one-dimensional spaces
 with a differentiable structure called \emph{train tracks}.
When representing trellises, we also need to take into account the topology of the stable and unstable curves.
By cutting along the unstable curves, we introduce new loops in the graph. 
The stable curves are represented by special edges called \emph{control edges},
 since these control the behaviour of the one-dimensional representative of a trellis map.


\subsection{Combinatorics of graph maps}
\label{sec:graphcombinatorics}

We first give some standard definitions concerning graphs embedded in compact surfaces.

\begin{definition}[Graph]
A \emph{graph} $G$ is an one-dimensional CW-complex.
The \emph{vertices} of $G$ are the zero-dimensional cells, and the \emph{edges} are the one-dimensional cells.
The edges are oriented, and the reverse of the edge $e$ is denoted $\bar{e}$, so $\bar{\bar{\mbox{$e$}}}=e$.
The initial vertex of a oriented edge $e$ will be denoted $\init{e}$.
The \emph{valence} of a vertex $v$ is the number of directed edges for which $v$ is the initial vertex.
A graph can be described up to homeomorphism by giving its vertices and the incident edges at each vertex.
\end{definition}

The \emph{Euler characteristic} $\chi$ of a graph is $\card{\vertex(G)}-\card{\edge(G)}$
 where $\card{\vertex(G)}$ is the total number of vertices and $\card{\edge(G)}$ is the total number of edges.
$\chi$ is an homotopy type invariant.

\begin{definition}[Edge-path]
An \emph{edge-path} is a list $e_1\ldots e_n$ of oriented edges of $G$ such that $\final{e_i}=\init{e_{i+1}}$ for $1\leq i < n$.
An \emph{edge loop} is a cyclically-ordered list of edges.
The \emph{trivial} edge-path contains no edges and is denoted $\cdot$.
An edge-path $e_1\ldots e_n$ \emph{back-tracks} if $e_{i+1}=\bar{e}_i$ for some $i$, otherwise it is \emph{tight}.
\end{definition}

\begin{definition}[Surface embedding]
We will always consider a graph embedded in a surface by an embedding $i$.
This induces a natural cyclic order $\rhd$ on the oriented edges starting at each vertex.
The graph can be described up to local ambient homeomorphism by giving the cyclic ordering of the incident edges at each vertex.
\end{definition}
\begin{definition}[Turn]
A pair of edges $(e_1,e_2)$ is a \emph{turn} in $G$ at at vertex $v$ if $v=\init{e_1}=\init{e_2}$ and $e_2\rhd e_1$, so $e_2$ immediately follows $e_1$ in the cyclic order at $v$.
\end{definition}

\begin{definition}[Peripheral loop]
An edge-loop $\pi=\ldots,p_1,p_2,\ldots,p_n,\ldots$ is \emph{peripheral} in $G$ if $(p_{i+1},\bar{p}_i)$ is a turn in $G$ for all $i$.
A peripheral loop $\pi$ is \emph{simple} if $\pi$ is a simple closed curve, or equivalently, if $\init{p_i}\neq\init{p_j}$ for $i\neq j$.
Peripheral loops which are not simple may repeat vertices, or even edges (with opposite orientation).
Peripheral loops define anti-clockwise curves around boundary components of the surface $M$.
\end{definition}

\begin{definition}[Frontier and link]
If $H$ is a subgraph of $G$, an edge-loop $\pi$ is a \emph{frontier loop} of $H$ if $\pi$ is peripheral in $H$ but not in $G$.
If $H$ is a subgraph of $G$, then the \emph{link} of $H$, denoted $\link{H}$ is the set of edges $e$ of $G$ with $\init{e}\in H$, but $e\not\in H$.
\end{definition}

\begin{definition}[Graph maps]
A \emph{graph map} $g$ is a self-map of $G$ taking a vertex to a vertex, and edge $e$ to an edge-path $e_1\ldots e_k$
 such that $\init{e_1}=g(\init{e})$ for all directed edges $e$.
If $g$ is a graph map on $G$, where $G$ is embedded in $M$ by an embedding $i$,
 then we say $g$ is \emph{embeddable} if there are arbitrarily small perturbations of $i\circ g$
 which are embeddings and preserve the cyclic order at any vertex.
(A combinatorial definition of an embeddable graph map also exists.)
\end{definition}
\begin{definition}[Derivative map]
The \emph{derivative} map $\partial g$ takes oriented edges to oriented edges or $\cdot$, with $\partial g(e_i)=e_j$
 if $g(e_i)=e_j\ldots$ and $\partial g(e_i)=\cdot$ if $g(e_i)=\cdot$.
\end{definition}

\begin{definition}[Peripheral subgraph]
The \emph{peripheral subgraph} $P$ of $g$ is a maximal invariant subset of $G$ consisting of simple peripheral loops.
Edges of $P$ are called \emph{peripheral edges}.
If $g^n(e)\subset P$ for some $n$, then $e$ is \emph{pre-peripheral}.
The set of pre-peripheral edges is denoted $\pre{P}$, and contains $P$.
\end{definition}

The dynamics of $g$ can be split into \emph{transitive components}.
\begin{definition}[Transitive and irreducible component]
Define a transitive relation on the edges $E$ of $G$ by $e_1\geq e_2$ if there exists an $n$ such that $g^n(e_1)$ contains the edge $e_2$,
 and a stronger transitive relation  by $e_1\succcurlyeq e_2$ if there exists $n$ such that $g^m(e_1)$ contains the edge $e_2$ for all $m\geq n$.
A \emph{transitive component} of $g$ is an equivalence class under $\geq$, and an \emph{irreducible component} of $g$ is an equivalence class under $\succeq$.
\end{definition}
The quotient of these relation by the equivalence classes give a partial orders on the transitive/irreducible components.
A transitive component of a graph map $g$ is a union of irreducible components which are cyclically permuted by $g$.

\begin{definition}[Transition matrix and growth rates]
The \emph{transition matrix} of $g$ is the matrix $A=(a_{ij})$ where $a_{ij}$ is the number of times the edge $e_j$ appears in the image path of edge $e_i$.
The largest eigenvalue of $A$ is the \emph{growth rate} $\lambda$ of $g$, and the logarithm of the growth rate gives the topological entropy of $g$.
\end{definition}
We can compute growth rates for any transitive component, and the growth rate for the entire graph map is the maximum of the growth rates on each component.

\begin{example}
\figref{graph} shows a graph embedded in a disc with four punctures.
The peripheral loops are $p_1$, $p_2$, $p_3$, $p_4p_5$ (surrounding the punctures) and
 $a\bar{p}_1\bar{a}\bar{p}_4c\bar{p}_3\bar{c}\bar{p}_5b\bar{p}_2\bar{b}$.
\fig{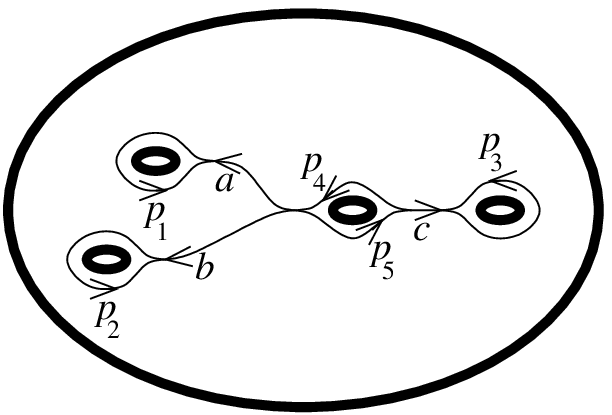}{A graph embedded in a disc with four punctures.}

There is a graph map $g$ such that
\[ \begin{array}{c} g(p_1)=p_2,\ g(p_2)=p_3,\ g(p_3)=p_1,\ g(p_4)=p_5,\ g(p_5)=p_4, \\[1ex]
     g(a) = c\bar{p}_3\bar{c}\bar{p}_5b,\ g(b)=c\bar{p}_3\bar{c}\bar{p}_5b\bar{p}_2\bar{b}\bar{p}_5c,\ g(c)=a. \end{array} \]
The peripheral subgraph of $g$ is $\{ p_1,p_2,p_3,p_4,p_5 \}$.
The link of the peripheral subgraph consists of the oriented edges $\{a,\bar{a},b,\bar{b},c,\bar{c}\}$.
\end{example}

When relating graph maps to trellis types in Section~\ref{sec:compatiblegraph} we need to consider topological pairs to capture the configuration of $T^S$.
\begin{definition}[Control edges]
Let $(G,W)$ be a topological pair, where $W$ is a finite set of points such that every edge of $G$ contains at most one point of $W$.
Then the edges of $G$ containing a point of $W$ are called \emph{control edges}; other edges of $G$ are \emph{free edges}.
The set of control edges of $G$ is denoted $Z(G)$.
The pair $(G,W)$ is a \emph{controlled graph}.
A graph map $g:(G,W)\fto(G,W)$ is a \emph{controlled graph map} if for every control edge $z_0$, the image of $g(z_0)$ is a control edge $z_1$.
We denote a controlled graph map $g$ of $(G,W)$ by $(g;G,W)$.
\end{definition}
A vertex of $G$ which is the endpoint of a control edge is called a \emph{control vertex}.
All other vertices are called \emph{free vertices}, and edges which are not control edges are called \emph{free edges}.
Free edges which are neither peripheral nor pre-peripheral are called \emph{expanding edges.}

A controlled graph is \emph{proper} if it is connected and free vertex has valence at least $3$.
A proper controlled graph has at most $3\card{\control(G)}-3\chi$ edges, and $3\card{\control(G)}-2\chi$ vertices, where $\card{\control(G)}$ is the number of control edges.


\subsection{Compatible graph maps}
\label{sec:compatiblegraph}

Our main tool for computing and describing the dynamics forced by a trellis mapping class or trellis type is to relate the trellis map to a graph map.
We do not relate a trellis mapping class directly to a graph map, but instead use the map $\cut{f}$ induced by cutting along the unstable curves.
To relate this map to a map on a different space, we introduce a new concept of equivalence, that of \emph{exact homotopy equivalence},
 which is the natural extension of the concept of homotopy equivalence to maps of pairs.

We first need to define what it means for a map of pairs to be \emph{exact}.
\begin{definition}[Exact map of pairs]
\label{defn:exactmap}
A map of topological pairs $f:(A,B)\fto (X,Y)$ is \emph{exact} if $f^{-1}(Y)=B$.
(Equivalently, $f$ is exact if $f(B^C)\subset Y^C$.)
\end{definition}
Note that this is equivalent to saying that $f$ is a map of triples, $f:(A,B,B^C)\fto(X,Y,Y^C)$.
However, we do not use this notation, as we treat exact maps as a special class of map in the category of topological pairs.
In particular, we will often take homotopies between exact maps and non-exact maps in the category of topological pairs.

We can compare systems on different spaces via exact homotopy equivalence.
\begin{definition}[Exact homotopy equivalence]
\label{defn:exacthomotopy}
Topological pairs $(A,B)$ and $(X,Y)$ are \emph{exact homotopy equivalent} if there are exact maps $p:(A,B)\exto(X,Y)$ and $q:(X,Y)\exto(A,B)$
 such that $q\circ p\homotopic \id_{A}$ and $p\circ q\homotopic \id_{X}$. 
$p$ and $q$ are called \emph{exact homotopy equivalences}.
Maps $f:(A,B)\fto(A,B)$ and $g:(X,Y)\fto(X,Y)$ are \emph{exact homotopy equivalent} if there is an exact homotopy equivalence $p:(A,B)\exto(X,Y)$ such that $p\circ f=g\circ p$.
\end{definition}
Note that the homotopies used in the definition of exact homotopy equivalence can be taken through any map of pairs and not just exact maps. 

A graph map representing the topology of a trellis via exact homotopy equivalence is called \emph{compatible} with the trellis.
\begin{definition}[Compatible graph map]
Let $(G,W)$ be a controlled graph.
The $(G,W)$ is \emph{compatible} with a transverse trellis $T$ if $(G,W)$ and $(\cut[U]{T},T^S)$ are exact homotopy equivalent by an embedding $i:(G,W)\exto(\cut[U]{T},T^S)$, and $G$ crosses $T^S$ transversely.
A controlled graph map $g$ of $(G,W)$ is \emph{compatible} with the trellis type $[f;T]$ if the embedding $i$ is an exact homotopy equivalence between $g$ and $\cut{f}$, and $g$ and $\cut{f}$ have the same orientation at points of $W$.
\end{definition}
Note that the inclusion $i$ induced a bijection between the regions of $\cut{T}$ (and hence of $T$) with the regions of $(G,W)$,
 and that all compatible graphs are exact homotopy equivalent.
We restrict to transverse trellises since a trellis with tangencies may not have a compatible controlled graph.
It is always possible to find a topological pair $(G,H)$ which is exact homotopy equivalent to $(\cut[U]{T},T^S)$ for which $G$ is a graph,
 but it may not be possible to take $H$ to be a finite set of points.

\fig{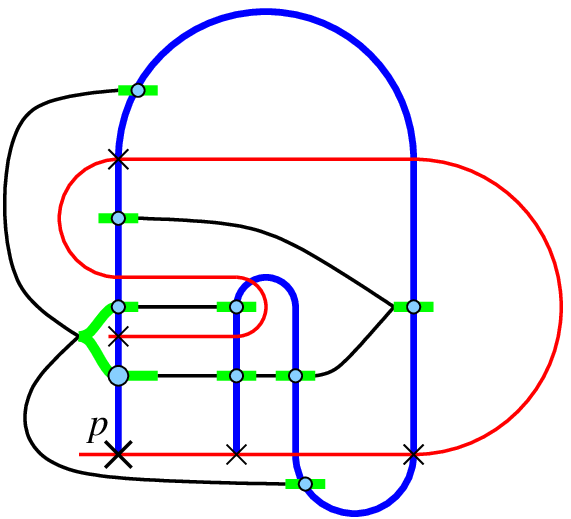}{The controlled graph $(G,W)$ is compatible with the trellis $T$.}
The controlled graph shown in \figref{henonembedded} is compatible with the trellis $T$.
The control edges are shown as thick green lines.

There are many controlled graph maps compatible with a trellis type $[f;T]$.
To use graph maps to describe the dynamics forced by $[f;T]$, we need to define a subclass of controlled graph maps which have minimal entropy in the exact homotopy class.
These graph maps are called \emph{efficient}, by analogy with Nielsen-Thurston theory.
Of the efficient graph maps, we show in Section~\ref{sec:graphuniqueness} that there is at most one which is \emph{optimal}, giving a canonical \emph{graph representative} for the trellis type.

Efficient and optimal graph maps can be defined in terms of their actions on the turns of $G$.
Some of following notions are used in the definitions, and all are useful in describing the algorithms.
\begin{definition}[Bad turns]
Let $g$ be a graph map of $(G,W)$, and $(e_0,e_1)$ a turn of $G$. 
Then
\begin{enumerate}\isz\vspace*{-\baselineskip}
\item If $\partial g(e_0)\neq\partial g(e_1)$, the turn is a \emph{good} turn.
\item If $\partial g(e_0)=\partial g(e_1)$ and at least one of $e_0$ and $e_1$ is a control edge, the turn is a \emph{controlled} turn. 
  If both $e_0$ and $e_1$ are control edges, the turn is \emph{fully controlled}, if only one is, the turn is \emph{half controlled}.
\item If $\partial g(e_0)=\partial g(e_1)$ and both $e_0$ or $e_1$ are free edges, the turn is a \emph{bad} turn.
\item A bad turn $(e_0,e_1)$ is \emph{inefficient} if in addition there is an edge $e$ such that $g^n(e)=\ldots \bar{e}_0e_1\ldots$.
\end{enumerate}
\end{definition}

In order to find efficient graph maps, we pass through intermediate graphs where there are no obvious ways of reducing the dynamics locally.
Such maps are called \emph{tight}.
\begin{definition}[Tight graph map]
Let $g$ be a graph map of a controlled graph $G$.
Then $g$ is \emph{vertex tight} if there is no free vertex $v$ such that the derivative map $\partial g$ takes the same value for all oriented edges $e$ with $\init{e}=v$,
and $g$ is \emph{edge tight} if for every edge $e$, the edge path $g(e)$ is nontrivial and does not back-track.
$g$ is \emph{tight} if it is both vertex tight and edge tight.
\end{definition}
Note that $\partial g$ is constant at a vertex $v$ if there is an edge $e$ such that $g(e_i)=e\ldots$ for all edges $e_i$ at $v$.

\begin{definition}[Efficient and optimal graph maps]
A controlled graph map $g$ is \emph{efficient} if there are no inefficient turns, and \emph{optimal} if there are no bad turns, and every invariant forest contains a control edge.
\end{definition}

A \emph{graph representative} of$[f;T]$ is then an optimal controlled graph map compatible with $([f];T)$.
\begin{definition}[Graph representative]
A controlled graph map $(g;G,W)$ is a \emph{graph representative} of an transverse trellis type $[f;T]$ if $g$ is an optimal graph map which is compatible with $([f];T)$.
\end{definition}

The main theorem concerning the representation of trellis types by graph maps is stated below.
\begin{theorem}[Existence and uniqueness of graph representatives]
Let $[f;T]$ be a proper trellis type with no invariant curve reduction.
Then $[f;T]$ has a unique graph representative $(g;G,W)$.
Further, if $[f_0;T_0]$ and $[f_1;T_1]$ are trellis types, the graph representatives $(g_0;G_0,G_0)$ and $(g_1;G_1,W_1)$ are homeomorphic, then $[f_0;T_0]=[f_1;T_1]$.
\end{theorem}
We prove the existence of a graph representative in Section~\ref{sec:algorithm} by giving an algorithm to compute it, and uniqueness in Section~\ref{sec:graphuniqueness}.
Since the maps $\cut{f_0}$ and $\cut{f_1}$ for different trellis types are not exact homotopy equivalent, it is trivial that different trellis types have different graph representatives.
This means that the graph representative provides a convenient way of specifying a trellis type.

If we are interested in the dynamics forced by the trellis mapping class, we can restrict attention to the subgraph of the graph representative $(g;G,W)$ containing the nontrivial dynamics.
This subgraph is called the \emph{essential graph representative}.
\begin{definition}[Essential graph representative]
Let $(g;G,W)$ be the graph representative of a trellis type $[f;T]$.
Then the \emph{essential graph representative} of $[f;T]$ is the controlled graph map $(g;\overline{G},\overline{W})$, where $\overline{G}=\bigcap_{n=0}^{\infty} g^n(G)$, and $\overline{W}=\overline{G}\cap W$.
\end{definition}
We can further simplify the representation by collapsing all control edges to points.
This gives the \emph{topological graph representative}.
The topological graph representative is an invariant of the forcing orbits.
\begin{definition}[Topological graph representative]
Let $(g;G,W)$ be the graph representative of a trellis type $[f;T]$.
Then the \emph{topological graph representative} of $[f;T]$ is the topological conjugacy class of graph map obtained by collapsing all control edges of the essential graph representative to points.
\end{definition}

\begin{example}
The graph representatives for the trellis type shown in \figref{henonembedded} are shown in \figref{henongraph}.
\fig{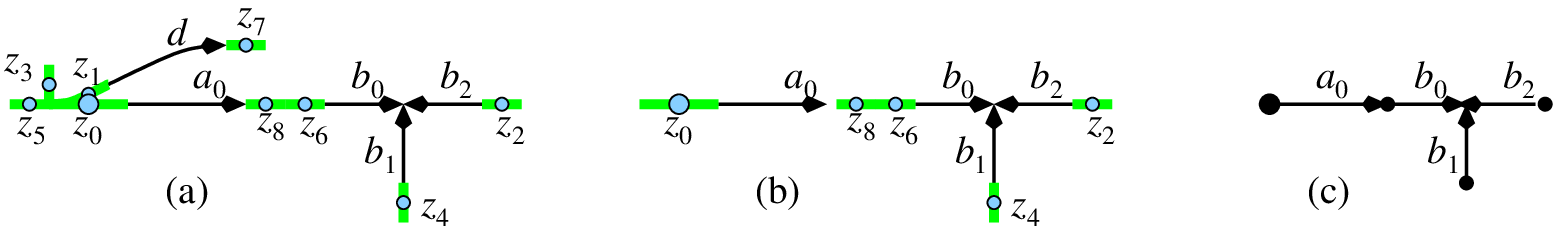}{Graph representatives for the trellis type shown in \figref{henonembedded}. (a) is the graph representative, (b) the reduced graph representative and (c) the topological graph representative.}

The control edges map
\[ z_0,z_1,z_2\mapsto z_0; \quad z_3\mapsto z_1; \quad z_4\mapsto z_2; \quad z_5\mapsto z_3; \quad z_6,z_7,z_8\mapsto z_4; \]
and the expanding edges map 
\[ a_0\mapsto a_0\bar{z}_8z_6b_0\bar{b}_1; \quad b_0\mapsto b_1; \quad b_1\mapsto b_2; \quad b_2\mapsto a_0\bar{z}_8z_6b_0; \quad d\mapsto a_0\bar{z}_8z_6b_0\bar{b}_1. \]
The edges $z_1$, $z_3$, $z_5$, $z_7$ and $d$ do not lie in $g^2(G)$ so are not part of the essential graph representative shown in \figrefpart{henonembedded}{b}.
Collapsing the control edges gives the graph shown in \figrefpart{henonembedded}{c}, with
\[ a_0\mapsto a_0b_0\bar{b}_1; \quad b_0\mapsto b_1; \quad b_1\mapsto b_2; \quad b_2\mapsto a_0b_0. \]
Joining the edges $a_0$ and $b_0$ to give a single edge $e_0$, and letting $e_1=b_1$ and $e_2=b_2$ gives the topologically conjugate graph map
\[ e_0\mapsto e_0\bar{e}_1e_1; \quad e_1\mapsto e_2; \quad  e_2\mapsto e_0. \]
The transition matrix for this graph map is 
\[ \left(\begin{array}{ccc}1&2&0\\0&0&1\\1&0&0\end{array}\right) \]
and the growth rate is given by the largest root of the characteristic polynomial $\lambda^3-\lambda^2-2=0$.
Numerically, the growth rate is $1.696$, giving topological entropy $\htop(g)\approx 0.528$.
\end{example}
  

\subsection{Properties of graph representatives}
\label{sec:graphproperties}

We now give some elementary properties concerning relationship between a trellis and a compatible graph, and of the dynamics of a graph representative.
Our main tool will be to look at curves in $\cut{T}$ with endpoints in $T^S$ which cross $T^S$ transversely.
These curves can be projected to $(G,W)$ by the exact homotopy retract $r:(\cut[U]{T},T^S)\exto(G,W)$.
Since the inclusion $i:(G,W)\exto(\cut[U]{T},T^S)$ is an exact homotopy inverse for $r$, there is a natural equivalence between exact homotopy classes of curves
 in $(\cut[U]{T},T^S)$ and in $(G,W)$.

\fig{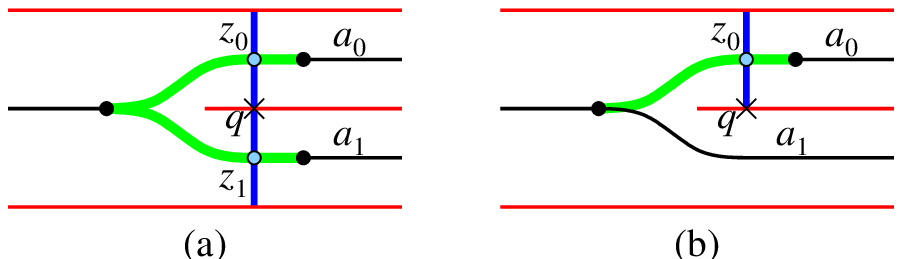}{Cusps at the end of an unstable curve. (a) shows a cusp where the last intersection is a transverse crossing.
   (b) shows a cusp where the last intersection is an endpoint of $T^S$.}
The ends of $T^U$ give rise to \emph{cusps} in $(G,W)$, as shown in \figref{cusp}. 
In \figrefpart{cusp}{a}, the control edges $z_0$ and $z_1$ surround an end of $T^U$.
Then $g(z_0)$ and $g(z_1)$ must be the same, unless the intersection point $q$ is a point of $T^P$.
In \figrefpart{cusp}{b}, the point $q$ is an endpoint of $T^S$.
In this case, the initial edge of $g(a_1)$ is $g(z_1)$, and the turn is only half controlled.

\fig{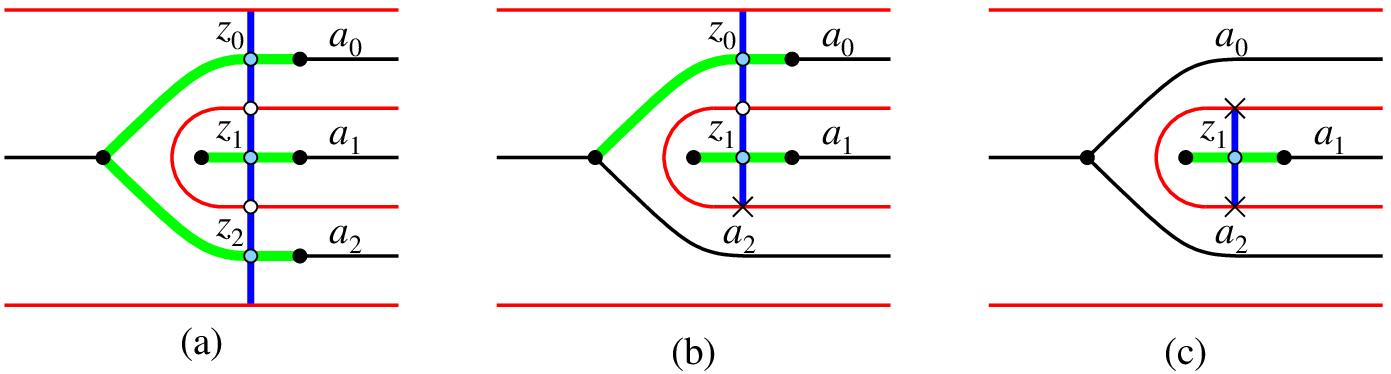}{Turns surrounding a bigon.}
If $z_0$ and $z_1$ form a cusp, then unless $z_0$ and $z_1$ are periodic, in which case the cusp forms at the end of an open branch of $T^U$, we have $g(z_0)=g(z_1)$.
The only other possibility for a turn where both edges are control edges is shown in \figrefpart{turn}{a}.
Here, the control edges $z_0$ and $z_2$ are parallel to the unstable boundary $B^U$ of a bigon $B$ whose stable boundary crosses a control edge $z_1$.
If $f(B^U)\cap T^U=\emptyset$, then $g(z_0)=g(z_1)=g(z_2)$.
Removing the stable boundary segments $S_2$ crossing $z_2$ as shown in \figrefpart{turn}{b} does not significantly change the topology of the graph;
 in particular, the valence-$3$ vertex remains a control vertex.
However, if we also remove the stable boundary segment $S_0$ crossing $z_0$, this vertex is no longer a control vertex, so the resulting graph map is no longer optimal.
Indeed, the graph representative for the resulting trellis may not be surface embeddable, and its entropy need not equal $\htop[f;T]$.
This is a consequence of the failure of Lemma~\ref{lem:wellformed} for minimal iterates of curves.

\fig{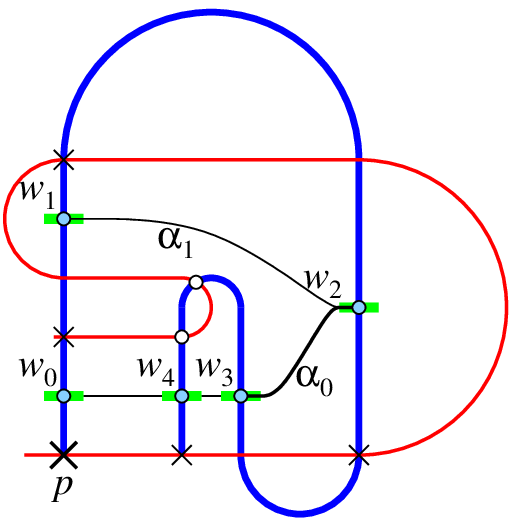}{A minimal iterate of a curve in a graph.}
Now consider the graph shown in \figref{henonembeddediterate}.
The curve $\alpha_0$ from $w_3$ to $w_2$ maps to $\alpha_1$ from $w_2$ to $w_1$.
We shall see that if the curve $\alpha_0$ does not back-track except in control edges, and $g$ is an efficient graph map, then $g\circ\alpha_0$ is a minimal iterate.
This will allow us to prove the existence of (bi)asymptotic orbits directly from the graph map.
Also, by relating minimal extensions to the graph map using curves isotopic to $T^U$, we can easily find properties of minimal extensions.

For each point $p$ of $T^P$, the boundary component obtained by cutting along $T^U(p)$ consists of two smooth pieces.
These pieces are homotopic to a curve in the graph called an \emph{unstable-parallel curve}.
\begin{definition}[Unstable-parallel curve]
A simple exact curve $\alpha:(I,J)\exto(\cut[U]{T},T^S)$ with endpoints in $T^S$ is \emph{unstable-parallel} if there is an exact isotopy $\alpha_t:(I,J)\exto(\cut[U]{T},T^S)$
 such that $\alpha_0=\alpha$ and $\alpha_1$ lies in $T^U$.
\end{definition}

\fig{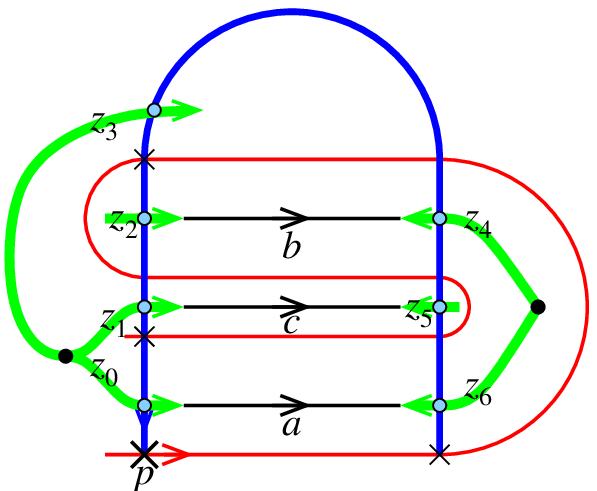}{Graph compatible with the Smale horseshoe trellis.}
In \figref{horseshoeembedded} we show a graph compatible with the Smale horseshoe trellis.
Passing anticlockwise around the unstable $T^U(p)$, we have two smooth curves, one curve $\alpha$ with the same orientation as $T^U(p)$ and lying to its right,
 and one curve $\beta$ with opposite orientation to its left.
These curves can be projected onto the graph representative $(G,W)$ giving edge-paths $\alpha= z_3\bar{z}_3z_1c\bar{z}_5z_5\bar{c}\bar{z}_1$
 and $\beta= z_0a\bar{z}_6z_4\bar{b}\bar{z}_1z_1b\bar{z}_4z_6\bar{a}\bar{z}_0$.
Note that each edge is traversed twice.
This is always the case if we traverse all boundary components of $\cut{T}$, including those which are boundary components of $M$.

\fig{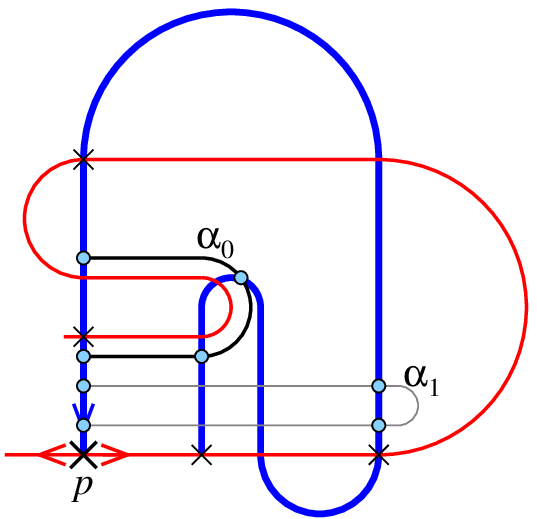}{Curve $\alpha_1$ is unstable isotopic to $T^U$, and $\alpha_1$ follows a minimal extension.}
One important feature about unstable-parallel curves is that if $([\widetilde{f}];\widetilde{T})$ is a minimal unstable extension of $([f];T)$
 and $\alpha$ is an unstable-parallel curve for $T$, then $\tilde{f}_{\min}(\alpha)$ is unstable-parallel curve for $\widetilde{T}$.
Therefore, we can use unstable parallel curves to construct minimal extensions, as shown in \figref{unstableisotopic}.
An unstable-parallel curve in a graph $G$ compatible with a trellis $T$ is a subpath of a peripheral loop.


\subsection{Reducibility of graph representatives}
\label{sec:graphreducibility}

Just as for trellis mapping classes, some graph maps are \emph{reducible} and can be split into simpler pieces.
\begin{definition}[Reducible graph map]
A controlled graph map $g$ of $(G,W)$ is \emph{reducible} if $G$ has an invariant subgraph $H$ such that $H$ does not contain any control edges and either
\begin{enumerate}\isz
  \item $H$ has negative Euler characteristic. We call this an \emph{$n$-component separating} reduction, where $n$ is the number of components of $H$.
  \item $H$ is a union of $n$ simple closed curves but is not a subset of $P$. We call this an \emph{$n$-curve non-separating} reduction.
\end{enumerate}
$g$ has an \emph{attractor-repellor} reduction if there is an invariant subgraph $H$ such that 
\begin{enumerate}
\setcounter{enumi}{2}
\item $H$ and $H^C$ both contain control edges.
\end{enumerate}
\end{definition}
If a controlled graph map compatible with an trellis mapping class has a reduction, we cut out the invariant subgraph and consider the resulting pieces separately.
The reduction algorithm is described in Section~\ref{sec:reductionalgorithms}.

If $[f;T]$ has a separating invariant curve reduction which gives rise to an invariant subsurface with no components of $T$, then there cannot be a graph representative,
 since we cannot move bad turns in the corresponding invariant subgraph to control edges.
Therefore, we can only hope to find graph representatives for a trellis type without invariant curve reductions.
The next result gives a condition for a trellis mapping class to have an attractor-repellor reduction.
\begin{theorem}[Reductions of compatible graph maps]
Let $g$ be a proper graph map compatible with a well-formed trellis type $[f;T]$.
If $g$ has a proper invariant subgraph $H$ which is disjoint from $W$ and either has negative Euler characteristic, or is a union of non-peripheral circles.
Then $[f;T]$ has an invariant curve reduction.
If $g$ has any other proper invariant subgraph which either has negative Euler characteristic or is a union of non-peripheral circles,
 then $g$ has an attractor-repellor decomposition.
\end{theorem}

\begin{proof}
If $H$ is in invariant subgraph which is disjoint from $W$, then we can take the boundary curves of $H$ to be a set of reducing curves.
If $H$ is any other invariant subgraph, let $\overline{H}=\bigcup_{n=0}^\infty g^{-n}(H)$ and let $\widehat{H}$ be the set of components of $H$ which
 either contain points of $H$ or contain control edges mapping into $H$.
Since $g(H)\subset H$, we must have $g(\widehat{H})\subset\widehat{H}$, and further, $g(\partial\widehat{H})\subset\partial\widehat{H}$.
Therefore $\partial\widehat{H}$ is contained in free edges. 
We introduce control edges at points of $\partial\widehat{H}$, which corresponds to introducing extra stable curves $T^S_N$ curves in $T$.
The region bounded by $T^U$ and $T^S_N$ is then an attractor, which we cut out to give one component of the reduction.
The complement of this region is a repellor. 
We can find the dynamics of this repellor by removing all control edges in $\widehat{H}$, which corresponds to removing all stable curves in $A$.
We cut along these curves to obtain the attractor.
\end{proof}

\begin{example}
\fig{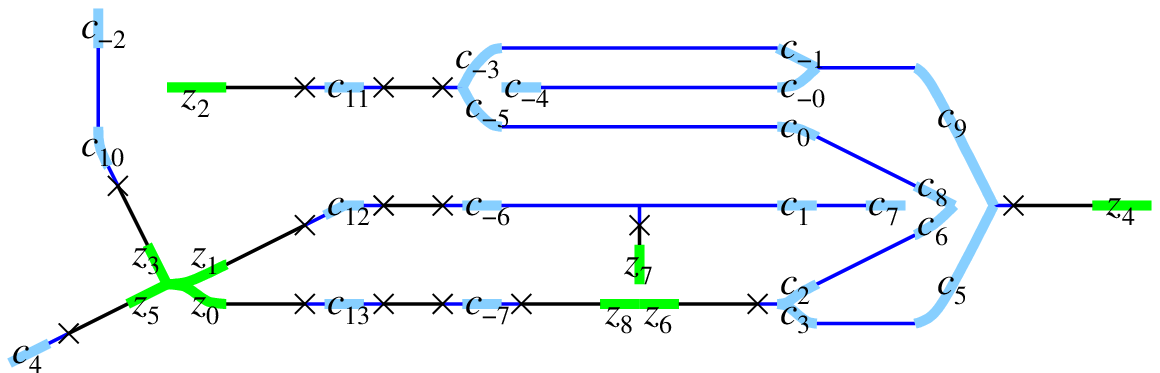}{Graph of a trellis with an attractor-repellor reduction.}
\figref{arreductiongraph} shows the graph of a trellis with an attractor-repellor reduction.
The control edges $z_i$ map
\[ 
  z_0,z_1,z_2\mapsto z_0; \quad z_3\mapsto z_1; \quad z_4\mapsto z_2; \quad z_5\mapsto z_3; \quad z_6,z_7,z_8\mapsto z_4;
\]
and the other control edges map
\[ \begin{array}{c}
     c_{-7},c_{-6}\mapsto c_{3}; \quad c_{-5},c_{-4},c_{-3}\mapsto c_{2}; \quad c_{-2}\mapsto c_{1};
       \quad c_{-1},c_{-0}\mapsto c_{0}; \\[1ex]
     c_{0},c_{1},c_{2}\mapsto c_{-0}; \quad c_{ 3}\mapsto c_{^-1}; \quad c_{4}\mapsto c_{-2}; \quad c_{5}\mapsto c_{-3}; 
       \quad c_{6},c_{7},c_{8}\mapsto c_{-4}; \quad c_{9}\mapsto c_{-5}; \\[1ex]
       \quad c_{10}\mapsto c_{-6}; \quad c_{11},c_{12},c_{13}\mapsto c_{-7}.
   \end{array}
\]
The boundary points of the invariant set $A$ are marked with crosses, and the control edges not in $A$ are those labelled $z_i$, and the control edges in $A$ are labelled $c_i$.
The trellis with this graph representative is the one given in \figref{arreduction}.
\end{example}

We now show that the every edge eventually maps over a peripheral edge or a control edge, or $g$ is reducible.
\begin{lemma}
\label{lem:expandingimage}
Let $(g;G,W)$ be the graph representative of an irreducible trellis type $[f;T]$.
Then if $e$ is any edge of $g$, there exists $n$ such that $g^n(e)$ contains a control edge or a peripheral edge.
\end{lemma}

\begin{proof}
Let $H=\bigcup_{n=0}^\infty g^n(e)$, an invariant subgraph of $G$.
Suppose $H$ contains no control edges or peripheral edges.
Then either $H$ is a forest, or $H$ must contain a component with negative (or zero) Euler characteristic, which means $g$ is reducible.
Both of these are contradictions.
\end{proof}


\subsection{Uniqueness of graph representatives}
\label{sec:graphuniqueness}

It is clear that given an efficient graph map compatible with a trellis map, the trellis mapping type can be reconstructed from the graph map.
The converse is also true; up to isomorphism there is only one efficient controlled graph map compatible with a given trellis mapping class.

\begin{theorem}[Uniqueness of graph representative]
\label{thm:efficientgraph}
Let $[f;T]$ be an irreducible trellis type.
Then, up to isomorphism, there is a unique efficient graph map $(g;G)$ compatible with $[f;T]$.
\end{theorem}

We first give a verifiable criterion for an exact homotopy equivalence between two graphs to be homotopic to a homeomorphism.
\begin{lemma}
\label{lem:graphhomeomorphism}
Suppose $(G_1,W_1)$ and $(G_2,W_2)$ are controlled graphs, and $p$ is an exact homotopy equivalence between $G_1$ and $G_2$.
Suppose there are strictly positive real functions $l_1$ and $l_2$ on the edges of $G_1$ and $G_2$ such that whenever $\alpha_1$ and $\alpha_2$ are tight paths joining points of $W$
 with $p\circ\alpha_1\homotopic\alpha_2$ we have $l_1(\alpha_1)=l_2(\alpha_2)$.
(We take the distance between a point of $W$ and a control vertex to be half the length of the control edge.)
Then $p$ is homotopic to a homeomorphism from $(G_1,W_1)$ to $(G_2,W_2)$.
\end{lemma}

\begin{proof}
It is sufficient to show that the universal covers $(\widetilde{G}_1,\widetilde{W}_1)$ and  $(\widetilde{G}_2,\widetilde{W}_2)$ are homeomorphic.
Since $W_1$ and $W_2$ are homeomorphic, and $G_1$ and $G_2$ are homotopy-equivalent, the sets $\widetilde{W}_1$ and $\widetilde{W}_2$ are homeomorphic, and we call them both $\widetilde{W}$.
The functions $l_1$ and $l_2$ lift naturally to $\widetilde{G}_1$ and $\widetilde{G}_2$.
For any four points $\widetilde{w}_1$, $\widetilde{w}_2$, $\widetilde{w}_3$ and $\widetilde{w}_4$ in $\widetilde{W}$. there are only two possibilities for the span; either it contains a valence-$4$ vertex or two valence-$3$ vertices.
Let $\alpha_{ij}$ be the curve from $w_i$ to $w_j$.
\fig{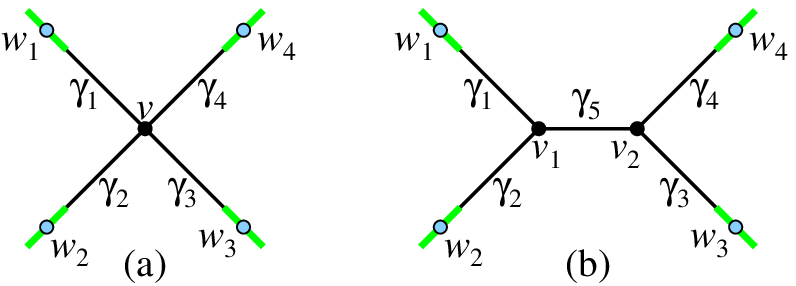}{The span of four points of $W$ in a simply connected region can either have a valence-$4$ vertex (a) or two valence-$3$ vertices (b).}

Suppose the span of $\{w_1,w_2,w_3,w_4\}$ contains a valence-$4$ vertex $v$ as shown in \figrefpart{grapheqn}{a}
Then for $i=1,2$,  
\[ l_i(\alpha_{12})+l(\alpha_{34})=l_i(\alpha_{13})+l_i(\alpha_{24})=l_i(\alpha_{14})+l_i(\alpha_{23})=l_i(\gamma_{1})+l_i(\gamma_{2})+l_i(\gamma_{3})+l_i(\gamma_{4}) \]
However, if the span of $G\{w_1,w_2,w_3,w_4\}$ contains a valence-$3$ vertices $v_1$ and $v_2$ as shown in \figrefpart{grapheqn}{b}, we have
\[ l_i(\alpha_{12})+l(\alpha_{34})=l_i(\gamma_{1})+l_i(\gamma_{2})+l_i(\gamma_{3})+l_i(\gamma_{4}) \]
but
\[ l_i(\alpha_{13})+l_i(\alpha_{24})=l_i(\alpha_{14})+l_i(\alpha_{23})=l_i(\gamma_{1})+l_i(\gamma_{2})+l_i(\gamma_{3})+l_i(\gamma_{4})+2l_i(\gamma_{5}), \]
so 
\[ 
  l_i(\alpha_{12})+l(\alpha_{34})\leq l_i(\alpha_{13})+l_i(\alpha_{24})=l_i(\alpha_{14})+l_i(\alpha_{23})
\]
Hence the universal covers are homeomorphic by a homeomorphism $\widetilde{h}$ homotopic to $\widetilde{p}$.
Further, since the computations given above are equivariant, the homeomorphism $\widetilde{h}$ is equivariant, so projects to a homeomorphism $h:(G_1,W_1)\fto(G_2,W_2)$.
\end{proof}

The proof of the main theorem is as follows:
\begin{proof}
Let $p:(g_1;G_1)\fto (g_2;G_2)$ be an exact homotopy equivalence with homotopy inverse $q$.
Thus $q\circ p\homotopic\id$, $p\circ q\homotopic\id$ and $g_1\homotopic q\circ g_2\circ p$, taking homotopies relative to $W$.
It is clear that we can choose $p$ and $q$ to be mutual inverses on the control and peripheral edges, so merely need to consider the expanding edges.
We therefore make no distinctions between control edges in $G$ and $G_2$.

First, consider a tight edge-path $\gamma_1$ in $G_1$ which has endpoints in $W$ but no interior points in $W$.
Let $\gamma_2$ be the tight edge-path homotopic to $p(\gamma_1)$.
Consider $\gamma_1$ and $\gamma_2$ as exact curves $(I,J_0)\exto(G,W)$, and take minimal iterates under $g_1$ and $g_2$.
Since $G_1$ and $G_2$ are graphs, these minimal iterates can each be represented by a single tight edge-path.
Then $\gamma_1\homotopic q\circ p(\gamma_1)\homotopic q(\gamma_2)$, so $p({(g_1)}_{\min}(\gamma_1))\homotopic{(g_2)}_{\min}(\gamma_2)$.
Taking further iterates gives $p({(g_1)}^n_{\min}(\gamma_1))\homotopic{(g_2)}^n_{\min}(\gamma_2)$ for all $n$.
Thus ${(g_1)}^n_{\min}(\gamma_1)$ and ${(g_2)}^n_{\min}(\gamma_2)$ contain the same number of control and peripheral edges for all $n$.
Since by Lemma~\ref{lem:expandingimage}, this number must be positive for all edges for some $n$, the graphs $G_1$ and $G_2$ are homeomorphic by a graph map $h$ homotopic to $p$ by Lemma~\ref{lem:graphhomeomorphism}.
This homeomorphism gives a topological conjugacy between $g_1$ and $g_2$.
\end{proof}


\section{Computing a Graph Representative}
\label{sec:algorithm}

We now give an algorithm for obtaining an optimal controlled graph map compatible with a given trellis mapping class.
If the trellis mapping class has is reducible, the algorithm may instead find a reduction, but can then be used to compute a graph representative for the irreducible components.
The moves of the algorithm are based on the algorithm of Bestvina and Handel \cite{BestvinaHandel95}, but the method is quite different.
Instead of reducing the topological entropy at each step, the algorithm reduces a zeta function, which gives a more precise measure of the growth based on the control and peripheral edges.
This results in an algorithm where we need only consider the local action of the graph map.

In the case of a surface mapping class, there are no control edges, and no optimal graph map.
We use a slightly different algorithm to find an efficient graph map, and consider a zeta function based on the peripheral edges to show that the algorithm terminates.
If there are no peripheral edges, we create temporary peripheral edges by puncturing at a periodic orbit.


\subsection{Graph moves}
\label{sec:graphmoves}

We now describe the moves we need for the algorithm.
Each move is an exact homotopy equivalence; in particular, the Euler characteristic and number of control edges are unchanged.
The property of being surface embeddable need not be preserved in general;
 however, performing any move followed by all possible vertex and edge tightenings does preserve this property.

The movers of edge tightening and vertex homotopy do not affect the graph $(G,W)$, but change the map $g$ to a new map $g\p$.
\begin{description}
\item[Edge tightening]
  Suppose $e_0$ be an edge of $G$ with $g(e_0)=e_1\ldots e_ie_{i+1}\ldots e_n$ and $e_{i+1}=\bar{e}_i$.
  Homotope $g$ to $g\p$ with $g\p(e_0)=e_1\ldots e_{i-1}e_{i+2}e_n$.
  \fig{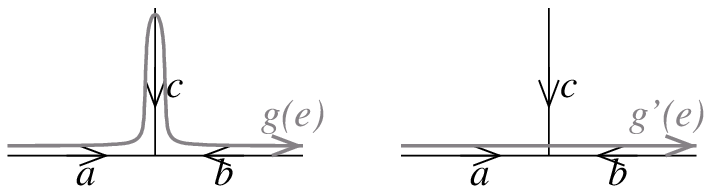}{Tightening the edge $e$.}
\item[Vertex homotopy]
  Let $v$ be a vertex of $G$ with incident edges $e_1,\ldots,e_k$, and let $\beta$ be an edge path with $\init{\beta}=g(v)$.
  Homotope $g$ to a new map $g\p$ such that $g\p(e_i)=\bar{\beta}g(e_i)$ for $i=1,\ldots,e_k$, unless $g(e_i)$ can be written $\beta\epsilon_i$ for some edge-path $\epsilon_i$,
    in which case we take $g\p(e_i)=\epsilon_i$.
  We say this move homotopes $v$ \emph{across} the edge-path $\beta$.
\end{description}

The moves of edge collapsing and vertex splitting change a graph $G$ to a graph $G\p$.
\begin{description}
\item[Edge collapsing]
  Suppose $e_0$ is an edge of $G$ with initial vertex $v$ and final vertex $w$ such that $g(e_0)=\cdot$.
  Suppose the cyclic ordering of edges at $v$ is $\lhd e_0\lhd a_1\lhd \cdots\lhd a_j\lhd e_0\lhd$ and the cyclic ordering of edges at $w$ is $\lhd\bar{e}_0\lhd b_1\lhd \cdots b_k\lhd\bar{e}_0\lhd$.
  Let $G\p$ be the graph formed by removing the vertex $w$ and the edge $e_0$ from $G$ such that the cyclic ordering of edges at $v$ is $\lhd a_1\lhd \cdots a_j\lhd b_1\lhd\cdots b_k\lhd a_1\lhd$.
  Let $g\p$ be the graph map obtained by removing all occurrences of $e_0$ from the image $g(e)$ of any edge $e$.
\item[Vertex splitting]
  Suppose $v$ is a vertex of $G$, and the cyclic ordering of edges at $v$ is $\lhd a_1\lhd\cdots\lhd a_j\lhd b_1\lhd\cdots\lhd b_k\lhd$.
  Let $G\p$ be a graph formed by introducing a new vertex $w$ and new edge $e_0$ to $G$ such that the cyclic ordering of edges at $v$ is $\lhd e_0\lhd a_1\lhd \cdots\lhd a_j\lhd e_0\lhd$ and the cyclic ordering of edges at $w$ is $\lhd\bar{e}_0\lhd b_1\lhd \cdots b_k\lhd\bar{e}_0\lhd$.
  Let $g\p$ be the graph map obtained by taking $g\p(e_0)=\cdot$, and by replacing every occurrence of an edge $b_i$ for $1\leq i\leq k$ in the image $g(e)$ with the edge-path $e_0b_j$ (and every occurrence of $\bar{b}_i$ with $\bar{b}_i\bar{e}_0$).
  \fig{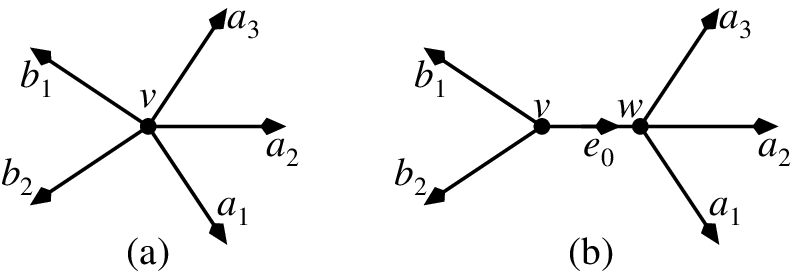}{(a) Before and (b) after a vertex splitting.}
\end{description}
Edge collapsing and vertex splitting are inverse operations.

The following moves can all be expressed in terms of the basic moves described above.
\begin{description}
\item[Vertex tightening]
  If $v$ is a vertex of $G$ with incident edges $e_1,\ldots,e_k$ and there is an edge-path $\beta$ such that $g(e_i)=\beta\epsilon_i$ for all edges $e_i$ with initial vertex $v$, 
  homotope $v$ across the edge-path $\beta$.
  The resulting graph map $g\p$ has $g\p(e_i)=\epsilon_i$ for all edges $e_i$ with initial vertex $v$.
\item[Collapsing an invariant forest]
  If $H\subset G$ is a set of expanding edges with $g(H)\subset H$ such that every component of $H$ is simply connected (i.e. $H$ is a forest),
   collapse every edge of $T$ to a point.
  This move can be realised as a combination of homotopies and edge collapsing (though at intermediate stages, the set of control edges may not be invariant).
  The move is useful since the topology and dynamics of an invariant forest are trivial, and are removed to preserve transitivity of the graph map.
\item[Tidying]
  Perform edge tightenings, vertex tightenings, and collapsing invariant forests until no further tightening or collapsing is possible,
    with the exception that we never collapse peripheral edges at this stage (even if they have a trivial image).
\item[Valence-$3$ homotopy]
  Suppose $v$ is a valence-$3$ vertex with incident edges $e_1$, $e_2$ and $e_3$ such that $g(e_i)=\beta\epsilon_i$ for $i=1,2$ and $g(e_3)=\epsilon_3$.
  Homotope $v$ across $\beta$ and to give $g\p(e_i)=\epsilon_i$ for $i=1,2$ and $g\p(e_3)=\bar{\beta}\epsilon_3$.
\item[Folding]
  Suppose there are edges $e_1\lhd e_2$ at a vertex $v$ with $g(e_1)=\epsilon_0\epsilon_1$ and $g(e_2)=\epsilon_0\epsilon_2$.
  First split off $e_1$ and $e_2$ to a new vertex $w$ with new edge $e_0$ from $v$ to $w$.  
  We then perform any possible edge tightenings; note that since $e_0$ maps to $\cdot$, we cannot tighten at $w$.
  This gives a graph map $g\p$ with $g\p(e_0)=\cdot$, $g\p(e_1)=\epsilon_0\p\epsilon_1\p$ and  $g\p(e_2)=\epsilon_0\p\epsilon_2\p$.
  Finally, we homotope at $w$ across the edge-path $\epsilon_0\p$ to obtain a graph map $g\pp$ with $g\pp(e_0)=\epsilon_0$, $g\pp(e_1)=\epsilon_1$ and $g\pp(e_2)=\epsilon_2$.
\end{description}
The analysis of tightening moves is simple, since they always decrease the complexity of the graph map.
Similarly, collapsing invariant forests also decreases the complexity.
Performing the tidying operation leave a nice graph map to use as a starting point for the more complicated moves of vertex homotopy and edge splitting.
\begin{definition}[Tight and tidy graph maps].
We say a graph map is \emph{tight} if no edge or vertex tightenings are possible.
A tight graph map is \emph{tidy} if it has no invariant forest of expanding edges, and every valence-$2$ vertex is an endpoint of a control edge.
\end{definition}

\fig{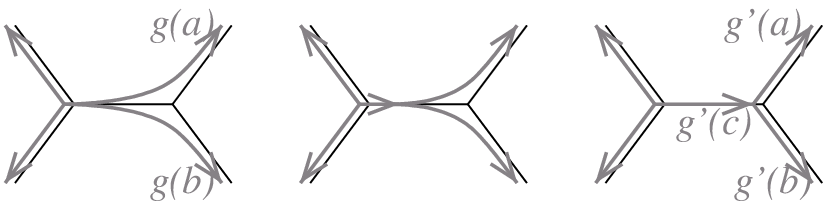}{Folding edges $a$ and $b$.}


\subsection{Effect of graph moves}
\label{sec:moveeffect}

We measure the complexity of a graph map by its \emph{zeta functions}.
\begin{definition}[Zeta function]
Let $g$ be a graph map of a graph $G$ and $H$ an invariant subgraph.
The \emph{$H$-length} of an edge path $\epsilon$, denoted $\len[H]{\epsilon}$ is the number of occurrences of an edge of $H$ in $\epsilon$.
The \emph{$H$-norm} of a graph map $g$ is given by 
\[
  \norm[H]{g} = \sum_{e\in \edge(G)}\len[H]{g(e)}
\]
and the \emph{control zeta function} of $g$ is 
\[
  \zeta_{g;H}(t) = \sum_{n=0}^\infty \norm[H]{g^n}\,t^n .
\]
\end{definition}
It is straightforward to show that if every invariant subgraph of $G$ contains an edge in $H$, then the logarithm of the asymptotic growth rate of $\norm[H]{g^n}$ equals $\htop(g)$.
If $R$ is the radius of convergence of $\zeta_{g;H}$, the also $\htop(g)=\log(1/R)$, where $R$ is

Zeta functions are ordered by the ordering $\sum_{n=0}^{\infty} a_n\,t^n<\sum_{n=0}^{\infty}b_n\,t^n$ if for some $k$ we have $a_k<b_k$, but $a_i=b_i$ for $i<k$.
If $k\leq N$, we say $\sum_{n=0}^{\infty} a_n\,t^n<_N\sum_{n=0}^{\infty}b_n\,t^n$.
Equivalently, $\zeta_1<\zeta_2$ if there exists $\epsilon>0$ such that $\zeta_1(t)<\zeta_2(t)$ for all $t\in(0,\epsilon)$, and $\zeta_1<N\zeta_2$ if $\zeta_1+at^n<\zeta_2$ for all $a\geq 0$ and $n>N$.

The following lemma gives some simple results on $H$-norms. 
\begin{lemma}
Let $g$ be a graph map and $H$ an invariant subgraph of $G$. Then
\begin{enumerate}
\item If $g\p$ is obtained from $g$ by edge tightening or vertex tightening, then $\norm{(g\p)^n}\leq \norm[H]{g^n}$.
\item If $g\p$ is obtained from $g$ by collapsing an edge $e$, then $\norm[H]{(g\p)^n}<\norm[H]{g^n}$ if $e\in H$ and $\norm[H]{(g\p)^n}=\norm[H]{g^n}$ otherwise.
\item If $g\p$ is obtained from $g$ by collapsing an invariant forest $F$, then $\norm[H]{(g\p)^n}<\norm[H]{g^n}$ if $F\cap H\neq\emptyset$, and $\norm[H]{(g\p)^n}=\norm[H]{g^n}$ otherwise.
\item If $g\p$ is obtained from $g$ by a vertex splitting introducing an edge $e$ to give a new graph $G\p$, then $\norm[H]{(g\p)^n}=\norm[H]{g^n}$.
\end{enumerate}
\end{lemma}

Valence-$3$ homotopies are more complicated, but the following result shows that $H$-norms are still reduced.
\begin{lemma}
\label{lem:valencethreehomotopy}
Let $g$ be a graph map with a valence-$3$ vertex $v$ which has incident edges $e_1$, $e_2$ and $e_3$ such that $g(e_i)=\beta\epsilon_i$ for $i=1,2$ and $g(e_3)=\epsilon_3$, and that $g\p$ is obtained from $g$ by homotoping $v$ across $\beta$.
Suppose further that $g^{n-2}(\beta)$ does not contain any edges of $H$, but $g^{n-1}(\beta)$ does.
Then
\begin{enumerate}\isz
\item\label{item:edgehnorm} For any $k< n$ and any edge $e$ of $G$, $\len[H]{{(g\p)}^{k}(e)}=\len[H]{{g}^{k}(e)}$, and $\norm[H]{{(g\p)}^k}=\norm[H]{g^k}$.
\item\label{item:graphhnorm} $\norm[H]{{(g\p)}^{n}}<\norm[H]{g^{n}}$
\end{enumerate}
\end{lemma}

\begin{proof}
Since $g^k(\beta)$ does not contain any control edges for $k<n-1$, so $\norm[H]{g^k(\beta)}=0$ for these values of $k$.
For part~\ref{item:edgehnorm}, we use induction on $k$.
The result is trivial for $k=0$.
Assume the result is true for some $k<n$.
Then 
\begin{eqnarray*}
  \norm[H]{{(g\p)}^{k+1}(e_i)} 
    & = & \norm[H]{{(g\p)}^k(\epsilon_i)} = \norm[H]{g^k(\epsilon_i)} \textrm{ (by the inductive hypothesis)} \\
    & = & \norm[H]{g^k(\beta)g^k(\epsilon_i)}-\norm[H]{g^k(\beta)} = \norm[H]{g^k(\beta\epsilon_i)}-\norm[H]{g^k(\beta)} \\
    & = & \norm[H]{g^{k+1}(e_i)}-\norm[H]{g^k(\beta)} \mbox{ for } i=1,2 \\ 
 \norm[H]{{(g\p)}^{k+1}(e_3)}
    & = & \norm[H]{{(g\p)}^k(\beta\epsilon_3)} = \norm[H]{g^k(\beta\epsilon_3)} = \norm[H]{g^k(\beta)}+\norm[H]{g^k(\epsilon_3)} \\
    & = & \norm[H]{g^k(e_3)}+\norm[H]{g^k(\beta)} \\
 \norm[H]{{(g\p)}^{k+1}(e_i)} & = & \norm[H]{{(g\p)}^k(\epsilon_i)} = \norm[H]{g^k(\epsilon_i)} = \norm[H]{g^{k+1}(e_i)} \mbox{ otherwise} .
\end{eqnarray*} 
Clearly then, $\norm[H]{{(g\p)}^k} = \sum_{e\in G}\norm[H]{{(g\p)}^k(e)} = \sum_{e\in G}\norm[H]{{(g\p)}^k(e)} = \norm[H]{g^k}$ for $k<n$.
Part~\ref{item:graphhnorm} follows from the above calculations since $\norm[H]{{(g\p)}^{n}}=\norm[H]{g^{n}}-\norm[H]{g^{n-1}(\beta)}<\norm[H]{g^n}$.
\end{proof}

The following lemma shows that if we can find a bound on the number of edges of a graph, we can control which $H$-norm of $g$ decreases in a valence-$3$ homotopy.
\begin{lemma}\label{lem:irreduciblehomotopy}
Suppose $g$ is a map of a graph $G$ with invariant subgraph $H$, and there are no invariant subgraphs of $g$ which do not contain an edge in $H$.
Suppose $g\p$ is obtained from $g$ by a valence-$3$ homotopy.
Then there exists $n\leq |G|-|H|$ such that $\norm[H]{{(g\p)}^k} =\norm[H]{{g}^k}$ for $k< n$ and $\norm[H]{{(g\p)}^{n}} = \norm[H]{{g}^{n}}$
In particular, $\zeta_{g\p;H}<_{|G|} \zeta_{g;H}$.
\end{lemma}

\begin{proof}
Let $K_0$ be the subgraph containing the edges of $\beta$, and define $K_i$ by $K_i=K_{i-1}\cup g(K_{i-1})$.
If $K_{n}$ contains an edge of $H$ for some least $n$, then $g^n(\beta)$ contains an edge of $H$ and we are done.
If not, then $K_n$ has more edges than $K_{n-1}$, and hence has at least $n+1$ edges.
Therefore, $K_n$ contains an edge of $H$ for some $n\leq |G|-|H|$, the number of edges in $G\setminus H$.
\end{proof}

The following lemma is used to show that our algorithm terminates.
\begin{lemma}\label{lem:finiteprocess}
Suppose an algorithm consists of graph moves such that there exists $N\in\N$ such that at each stage $\zeta_{g\p;H}<_N \zeta_{g;H}$.
Then the algorithm terminates.
\end{lemma}

\begin{proof}
We use induction on $N$.
Clearly, the result is true for $N=0$, since $\norm[H]{g^0}=|H|$ is a non-negative integer.

Suppose the result holds for $N-1$.
Any move either decreases $\norm[H]{g^N}$ or $\norm[H]{g^k}$ for $k\leq N-1$, but since $\norm[H]{g^N}$ is a non-negative integer,
 there is no infinite sequence of moves only decreasing $\norm[H]{g^N}$. 
By the inductive hypothesis, there can only be finitely many moves for which $\norm[H]{g^k}$ with $k<N$ is decreased, and the result follows.
\end{proof}


\subsection{Algorithm to compute the graph representative}
\label{sec:optimalalgorithms}

We now give the algorithm to compute a graph representative, or find a reduction.
An algorithm of Bestvina and Handel for free groups \cite{BestvinaHandel92} considers the reducible case as a single graph.

\begin{algorithm}
\label{alg:main}
Suppose $g$ is a controlled graph map.
Perform the following moves, followed by tidying (noting that tidying does not affect the peripheral subgraph), until the graph map is optimal or there is an invariant subgraph of free edges which is not a forest.
\begin{enumerate}\isz
\item\label{item:algfold} Fold a bad turn at a vertex of valence greater than $3$, or at a control vertex.
\item\label{item:alghomotope}  Homotope at a bad turn at a free vertex of valence $3$.
\item\label{item:algcollapse}  If no other moves are possible, collapse any peripheral edges which have trivial image.
\end{enumerate}
\end{algorithm}

The proof of termination considers zeta functions of the graph map.
We need to consider both the control zeta function $\zeta_{g;C}$ and the peripheral zeta function $\zeta_{g;C}$, and also the number of peripheral edges $|P|$.
The standard procedure for dealing with the peripheral subgraph $P$ is \emph{absorbing},
 which yields a graph map for which there is no invariant set which deformation-retracts onto $P$ other than $P$ itself.
However, absorbing into the peripheral subgraph can be accomplished by the moves of Algorithm~\ref{alg:main}.
A fundamental observation is that $G$ has at most $|G|=2\card{C}-3\chi$ edges, where $\chi$ is the Euler characteristic, and hence any decrease in a zeta function must occur in at most $|G|$ coefficients.

\begin{theorem}[Computation of the graph representative]
\label{thm:termination}
Algorithm \ref{alg:main} terminates at a reduction or an optimal graph map.
\end{theorem}

\begin{proof}
If there is an invariant subgraph which does not contain any control or peripheral edges and is not a forest, $g$ is irreducible and the algorithm has found a reduction.
Further, at any stage there are at most $|E|=2|C|-3\chi$ edges of the graph, which means that any reduction of a zeta function must occur at most in the $|E|$th coefficient.
Therefore, every move apart from collapsing in the peripheral subgraph either reduces the control zeta function, increases the number of peripheral edges, or decreases the peripheral zeta function.

After collapsing peripheral edges in step \ref{item:algcollapse}, the link of the peripheral subgraph must be invariant under $\partial g$.
If there is no reduction, every expanding edge in this subgraph must eventually map to a control edge, or else it would be contained in an invariant subgraph of free edges with negative Euler characteristic.
We may perform homotopies across pre-peripheral edges disjoint from $P$ which do not affect the control zeta function, but cannot introduce new peripheral edges without homotoping across an edge in $\link{P}$
At this homotopy, the control zeta function decreases, so there is no infinite sequence of moves which does not decrease the control zeta function by Lemma~\ref{lem:finiteprocess}.
A further application of Lemma~\ref{lem:finiteprocess} shows that there are at most finitely many moves which decrease the control zeta function.
Hence the algorithm either finds a reduction, or terminates at an optimal graph map.
\end{proof}

If there are no control edges, we are in the realm of Nielsen-Thurston theory.
This case cannot be treated directly by our algorithm, since there is no optimal graph map in the homotopy class.
The proof of Theorem~\ref{thm:termination} fails since after collapsing peripheral edges, all edges in $\link{P}$ map to cover the entire graph, which is itself an invariant subgraph of negative Euler characteristic which does not include any control edges.
Instead we consider the entropy in addition to the peripheral zeta function in order to show that the algorithm finds an efficient graph representative.

\begin{algorithm}\label{alg:peripheral}\ 
Suppose $g$ is a graph map compatible with a surface mapping class, and $g$ has nonempty peripheral subgraph $P$.
\begin{enumerate}\isz
\item\label{item:algmain} Perform Algorithm~\ref{alg:main}, stopping after collapsing peripheral edges in step~\ref{item:algcollapse}.
\item If the graph map is efficient or has a reduction, the algorithm terminates. Otherwise, repeat step~\ref{item:algmain}.
\end{enumerate}
Perform Algorithm~\ref{alg:main}, checking after step~\ref{item:algcollapse} whether the graph map is efficient, at which point the algorithm terminates.
\end{algorithm}
If there are no peripheral edges, we need to artificially introduce peripheral edges by puncturing at a periodic orbit.
\begin{algorithm}\label{alg:noperipheral}\ 
Suppose $g$ is a map of a graph $G$ with no peripheral curves.
Perform the following moves until an efficient graph map or a reduction is found.
\begin{enumerate}\isz
\item\label{item:algpuncture} Let $P$ be a periodic orbit of $g$. Puncture at points of $P$ by introducing artificial peripheral loops at $P$.
\item Perform the Algorithm~\ref{alg:peripheral} for graph maps with peripheral loops.
\item Collapse the artificial peripheral loops at $P$ and tidy. If the resulting graph map is efficient or has a reduction, the algorithm terminates.
  If the resulting graph map has peripheral loops, perform Algorithm~\ref{alg:peripheral}.
  Otherwise, return to step~\ref{item:algpuncture}
\end{enumerate}
\end{algorithm}

\begin{theorem}
Algorithms~\ref{alg:peripheral} and \ref{alg:noperipheral} terminate at a reduction or an efficient graph map.
\end{theorem}

\begin{proof}
We first consider Algorithm~\ref{alg:peripheral}.
If $g$ is not efficient, there is a turn in the link of $P$ which is inefficient.
Folding this turn decreases entropy.
Although other moves may increase entropy, the graph map obtained by folding and homotoping can also be obtained by the moves of the Bestvina-Handel algorithm, so the entropy must strictly decrease.
Since the entropy is bounded above, and the number of edges of $g$ is bounded by $-3\chi$, there are only finitely many possible values for the entropy.
Hence Algorithm~\ref{alg:peripheral} terminates.

Algorithm~\ref{alg:noperipheral} terminates for similar reasons, since the entropy is bounded, and must strictly decrease.
\end{proof}

We now give an example of computing an optimal graph map using Algorithm~\ref{alg:main}.
\begin{example}
\label{ex:graphalgorithm}
Consider the trellis type $[f;T]$ shown in \figref{algorithm2trellis}.
\fig{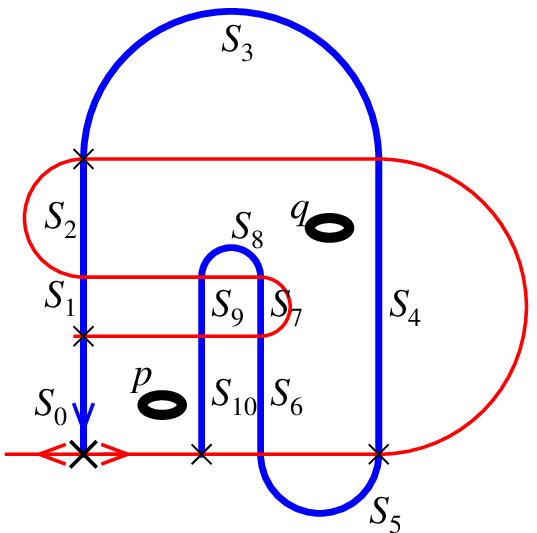}{A trellis with a period-$2$ orbit.}

There are eleven stable segments, and the surface has two punctures, $p$ and $q$, which are permuted by the diffeomorphism $f$.

\fig{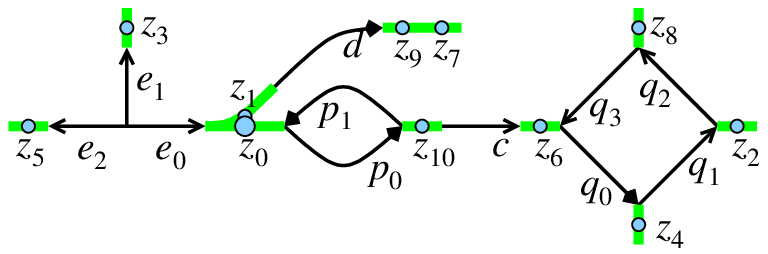}{An initial compatible graph}
An initial compatible graph is shown in \figref{algorithm2initial}.
The edge-loop $p_0p_1$ surrounds the puncture $p$, and the edge-loop $q_0q_1q_2q_3$ surrounds the puncture $q$.
The control edges map under $g$ as follows:
\[ z_{0},z_{1},z_{2}\mapsto z_{0}; \quad z_{3}\mapsto z_{1}; \quad z_{4}\mapsto z_{2}; \quad z_{5}\mapsto z_{3}; \quad z_{6},z_{7},z_{8},z_{9},z_{10}\mapsto z_{4} \]
We compute $g$ on the expanding edges, which map:
\[ \begin{array}{c}
     p_0\mapsto p_0\bar{z}_{10}cz_{6}q_0; \quad p_1\mapsto q_1q_2q_3\bar{z}_{6}\bar{c}z_{10}\bar{p}_0; \\[1ex]
     q_0\mapsto q_1; \quad q_1\mapsto q_2q_3\bar{z}_{6}\bar{c}z_{10}p_1; \quad q_2\mapsto p_0\bar{z}_{10}cz_{6}\bar{q}_3\bar{q}_2\bar{q}_1; \quad q_3\mapsto \cdot; \\[1ex]
     c\mapsto \cdot; \quad d\mapsto p_0\bar{z}_{10}cz_{6}\bar{q}_3\bar{q}_2\bar{q}_1; \quad e_0\mapsto e_0; \quad e_1\mapsto e_0; \quad e_2\mapsto e_1. 
   \end{array}
\]
Note that the edges $p_i$ and $q_i$ are not peripheral edges, since they do not form an invariant set. 
However the edge-loop $p_0p_1$ maps to the edge loop $p_0\bar{z}_{10}cz_{6}q_0q_1q_2q_3\bar{z}_{6}\bar{c}z_{10}\bar{p}_0$ which is freely homotopic to $q_0q_1q_2q_3$.

\fig{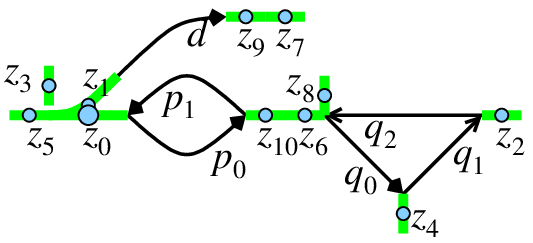}{Graph obtained by collapsing the invariant forest.}
The edges $q_3$, $c$, $e_1$, $e_2$ and $e_3$ form an invariant forest which does not contain any control edges.
Each of these edges can be collapsed to a point to give the graph in \figref{algorithm2conditioned}.
The expanding edges now map
\[ \begin{array}{c}
     p_0\mapsto p_0\bar{z}_{10}z_{6}q_0; \quad p_1\mapsto q_1q_2\bar{z}_{6}z_{10}\bar{p}_0; \\[1ex]
     q_0\mapsto q_1; \quad q_1\mapsto q_2\bar{z}_{6}z_{10}p_1; \quad q_2,d\mapsto p_0\bar{z}_{10}z_{6}\bar{q}_2\bar{q}_1.
   \end{array}
\]
  
\fig{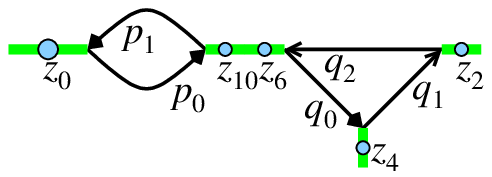}{Essential graph map.}
The edge $d$ has no preimage, so is not part of the essential graph representative.
We can remove $d$, and also $z_1$, $z_3$, $z_5$ and $z_8$ to obtain the graph map shown in \figref{algorithm2essential}
The control zeta function is $\zeta_{g;C}(t)=5+13t+25t^2+\cdots$.

\fig{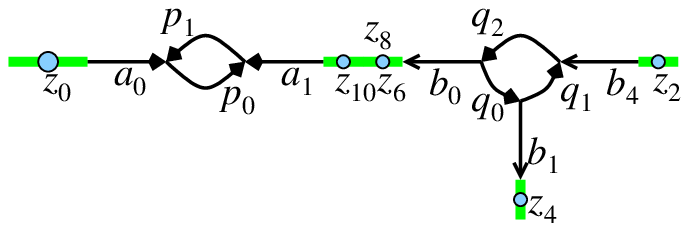}{Graph obtained after folding}
Performing five folding operations introduced edges $a_0$, $a_1$, $b_0$, $b_1$ and $b_4$, as shown in \figref{algorithm2folded}.
Before tightening, the graph map is
\[ \begin{array}{c}
     p_0\mapsto a_0p_0\bar{a}_1\bar{z}_{10}z_{6}\bar{b}_0q_0b_1; \quad p_1\mapsto \bar{b}_1q_1\bar{b}_4b_4q_2b_0\bar{z}_{6}z_{10}a_1\bar{p}_0\bar{a}_0; \\[1ex]
     q_0\mapsto \bar{b}_1q_1\bar{b}_4; \quad q_1\mapsto b_4q_2b_0\bar{z}_{6}z_{10}a_1p_1\bar{a}_0;
       \quad q_2\mapsto a_0p_0\bar{a}_1\bar{z}_{10}z_{6}\bar{b}_0\bar{q}_2\bar{b}_4b_4\bar{q}_1b_1.
   \end{array}
\]

On performing tightenings and valence-$3$ homotopies, we obtain the graph map
\[ \begin{array}{c} 
     p_0\mapsto q_0; \quad p_1\mapsto q_1q_2; \quad q_0\mapsto \cdot; \quad
       q_1\mapsto q_2b_0\bar{z}_{6}z_{10}a_1p_1; \quad q_2\mapsto p_0\bar{a}_1\bar{z}_{10}z_{6}\bar{b}_0\bar{q}_2; \\[1ex]
     a_0\mapsto a_0p_0\bar{a}_1\bar{z}_{10}z_{6}\bar{b}_0; \quad a_1\mapsto \bar{b}_1; \quad b_0\mapsto \bar{q}_1b_1; \quad b_1\mapsto \bar{b}_4; \quad b_4\mapsto a_0. 
    \end{array}
\]
 
\fig{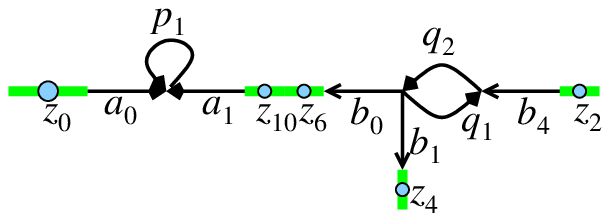}{Graph obtained after collapsing $q_0$ and $p_0$.}
Collapsing $q_0$ and then $p_0$ gives the graph map shown in \figref{algorithm2collapsed}, with
\[ \begin{array}{c} 
     \quad p_1\mapsto q_1q_2; \quad q_1\mapsto q_2b_0\bar{z}_{6}z_{10}a_1p_1; \quad q_2\mapsto \bar{a}_1\bar{z}_{10}z_{6}\bar{b}_0\bar{q}_2; \\[1ex]
     a_0\mapsto a_0\bar{a}_1\bar{z}_{10}z_{6}\bar{b}_0; \quad a_1\mapsto \bar{b}_1; \quad b_0\mapsto \bar{q}_1b_1; \quad b_1\mapsto \bar{b}_4; \quad b_4\mapsto a_0. 
    \end{array}
\]
and control zeta function $\zeta_{C}(t)=5+11t+19t^2+\cdots$.

\fig{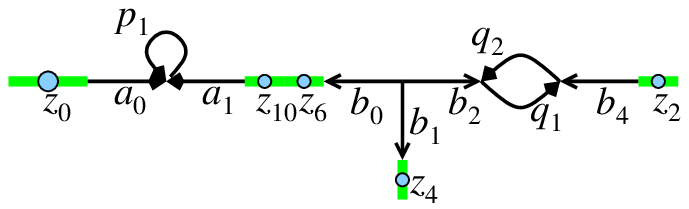}{Graph obtained after folding $q_1$ and $q_2$.}
Folding $q_1$ and $\bar{q}_2$ to form a new edge $b_2$ gives the graph map shown in \figref{algorithm2refolded}, with
\[ \begin{array}{c} 
     \quad p_1\mapsto b_2q_1q_2\bar{b}_2; \quad q_1\mapsto p_1; \quad q_2\mapsto \cdot; \\[1ex]
     a_0\mapsto a_0\bar{a}_1\bar{z}_{10}z_{6}\bar{b}_0; \quad a_1\mapsto \bar{b}_1; \quad 
     b_0\mapsto \bar{q}_1\bar{b}_2b_1; \quad b_1\mapsto \bar{b}_4; \quad b_2\mapsto q_2\bar{b}_2b_0\bar{z}_{6}z_{10}a_1; \quad b_4\mapsto a_0. 
    \end{array}
\]

\fig{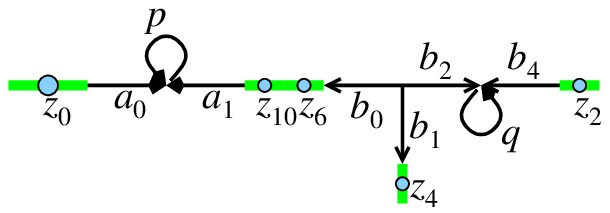}{Graph obtained after collapsing $q_2$.}
Since now $q_2$ has trivial image, we can collapse $q_2$ to obtain the graph map shown in \figref{algorithm2recollapsed}, with
\[ \begin{array}{c} 
     \quad p_1\mapsto b_2q_1\bar{b}_2; \quad q_1\mapsto p_1; \\[1ex]
     a_0\mapsto a_0\bar{a}_1\bar{z}_{10}z_{6}\bar{b}_0; \quad a_1\mapsto \bar{b}_1; \quad 
     b_0\mapsto \bar{q}_1\bar{b}_2b_1; \quad b_1\mapsto \bar{b}_4; \quad b_2\mapsto \bar{b}_2b_0\bar{z}_{6}z_{10}a_1; \quad b_4\mapsto a_0. 
    \end{array}
\]
and control zeta function $\zeta_{C}(t)=5+9t+17t^2+\cdots$.

\fig{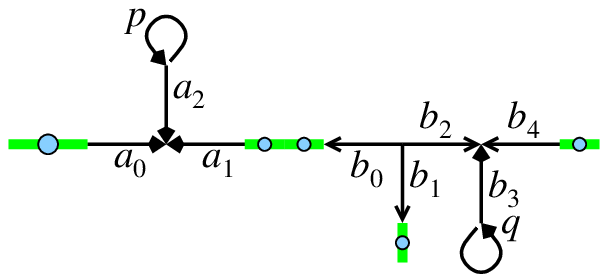}{The optimal graph representative.}
Finally, we fold up $p_1$ and $q_1$ to give new edges $a_2$ and $b_3$ as shown in \figref{algorithm2efficient}.
The graph map is now
\[ \begin{array}{c} 
     \quad p_1\mapsto q_1; \quad q_1\mapsto p_1; \\[1ex]
     a_0\mapsto a_0\bar{a}_1\bar{z}_{10}z_{6}\bar{b}_0; \quad a_1\mapsto \bar{b}_1; \quad a_2\mapsto b_3\bar{b}_2 \\[1ex]
     b_0\mapsto \bar{b}_3\bar{q}_1b_3\bar{b}_2b_1; \quad b_1\mapsto \bar{b}_4; \quad b_2\mapsto \bar{b}_2b_0\bar{z}_6z_{10}a_1; \quad b_3\mapsto a_2; \quad b_4\mapsto a_0.
    \end{array}
\]
which is optimal.
The control zeta function is $\zeta_{g;C}(t)=5+9t+15t^2+\cdots$.
The edges $p_1$ and $q_1$ are peripheral edges for the graph map.
\end{example}


\subsection{Algorithm for reduction}
\label{sec:reductionalgorithms}

We finally give algorithms for performing reductions.
We give two algorithms, one for separating reductions, and one for non-separating reductions.
These algorithms are modelled on those of Keil \cite{Keil97}, but are simpler to apply in practice, since they do not require $2$-simplices.
Since reduction is not the main focus of this paper, and the algorithms do not introduce any significant new ideas, we provide only a brief analysis.

\begin{algorithm}[Separating reduction]\ 
Let $g$ be a controlled graph map with an invariant subgraph $H$ consisting of $n$ components, all of which have negative Euler characteristic and contain no control edges.
Without loss of generality, we can assume the components of $H$ are cyclically permuted by $g$.
\begin{enumerate}\isz
\item\steplabel{separatingstart} Let $p_1$ be an edge in $H$ such that the preceding edge in the cyclic order of vertices at $\init{p_1}$ does not lie in $H$.
  Let $\pi$ be the edge path $p_1$.
\item\steplabel{separatingperipheralisecases} Let $p_k$ be the last edge of $\pi$, with $v_k=\final{p_k}$,
  and let $p_{k+1}$ be next edge in the cyclic order of edges at $v$ which lies in $H$.
  \begin{itemize}\isz
   \item If $v_k$ is not the initial vertex of an edge in $\pi$, adjoin $p_{k+1}$ to the end of $\pi$ and return to \stepref{separatingperipheralisecases}.
   \item If $p_{k+1}$ lies in $\pi$, this step terminates.
  \end{itemize}
  Otherwise, split off all edges from $\bar{p}_k$ to $p_{k+1}$ inclusive.
  \begin{itemize}\isz
    \item If $\bar{p}_{k+1}$ is not a directed edge in $\pi$, adjoin $p_{k+1}$ to the end of $\pi$ and return to \stepref{separatingperipheralisecases}.
    \item Otherwise, $\bar{p}_{k+1}=p_i$ in $\pi$.
      Homotope so that all edges at the new vertex $v_k\p$ so that $g\p(p_{k+1})=\cdot$, collapse $p_{k+1}$ and return to \stepref{separatingperipheralisecases}.
  \end{itemize}
\item Let $P$ be the frontier loop of $H$.
\item Restrict $g$ to a map $g\p$ of $G\p=G\setminus H \cup P$.
  Tighten this restricted map.
\item The output of the algorithm is a the pair of controlled graph maps $g^n\restrict{K}$, where $K$ is a component of $H$, and $g\p$ on $G\p$.
\end{enumerate}
\end{algorithm}

This algorithm works since once we have a situation where the frontier of $H$ consists of simple closed curves,
 tightening $g$ restricted to $G\p=G\setminus H\cup P$ gives a graph map $g\p$ of $G\p$ by virtue of the fact that $g$ is embeddable
 and acts as an automorphism of the fundamental group $\pi_1(G)$.

\begin{algorithm}[Non-separating reduction]\ 
\begin{enumerate}\isz
\item\steplabel{nonseparatingstart} Homotope to collapse all edges which eventually map into $H$.
\item\steplabel{nonseparatingtighten} Tighten $g$ restricted to $H$, homotoping vertices if necessary.
\item Consider a component of $H$ with vertices $v_0,v_1\ldots v_{k-1},v_k=v_0$ and edges $p_1,\ldots p_k$ with $\init{p_i}=v_{i-1}$ and $\final{p_i}=v_i$.
  Create two new components of $G$, one with vertices $v_i^+$ and edges $p_i^+$ from $v_{i-1}^+$ to $v_{i}^+$ and one with vertices $v_i^-$
   and edges $p_i^-$ from $v_{i}^-$ to $v_{i-1}^-$.
  Move all the edges at vertex $v_i$ in the interval $(\bar{p}_{i-1},p_i)$ to the interval $(\bar{p}_{i-1}^+,p_i^+)$ at $v_i^+$,
   and the edges in the interval $(p_i,\bar{p}_{i-1})$ to the interval $(\bar{p}_i^-,p_{i-1}^-)$ at $v_i^-$.
\item Let $p$ be an edge in $H$, and suppose $g$ is orientation-preserving.
  If $g(p)=q$ is a positively oriented edge, take $g\p(p^+)=q^+$ and $g\p(p^-)=q^-$, and if $g(p)=\bar{q}$, take $g\p(p^+)=q^-$ and $g\p(p^-)=q^+$.
  If instead $g$ is orientation-reversing, take $g\p(p^+)=\bar{q^-}$ and $g\p(p^-)=\bar{q}^+$ if $g(p)=q$, and $g\p(p^+)=\bar{q}^+$ and $g\p(p^-)=\bar{q}^-$ if $g(p)=\bar{q}$. 
\item Let $e$ be an edge not in $H$ with $g(e)=\epsilon\pi\ldots$ where $\epsilon$ is an edge path with final edge $e_k$ not in $H$,
    and $\pi=p_1\ldots p_j$ is an edge path in $H$.
  If $\final{e_k}=v^+$ for some vertex, take $g\p(e)=\epsilon\pi^+\ldots=\epsilon p_1^+p_2^+\ldots p_j^+\ldots$, and if $\final{e_k}=v^-$,
   take $g\p(e)=\epsilon {\bar{p}_1}^{\,-} \bar{p}_2^{\,-} \ldots \bar{p}_j^{\,-} \ldots$.
\item $H$ is not connected to the new graph $G\p$ and can be removed.
\item The output of the algorithm is a the controlled graph map $g\p$ of $G\p$, together with an identification of pairs of peripheral loops $P^\pm$.
\end{enumerate}
\end{algorithm}

This algorithm works since the construction of the splitting of $H$ into two copies of itself keeps track of the connections to the rest of the graph,
 and this allows us to easily construct a continuous graph map $g\p$.

After \stepref{nonseparatingstart}, all edges not in $H$ map to at least one edge in $H$.
After \stepref{nonseparatingtighten}, the edges in $H$ are permuted; in particular, $g$ is a homeomorphism on $H$.
The final graph map is continuous since the separate definitions of $g\p$ on $H\p$ and $G\p\setminus H\p$ match up by considering orientations.
In particular, note that since $g$ is a homeomorphism on $H$ and is embeddable, we cannot have an edge path $\epsilon$ in $G\setminus H$
 which enters a component of $H$ on one side and leave on the other.

\section{Symbolic Dynamics and Shadowing}
\label{sec:shadowing}

In this section, we discuss ways of representing the dynamics of a chaotic map.
There are two main methods, the standard technique of symbolic dynamics,
 and the more powerful notion of \emph{global shadowing}, which allows us to relate one map to another which is well understood.


\subsection{Symbolic dynamics}
\label{sec:symbolic}

\begin{definition}[Itinerary and code]
\label{defn:itinerary}
Let $f$ be a map of a space $X$, and $\{R_{\lambda}\}_{\lambda\in\Lambda}$ be a finite set of compact subsets of $X$.
Let $(x_i)$ be an orbit of $f$, so $f(x_{i-1})=x_i$ for all $i$.
Then a bi-infinite sequence $\ldots R_{\lambda_{-2}}R_{\lambda_{-1}}R_{\lambda_0}R_{\lambda_1}R_{\lambda_2}\ldots$ is an \emph{itinerary} for $(x_i)$
 if $x_i\in R_{\lambda_i}$ for all $i$.
If $x_0$ is a point of $X$, an $x_i=f^i(x_0)$ for all $i\geq0$, then the infinite sequence $R_{n_0}R_{n_1}R_{n_2}\ldots$
 is an itinerary for $x_0$ if $x_i\in R_{n_i}$ for all $i\geq0$. 
If $x_0$ is a period-$n$ point of $f$, then the word $R_{\lambda_0}R_{\lambda_1}\ldots R_{\lambda_{n-1}}$ is a \emph{code} for $x_0$ if $x_i\in R_{\lambda_i}$ for $0\leq i<n$.
\end{definition}
Note that if the sets $R_\lambda$ do not cover $X$, an orbit may not have any itinerary
 (if for some $i$, $x_i$ is not in $R_\lambda$ for any $\lambda$),
 and may have more than one itinerary if the sets $R_\lambda$ are not disjoint.
Typically, the sets $R_\lambda$ are chosen to have disjoint interiors and boundaries with measure zero,
 so that the set of points with more than one itinerary has measure zero.

If $f:(X,Y)\fto(X,Y)$ is a map of pairs, we typically take the sets $R_\lambda$ to have boundary in $Y$, and call such sets \emph{regions}.
If $(f;T)$ is a trellis map, we take itineraries of $\cut{f}$ using sets $R_\lambda$ which are regions of $(\cut[U]{T},T^S)$.
The itineraries of orbits of $\cut{f}$ then give itineraries of orbits of $f$.
One major advantage of trellis maps over other maps of surfaces is that the regions are specified in a natural way.

\begin{definition}[Shift space]
The \emph{shift space} $\Sigma_{\Lambda}(f)$ of a map $f$ on the regions $\{R_\lambda:\lambda\in\Lambda\}$ is the subset of $\{R_\lambda:\lambda\in\Lambda\}^\N$ consisting
 of itineraries of points under $f$.
\end{definition}

The following result shows how the symbolic dynamics of the graph representative of a trellis mapping class relates to the dynamics of a trellis map.
The proof is given in \cite{Collins99AMS}.
\begin{theorem}
\label{thm:symbolicdynamics}
Let $(g;G,W)$ be the graph representative of an irreducible trellis mapping class $([f];T)$, with $G$ embedded as a subset of $\cut[U]{T}$.
Then for every orbit $(y_i)$ of $g$, there is an orbit $(x_i)$ of $f$ with the same itinerary as $(y_i)$.
Further, if every region of $T$ is simply-connected, then there is at most one orbit of $g$ with a given itinerary (apart from orbits which enter the control set $C$).
\end{theorem}


\subsection{Global shadowing}
\label{sec:shadowingdefinition}

Global shadowing is an equivalence relation on orbits of a system which gives a condition under which two orbits can be considered to be close to each other
 over infinite time intervals.
The classical definition was introduced by Katok, and is given in terms of lifts of the map to the universal cover.
\begin{definition}
\label{defn:katokshadow}
Let $X$ be a Riemannian manifold with metric $g$ and distance function $d$, and $f:X\fto X$ a continuous map.
Let $\tilde{X}$ be the universal cover of $X$, let $\tilde{d}$ be the equivariant distance function obtained from the equivariant metric $\tilde{g}$ covering $g$, and let $\tilde{f}:\tilde{X}\fto\tilde{X}$ be a lift of $f$.
Then orbits $(x_i)$ and $(y_i)$ \emph{globally shadow} each other if they lift to orbits $(\tilde{x}_i)$ and $(\tilde{y}_i)$ of $\tilde{f}$
 such that $\{\tilde{d}(\tilde{x}_i,\tilde{y}_i)\}$ is bounded.
\end{definition}
This definition is independent of the metric used, and can be extended to certain equivariant metrics which are non-Riemannian.

As the concept is topological in spirit, it would be useful to have a purely topological classification.
We will also be concerned with spaces which are topological pairs, and seek a definition which extends naturally to this setting.
We would also like to talk about arbitrary sets of orbits as well as just pairs or orbits.
As a frequently-recurring theme in this work is that of homotopy classes of curves, a definition in terms of this formalism,
 rather than that of universal covers will also be useful.
The basic idea of our global shadowing relation is given below.
It relies on the definition of sets of curves of \emph{bounded lengths}.
\begin{definition}[Global shadowing]
\label{defn:globalshadow}
Orbits $(x_i)$ and $(y_i)$ of a map $f:(X,Y)\fto (X,Y)$ are said to \emph{globally shadow} each other if there are sets $J_i\subset I$ such that $J_i\subset J_{i+1}$ 
 and exact curves $\alpha_i:(I,J_i)\exto(X,Y)$ from $x_i$ to $y_i$ such that $f\circ\alpha_i$ is homotopic relative to endpoints to $\alpha_{i+1}$
 via curves $(I,J_{i+1})\fto(X,Y)$, and the curves $\alpha_i$ have \emph{bounded lengths}.
The curves $\alpha_i$ are called \emph{relating curves}.

Similarly, $(x_i)$ and $(y_i)$ \emph{forwards shadow} each other if there are shadowing curves for $i\geq0$,
 and \emph{backwards shadow} each other if there are shadowing curves for $i\leq0$.
\end{definition}
Note that we really need to allow the $J_i$ to be an increasing sequence of sets; for example, we could take $(x_i)$ and $(y_i)$ to be orbits which start in $X$ and map into $Y$.

There are many possible definitions of the term ``bounded length,'' and these mostly give identical definitions of global shadowing.
In the setting of smooth curves on Riemannian manifolds, we can use the natural definition of the length of a smooth curve.
In the setting of absolute neighbourhood retracts, we can use the following definition
\begin{definition}[Bounded lengths]
Let $X$ be a compact ANR. 
A set of curves $\Gamma$ has \emph{bounded lengths} if for any closed ANR subset $Z$ of $X$,
 there is a number $N(Z)$ such that each curve $\gamma\in\Gamma$ is homotopic to a curve $\alpha$ such that $\alpha(I)\cap Z$ has at most $N$ components.
\end{definition}

The following lemmas give some technical results which simplify the application of the global shadowing definition.
The proofs are straightforward, so are omitted.
\begin{lemma}
Let $X$ be a manifold, and let $Z$ be a set of disjoint hypersurfaces in $X$ such that $X\setminus Z$ simply-connected.
Suppose there is a number $N$ such that for any $\gamma\in\Gamma$, $\gamma$ has at most $n$ intersections with $Z$.
Then the curves in $\Gamma$ have bounded lengths.
\end{lemma}
\begin{lemma}
Let $\Gamma$ be a set of curves with the same endpoints which have bounded lengths.
Then there set of isotopy classes $\{ [\gamma] : \gamma\in\Gamma \}$ is finite.
\end{lemma}
\begin{lemma}
Let $U$ be a contractible open subset of a space $X$.
Then there is a set of curves $\Gamma$ of bounded lengths contained in $U$ such that for any two points $x_0,x_1$ in $U$,
 there is a curve in $\Gamma$ from $x_0$ to $y_0$.
\end{lemma}

Given a homotopy between two different maps, we can also define what it means for an orbit of one map to globally shadow an orbit of the other map.
\begin{definition}[Global shadowing for homotopic maps]
\label{defn:globalshadowhomotopic}
Let $f_t:(X,Y)\fto(X,Y)$ be a homotopy from $f_0$ to $f_1$, and $F:(X\times I, Y\times I)\fto(X\times I,Y\times I)$ be the fat homotopy of $f_t$ defined by $F(x,t)=(f_t(x),t)$.
Then an orbit $(x_i)$ of $f_0$ and an orbit $(y_i)$ of $f_1$ are said to globally shadow each other under the homotopy $f_t$
 if the orbits $(x_i,0)$ and $(y_i,1)$ of $F$ globally shadow each other under Definition~\ref{defn:globalshadow}.
\end{definition}


\subsection{Shadowing and Nielsen numbers}
\label{sec:nielsen}

One of the main tools we have to relate the dynamics of maps with the same homotopy or isotopy type is Nielsen periodic point theory.
\begin{definition}[Nielsen equivalence]
Let $f$ be a self-map of a topological pair $(X,Y)$.
Period-$n$ orbits $(x_i)$ and $(y_i)$ of $f$ are \emph{$n$-Nielsen equivalent} if there is a subset $J$ of $I$ and exact curves $\alpha_0,\ldots \alpha_{n-1}:(I,J)\fto(X,Y)$ from $x_i$ to $y_i$ such that $f\circ\alpha_i\homotopic\alpha_{i+1\;\mod\;n}$ (taking homotopies relative to endpoints).
\end{definition}

Nielsen equivalence is the same as global shadowing for periodic orbits.
\begin{lemma}
\label{lem:nielsenshadow}
Periodic orbits $(x_i)$ and $(y_i)$ globally shadow each other if and only if they are $n$-Nielsen equivalent for some $n$.
\end{lemma}

\begin{proof}
Clearly, if $(x_i)$ and $(y_i)$ are Nielsen equivalent, they globally shadow each other.
Conversely, if $(x_i)$ and $(y_i)$ globally shadow each other with shadowing curves $(\alpha_i)$, then the possible homotopy classes for the $\alpha_i$ are bounded and the
 number of components of the $J_i$ can be made to be bounded by suitable choices of the $\alpha_i$.
Hence, the $J_i$ can be taken to be fixed, say $J_i=J$ for sufficiently large $i$.
Then the homotopy classes of the $\alpha_i$ must eventually repeat, say $\alpha_{k+n}\homotopic\alpha_k$ for some $k\in\Z$ and $n>0$.
We can then take the curves to be $\beta_i=\alpha_{k+i}$ for $0\leq i<n$.
\end{proof}

It is easy to show \cite{Collins01TOPOA} that Nielsen equivalence classes are closed in $X$, and hence have a well-defined fixed point index.
\begin{definition}[Nielsen number]
An $n$-Nielsen class is \emph{essential} if its fixed-point index under $f^n$ is nonzero.
The $n$th Nielsen number of $f$, denoted $N_n(f)$ is the number of essential Nielsen classes of period-$n$ orbits of $f$.
If $R_{\lambda_0},R_{\lambda_1},\ldots,R_{\lambda_{n-1}}$ is a code, then the its Nielsen number is the number of essential Nielsen classes with orbits of the code.
\end{definition}

The most important property of exact homotopy equivalence is that is preserves the Nielsen numbers.
Relative Nielsen theory can therefore be used to find periodic orbits and an entropy bound for a map of pairs.
A complete exposition is given in \cite{Collins01TOPOA}.

The logarithm growth rate of the Nielsen numbers gives the \emph{Nielsen entropy}.
\begin{definition}[Nielsen entropy]
The \emph{Nielsen entropy} of a map $f$ is the equal to the logarithm of the growth rate of the $n$th Nielsen numbers,
\[
  \hniel(f)=\limsup_{n\tendsto\infty}\frac{\log N_n(f)}{n} \;.
\]
\end{definition}


\subsection{Shadowing by trellis orbits}
\label{sec:trellisshadowing}

It is important to know how the dynamics of the graph representative models that of the trellis map.
In Nielsen-Thurston theory, the orbits of the train-track map globally shadow those of the original surface homeomorphism.
The results stated in Theorem~\ref{thm:shadowing} are analogous to the shadowing results of Nielsen-Thurston theory.
We show that orbits of $g$ are shadowed by orbits of $f$, and further, that periodic, asymptotic and biasymptotic orbits of $g$ are shadowed by periodic, asymptotic and
 biasymptotic orbits of $f$.

We first need to clarify what it means for an orbit of $g$ to be backward asymptotic to the control set $C$.
\begin{definition}[Asymptotic graph orbits]
An orbit $(y_i)$ of $g$ is \emph{forward asymptotic} to $C$ if $y_j\in C$ for some $j$.
An orbit $(y_i)$ of $g$ is \emph{backward asymptotic} to $C$ if there exists $k\in\Z$ and exact curves $\mu_i:(I,\{0,1\})\exto(G,W)$ such that $\mu_i$ begins in a periodic point of $W$, $\mu_i$ is homotopic in $\cut{T}$ to a curve in $T^U$, $g\circ\mu_i\supset \mu_{i+1}$ for $i<k$ and $y_i\in\mu_i$.
\end{definition}
Notice that if $(y_i)$ is backward asymptotic to $T^P$, then the $\alpha$=limit set of $(y_i)$ lies in the control set $C$.
The condition that the curves are homotopic to $T^U$ is crucial; otherwise, we may have a periodic orbit of $g$ which has $\alpha$-limit set in $C$, but is not shadowed by an orbit of $\cut{f}$ in $T^U$.

\begin{theorem}[Shadowing by trellis orbits]
\label{thm:shadowing}
Let $(f;T)$ be a trellis map and $(g;G,W)$ an efficient controlled graph map compatible with $(f;T)$.
Then for any orbit $(y_i)$ of $g$ there is an orbit $(x_i)$ of $f$ which shadows $(y_i)$.
Further,
\begin{enumerate}\isz
\item\label{item:periodicshadowing} If $(y_i)$ is periodic, then $(x_i)$ can be chosen to be periodic.
\item\label{item:stableshadowing} If $(y_i)$ is forward asymptotic to $C$ for some $j$, then $(x_i)$ can be chosen in the stable set of $T^P$ (so $x_j\in T^S$ for some $j$).
\item\label{item:unstableshadowing} If $(y_i)$ backward asymptotic to $C$, then $(x_i)$ can be chosen in the unstable set of $T^P$.
\item\label{item:biasymptoticshadowing} If $(y_i)$ is both forwards and backwards asymptotic to $C$, then $(x_i)$ can be chosen to be biasymptotic to $T^P$.
\end{enumerate}
\end{theorem}
Note that due to the asymmetry in our treatment of stable and unstable curves, part~(\ref{item:stableshadowing}) and part~(\ref{item:unstableshadowing}) are genuinely different statements and require different proofs.
As an immediate corollary, we have the following result.
\begin{theorem}[Entropy of efficient graph maps]
If $g$ is an efficient controlled graph map compatible with $(f;T)$, then $\htop(g)\leq\htop(f)$.
\end{theorem}
We shall see in Section~\ref{sec:entropy} that this entropy bound is sharp if $T$ is well-formed, but otherwise need not be.

The remainder of this section is devoted to proving Theorem~\ref{thm:shadowing}.
We assume throughout that $(f;T)$ is a trellis map with graph representative $(g;G,W)$, and that $(G,W)$ is embedded as a subset of $(\cut[U]{T},T^S)$ by an embedding $i$.
We let $\pi$ be a deformation retract $(\cut[U]{T},T^S)\exto(G,W)$, and $f_t$ be a homotopy from $\cut{f}=f_0$ to $f_1=i\circ g\circ\pi$ with fat homotopy $F$.

The following theorem is contained in \cite{Collins99AMS}. 
It shows that periodic orbits of $g$ are globally shadowed by an orbit of $\cut{f}$, and forms the basis for all other shadowing results.
The result is proved by showing that all Nielsen classes of $g$ are essential.
\begin{theorem}[Shadowing by periodic orbits]
Let $(y_i)$ be a period-$n$ orbits of $g$.
Then there is a period-$n$ orbit $(x_i)$ of $\cut{f}$ which is Nielsen equivalent to $(y_i)$.
\end{theorem}

To show that all orbits of $g$ are shadowed by orbits of $f$, we need to find a uniform bound for the lengths of shadowing curves in $g$.
\begin{lemma}
If $(x_i)$ and $(y_i)$ are orbits of $g$ which globally shadow each other,
 there are relating curves $\alpha_i$ from $x_i$ to $y_i$ such that the $\alpha_i$ lie in $C$, $P$ or $\pre{P}$.
\end{lemma}

\begin{proof}
Suppose the curves $\alpha_i$ are tight and do not cross any control edges.
Then $\alpha_{i+1}=g\circ\alpha_i$ for all $i$, so the $\alpha_i$ are eventually periodic or disjoint.
Since $g$ has no invariant subgraphs apart from $P$, and $g$ is expanding outside $\pre{P}$,
 the $\alpha_i$ can only be eventually periodic if they lie in $\pre{P}$.

If the $\alpha_i$ do cross $C$, they must eventually do so in periodic control edges.
We can split the curves into pieces that cross $C$ and pieces that do not.
The pieces that do not must then lie in $\pre{P}$.
\end{proof}

Thus there is a uniform bound on the length of any relating curve in $G$.
The following result shows that all orbits of $g$ are shadowed by those of $f$.

\begin{lemma}
\label{lem:shadow}
Let $(y_i)$ be an orbit of $(g;G,W)$.
Then there is an orbit $(x_i)$ of $(f;T)$ which shadows $(y_i)$.
\end{lemma}

\begin{proof}
If $(y_i)$ is periodic, the result follows from the relative Nielsen theory.
If not, suppose there are periodic orbits $(y_i^j)$ such that for every $i$, $\lim_{j\tendsto\infty}y_i^j=y_i$.
We choose orbits $(x_i^j)$ of $f$ such that $(x_i^j,0)$ shadows $(y_i^j,1)$ under $F$.
Since $M$ is compact, by choosing a subsequence if necessary, 
 we can assume the sequence $x_0^j$ converges to a point $x_0$ as $j\tendsto\infty$,
 and this means that $x_i^j$ converges to $x_i=f_0^i(x_0)$ for all $i>0$.
Choosing further subsequences, we can assume further that $x_i^j$ converges to a point $x_i$ as $j\tendsto\infty$ for any $i<0$.

We therefore have an orbit $(x_i)$ such that $x_i=\lim_{j\tendsto\infty} x_i^j$ for all $i$.
We now let $\pi_i^j$ be the relating curves for $(x_i^j,0)\sim (y_i^j,1)$.
Take an open cover $\mathcal{U}$ of $(\cut[U]{T}\times I, T^S\times I)$ by sets $U$ such that the components of $U\setminus T^S$ are simply-connected.
Let $\beta_i^j$ and $\gamma_i^j$ be curves from $(x_i,0)$ to $(x_i^j,0)$ and from $(y_i,1)$ to $(y_i^j,1)$ in a component of $U$ and $V$ respectively,
 which must exist for large enough $j$.
Then $f_0\circ\beta_i^j\homotopic \beta_{i+1}^j$ and $f_1\circ\gamma_i^j\homotopic \gamma_{i+1}^j$ whenever $j$ is sufficiently large.

Now let $\alpha_i^j=\beta_i^j\cdot\pi_i^j\cdot\gamma_i^j$.
The curves $\alpha_0^j$ have bounded lengths, since they are joins from three sets of curves which each have bounded lengths.
Therefore, there must be finitely many homotopy classes $[\alpha_0^j]$, at least one of which must occur infinitely often.
By taking an subsequence if necessary, we can assume $[\alpha_0^j]$ has the same value for all $j$.
Now, 
  \[ F\circ\alpha_i^j = F\circ\beta_i^j \cdot F\circ\pi_i^j \cdot F\circ\gamma_i^j
       \homotopic \beta_{i+1}^j\cdot\pi_{i+1}^j\cdot\gamma_{i+1}^j = \alpha_{i+1}^j \]
and since this must be true whenever $\alpha_i^j$ is defined, we have $[\alpha_i^j]$ is independent of $j$ for fixed $i$.
By taking subsequences if necessary, we can also ensure this is the case for $i$ negative
 (or this follows automatically if $F$ has a homotopy inverse).
Hence, we have curves $\alpha_i$ from $x_i$ to $y_i$ of bounded lengths such that $F\circ\alpha_i\homotopic\alpha_{i+1}$.
\end{proof}


\subsection{Shadowing by stable and biasymptotic orbits}
\label{sec:asymptoticshadow}

We now show that orbits of $g$ which enter $C$ are shadowed by orbits of $f$ in $W^S(T^P)$, that orbits backward asymptotic to $C$ are shadowed by orbits of $f$ in $W^U(T^P)$,
 and that orbits of $g$ which are biasymptotic to $C$ are shadowed by orbits of $f$ in $W^U(T^P)\cap W^S(T^P)$.
These results have no corresponding statements in Nielsen-Thurston theory, since we cannot guarantee the existence of periodic saddle points in this case.
We first prove an initial result which shows that orbits of $g$ are contained in curves with endpoints in $W$ which map over each other under $g$.
If $\alpha_0:(I,J)\exto(G,W)$ is a tight exact curve in $G$ with endpoints in $W$, then we write $g(\alpha_0)\supset\alpha_1$ if $g(\alpha)$ contains $\alpha_1$ as a sub-curve.

\begin{lemma}
\label{lem:graphcurves}
Let $(y_i)$ be an orbit of $g$.
Then there is a bi-infinite sequence of exact curves $\mu_i:(I,\{0,1\})\exto(G,W)$ such that $y_i\in\mu_i(I)$ and $g\circ\mu_i$ contains $\mu_{i+1}$ as a sub-curve.
\end{lemma}

\begin{proof}
First suppose that every component of $G\setminus W$ is homotopy-equivalent to a point or circle.
Then every pair of points of $W$ in the same region, there are two tight curves between them (for an homotopy circle) and one for a homotopy point.
Let $\Gamma$ be the set of all these curves.
Let $\Sigma_n$, $n\geq 0$ be the set of bi-infinite sequences of elements of $\Gamma$ such that $g(\gamma_i)\subset\gamma_{i+1}$ for $i\geq-n$.
This is nonempty, since if $y_i\in\gamma_i$, then $y_{i+1}\in g(\gamma_i)$, and is compact since $\Gamma$ is finite.
Clearly $\Sigma_{n+1}\subset \Sigma_n$ so the $\Sigma_n$ form a sequence of compact nested sets, hence $\Sigma=\bigcap_{n=0}^{\infty}\Sigma_n$ is nonempty.
Then any sequence $(\mu_i)$ is a sequence in $\gamma$.

Now if $g$ is irreducible optimal graph map, there exists $N$ such that the complement of $g^{-N}(W)$ consists of homotopy points and circles.
We consider $g$ as a self-map of $(G,g^{-N}(W))$, and let $(\gamma_n)$ be a sequence of curves found by the previous arguments.
Then $g^N(\gamma_{i-N})$ is an exact curve $(I,\{0,1\})\exto(G,W)$ containing $(y_i)$, and maps across $g^N(\gamma_{i+1-N})$.
Therefore taking $\mu_i=g^N(\gamma_{i-N})$ gives the required bi-infinite sequence of curves.
\end{proof}

\begin{lemma}
\label{lem:stableshadow}
Let $(y_i)$ be an orbit of $g$ with $y_0\in C$. Then there is an orbit $(x_i)$ of $f$ which globally shadows $(y_i)$ with $x_0\in T^S$.
\end{lemma}

\begin{proof}
Let $(\mu_i)$be a sequence of exact curves $(I,\{0,1\})\fto(G,W)$ such that $g\circ\mu_i$ contains $\mu_{i+1}$ as a sub-curve.
 and $y_i\in\mu_i$.
Let $S_i$ be the stable segment containing $y_i$ for $i\geq0$; without loss of generality we assume $y_{-1}\not\in T^S$.

We now use recursively construct a sequence of segments $S_i$ of $W^S(T^P)$ for $i<0$.
Suppose therefore for $i>n$ (where $n$ is negative) we have a stable segment $S_i$ of $f^{-n}(T^S)$ with an essential intersection with $\mu_i$ at a point $z_i$
 and that $f(S_i)\subset S_{i+1}$.
Then $g\circ\mu_n$ has an essential intersection with $S_{n+1}$ at a point $z_{n+1}$, so $f\circ\mu_n$ has an essential intersection with $S_{n+1}$.
Therefore, one of the segments $S_n$ of $f^{-1}(S_{n+1})$ is such that $f\circ\mu_n$ contains an essential intersection with $f(S_n)$ in the same intersection class as $z_{n+1}$ under the homotopy $f_t$.
Let $f(z_n)$ be this essential intersection, so that $z_n$ is an essential intersection of $S_n$ with $\mu_n$, as shown in \figref{stableshadowing}.
\fig{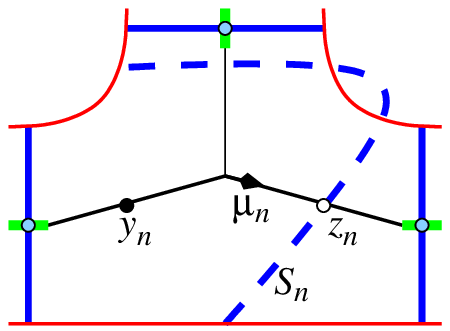}{Construction of $z_n$.}

The sequence of segments $S_n$ of $W^S(T^P)$ is such that $f(S_n)\subset S_{n+1}$ for $n<0$, so $f^{n}(S_n)\subset f^{n+1}(S_{n+1})\subset S_0$.
Hence the sets $f^n(S_n)$ are a nested sequence of compact intervals, so contain a limit point $x_0$.
We claim the orbit $(x_i)=f^i(x_0)$ globally shadows $(y_i)$ under the homotopy $f_t$.

We need to show that the orbits $(x_n,0)$ and $(y_n,1)$ globally shadow each other under $F$.
Let $\alpha_n$ be a curve from $(x_n,0)$ to $(z_n,0)$ in $S_n\times\{0\}$, $\beta_n$ a curve from $(z_n,0)$ to $(g^{-1}(z_{n+1}),1)$ in $\mu_{n+1}\times I$
 and $\gamma_n$ the curve from $(g^{-1}(z_{n+1}),1)$ to $(y_n,1)$ in $\mu_n(I)\times\{1\}$.
Then $F\circ\alpha_n$ is a curve in $T^S\times I$ and $F\circ\gamma_n$ is a curve in $\mu_{n+1}\times\{1\}$.
The curve $F\circ\beta_n$ goes from $(f(z_n),0)\in S_{n+1}\times\{0\}$ to $(z_{n+1},1)\in S_{n+1}\times\{0\}$,
 and is homotopic to a curve whose path lies in $S_{n+1}\times I$ since $z_{n+1}$ is an essential intersection of $S_n$ with $\mu_n$.
Hence $F\circ(\alpha_n\cdot\beta_n)$ is homotopic to a curve in $S_{n+1}\times I$ from $(x_{n+1},0)$ to $(z_{n+1},1)$ and $F\circ\gamma_n$ is homotopic to a curve from $(z_{n+1},1)$ to $(y_{n+1},1)$.
Since $\beta_{n+1}$ is homotopic to a curve in $\{z_{n+1}\}\times I$ joined with a curve in $\mu_{n+1}(I)\times\{1\}$, we have $F\circ(\alpha_n\cdot\beta_n\cdot\gamma_n)\homotopic\alpha_{n+1}\cdot\beta_{n+1}\cdot\gamma_{n+1}$ as required.
\end{proof}

\begin{lemma}
\label{lem:unstableshadow}
Let $(y_i)$ be an orbit of $g$ such that for $i\leq0$ there are exact curves  $\mu_i:(I,\{0,1\})\fto(G,W)$ containing $y_i$ such that $g(\mu_i)\supset \mu_{i+1}$, $\mu_i$ is unstable-parallel, has periodic initial vertex.
Then there is an orbit $(x_i)$ of $f$ which shadows $(y_i)$ such that $x_0\in T^U$.
\end{lemma}

\begin{proof}
For $i\geq0$, let $\sigma_i$ be a cross-cut in $\cut[U]{T}$ which does not intersect $T^S$, but which crosses $G$ exactly once in the edge containing $y_i$.
Let $U_0$ be the segment of $T^U$ which is exact homotopic to $\mu_i$.
Suppose for $i\leq n$, there is a segment $U_i$ of $W^U(T^P)$ such that $U_i\supset U_{i-1}$ and $U_i$ is exact homotopic to a tight curve $\mu_i$ in $G$ containing $y_i$.
Then $U_i$ has an essential intersection $z_i$ with $\sigma_i$.
Since $y_{i+1}\subset g\circ\mu_i$, there must be an essential intersection of $g\circ\mu_i$ with $\sigma_i$, and this intersection lies in the same intersection class as an intersection $z_{i+1}$ of $f(U_i)$ with $\sigma_{i+1}$.
Let $U_{i+1}$ be the segments of $W^U$ containing $z_{i+1}$, and let $\mu_{i+1}$ be the tight curve in $G$ which is exact homotopic to $U_i$, and so $y_{i+1}\in \mu_{i+1}$
Then the set $U=\bigcap_{n=0}^{\infty}f^{-n}(U_n)$ is a nonempty subset of $T^U$.
\fig{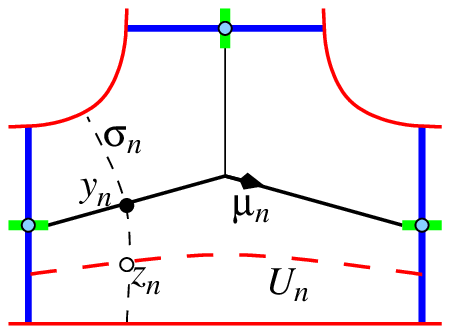}{Construction of $U_n$.}

We take $x_0\in U$, and $(x_i)=f^i(x_0)$ for $i\neq0$.
A similar analysis to that of Lemma~\ref{lem:stableshadow} shows that $(x_i)$ globally shadows $(y_i)$.
\end{proof}

\begin{lemma}
Let $(y_i)$ be an orbit of $g$ such that for $y_k\in C$ and for $i\leq0$ there are exact curves  $\mu_i:(I,\{0,1\})\fto(G,W)$ containing $y_i$ such that $g(\mu_i)\supset \mu_{i+1}$, $\mu_i$ is unstable-parallel, has periodic initial vertex.
Then there is an orbit $(x_i)$ of $f$ which shadows $(y_i)$ such that $x_0\in T^U$ and $x_k\in T^S$.
\end{lemma}

\begin{proof}
Without loss of generality we can assume that $y_k\in W$, but $y_{k-1}\not\in C$.
Take $\sigma_i$ as in the proof of Lemma~\ref{lem:unstableshadow} for $i<k$, and let $\sigma_k$ to be the stable segment $S_k$ containing $y_k$.
Construct $U_i$ as in the proof of Lemma~\ref{lem:unstableshadow} for $i\leq k$.
Then $z_k$ is a point in $U_k\cap S_k$.
We take $x_i=f^{i-k}(z_k)$.
Then a similar analysis to that of Lemma~\ref{lem:stableshadow} shows that $(x_i)$ globally shadows $(y_i)$.
Construct $S_i$ for $i\geq n$ as in the proof of Lemma~\ref{lem:stableshadow}.
\end{proof}
This completes the proof of Theorem~\ref{thm:shadowing}


\section{Entropy-Minimising Diffeomorphisms}
\label{sec:entropy}

In this section, we show that the entropy bounds obtained in Section~\ref{sec:shadowing} are sharp, at least for well-formed trellises.
That is, the topological entropy of the graph representative $g$ for a trellis type $[f;T]$ (which is the same as the Nielsen entropy $\hniel[f;T]$)
 is the infemum of the topological entropies of diffeomorphisms in the class.
The main results of this section are summarised in the following theorem.

\begin{theorem}[Existence of entropy minimising models]
\label{thm:mainentropy}
Let $([f];T)$ be a well-formed trellis mapping class.
Then for any $\epsilon>0$ there exists a diffeomorphism $\widehat{f}\in([f];T)$ such that $\htop(\widehat{f})<\hniel[f;T]+\epsilon$.
If there exists a diffeomorphism $\widehat{f}\in([f];T)$ such that any $\widehat{f}$-extension of $T$ is minimal, then there is a uniformly-hyperbolic diffeomorphism $\widetilde{f}\in([f];T)$ such that $\htop(\widetilde{f})=\hniel[f;T]$.
Further, if $([f];T)$ is irreducible and $\widetilde{f}\in([f];T)$ such that $\htop(\widetilde{f})=\hniel[f;T]$, then any $\widetilde{f}$-extension of $T$ is minimal.
\end{theorem}

We now give some examples which illustrate the hypotheses of the theorem.
The following example shows that the hypothesis that the trellis be well-formed is necessary.
\begin{example}
\label{ex:illformed}
The trellis mapping classes in \figref{illformed} are not well-formed.
\fig{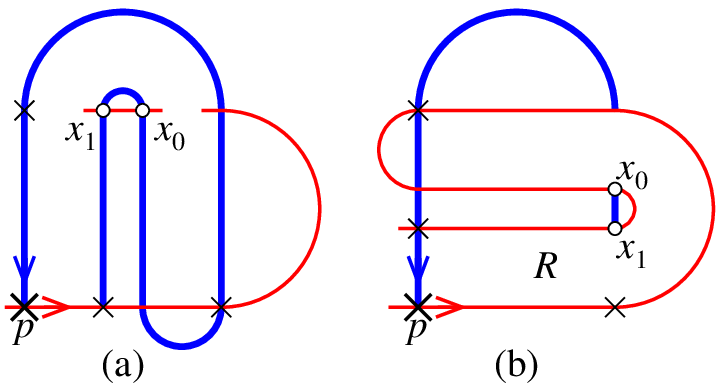}{Two ill-formed trellis mapping classes. The trellis in (a) has Nielsen entropy $\log2$, whereas the trellis in (b), which is the time reversal of that in (a), has Nielsen entropy zero.}
\par
The Nielsen entropy of the trellis mapping class in (a) is equal to $\log2$, so any diffeomorphism in the trellis mapping class must have topological entropy at least $\log2$.
Since the Smale horseshoe map has this trellis type, the topological entropy of the trellis type is exactly $\log2$.
The trellis mapping class of \figrefpart{illformed}{b} is conjugate to the time-reversal of the trellis mapping class in (a).
Since the topological entropy of a diffeomorphism is the same as that of its inverse, any diffeomorphism in this trellis mapping class must have topological entropy at least $\log2$.
However, all the edges of the graph representative are control edges, so the Nielsen entropy is zero.
\end{example}
The above example illustrates that a trellis which is not well-formed may have Nielsen entropy strictly less than the topological entropy, and may even have different Nielsen entropy from its time-reversal.
Even if a trellis mapping class is well-formed, it is not necessarily true that the Nielsen entropy is realised.
A trivial example is of a well-formed trellis type $[f;T]$ for which $\hniel[f;T]=\htop[f;T]$ by for which every diffeomorphism in $[f;T]$ has topological entropy greater than $\hniel[f;T]$ is the planar trellis type with a single transverse homoclinic intersection.
The Nielsen entropy of this trellis type is equal to zero, but every diffeomorphism with a transverse homoclinic point has strictly positive topological entropy.
However, it is simple to construct trellis maps with topological entropy arbitrarily close to zero. 

We now give a nontrivial example.
\begin{example}
\label{ex:nominimalentropy}
\fig{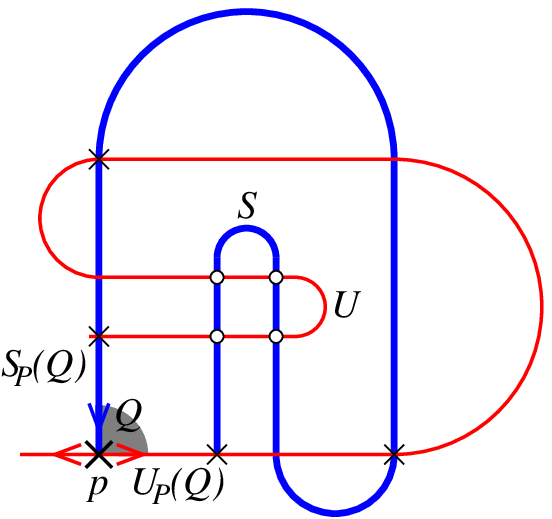}{A trellis for which the Nielsen entropy is not realisable.}
A trellis type $[f;T]$ for which the Nielsen entropy is not realisable is shown in \figref{nominimalentropy}.
Consider the segment $S$.
Taking backward minimal iterates of $S$ eventually yields a segment lying in the region $R(Q)$ with endpoints in the segment $U_P(Q)$.
Then, by the Lambda lemma, under any diffeomorphism $f$ in the trellis mapping class, $f^{-n}(S)$ tends to the closed branch of $T^S(p)$ as $n\tendsto-\infty$, so contains an intersection with $T^U$ for some $n$, even though any \emph{minimal} backward iterate of $S$ has no intersections with $T^U$.
Similarly, $f^n(U)$ must intersect $T^S$ for some $n$ even though any minimal iterate does not.
These observations can be used to show that the Nielsen entropy of the extension $[f;\widehat{T}]$ is greater than that of $[f;T]$, so $\htop(f)>\hniel[f;T]$
In general, the realisability of the entropy bound is closely related to the existence of a diffeomorphism for which every extension is a minimal extension.
\end{example}

In the case where this infemum is realised, we show how to construct a minimal-entropy uniformly-hyperbolic diffeomorphism in the trellis mapping class.
Otherwise, we show, for any $\epsilon>0$, how to construct a diffeomorphism whose entropy is within $\epsilon$ of the Nielsen entropy.
The following example gives a trellis mapping class for which the Nielsen entropy cannot be realised.


\subsection{Nielsen entropy of minimal trellises}
\label{sec:minimalentropy}

Most of the procedures we use to construct diffeomorphisms in a given trellis mapping class rely on extending the original trellis
 and introducing new branches in a controlled way.
The most important type of extension is a minimal extension, since, as noted in Section~\ref{sec:minimalextension}, we expect the Nielsen entropy of a minimal extension to be the same as that of the original trellis mapping class.
We show this in Theorem~\ref{thm:minimalextensionentropy}, using results on the graph representative of a minimal stable supertrellis given in Lemma~\ref{lem:minimalstableiterate}.
We also consider methods for introducing new points of $T^P$ without increasing the Nielsen entropy.
The only difficulty here is on finding the correct initial segment of a branch; once this has been achieved, we can take minimal iterates in the usual way.
Since the Nielsen entropy is computed via a graph representative, which induces an asymmetry between unstable and stable curves, we need to consider each separately.
We consider initial unstable segments in Lemma~\ref{lem:unstablesupertrellis}, and initial unstable segments in Lemma~\ref{lem:stablesupertrellis}.
In Theorem~\ref{thm:minimalsupertrellis} we apply these results to the cases we need to consider later.

\begin{lemma}
\label{lem:minimalstableiterate}
Let $([f];T)$ be an irreducible trellis mapping class with graph representative $(g;G,W)$, and let $([\widetilde{f}];\widetilde{T})$ be a minimal stable supertrellis with graph representative $(\widetilde{g};\widetilde{G},\widetilde{W})$.
Then there is an embedding $i:(G,W)\embed(\widetilde{G},\widetilde{W})$ such that $i\circ g=\widetilde{g}\circ i$, and $i$ restricts to a homeomorphism between the invariant subgraphs $\bigcup_{i=0}^{\infty}g^i(G)$ and $\bigcup_{i=0}^{\infty}\widetilde{g}^i(\widetilde{G})$
Further, if $l$ is a length function on $G$ such that $l(g(z))=l(z)$ for any periodic control edge, $l(g(p))=l(p)$ for any peripheral edge, and $\lambda l(e)\leq l(g(e))<\lambda_\epsilon l(e)$ for any other edge, there is a length function $\widetilde{l}$ on $\widetilde{G}$ with the same properties such that $\widetilde{l}(i(e))=l(e)$ for any edge of $G$.
\end{lemma}

\begin{proof}
First note that $(\widetilde{G};W)$ is exact homotopy equivalence to $(G;W)$, and this exact homotopy equivalence also gives an exact homotopy equivalence between $\widetilde{g}$ considered as a map of $(\widetilde{G};W)$ and $g$.
There is also a natural embedding $W\embed \widetilde{W}$ induced by the inclusion on the stable segments.

Let $\widetilde{\alpha}$ be a path in $\widetilde{G}$ with endpoints only in $W$ and minimal intersections with $\widetilde{W}$.
Since $\widetilde{f}$ is a minimal supertrellis, the path $\widetilde{g}(\alpha)$ also does not back-track.
Since any vertex of $(G,W)$ is contained in the span of the control set $C$, this is enough to prove that $i$ is homotopic to an embedding.
Therefore all control edges of $\widetilde{G}$ are either in the span of $C$, or in a forest which is attached to the span of $C$.
Since all control edges of this forest map into $C$, and have no control edge preimages, all edges of $\widetilde{G}$ map into the span of $C$ under $\widetilde{g}$.

Any edge of $\widetilde{G}$ maps to an edge-path which is the embedding of an edge of $G$ under $\widetilde{g}$.
We take the length function $\widetilde{l}$ in the span of $C$ to be given by $\widetilde{l}(\widetilde{e})=l(g(\widetilde{e}))*(l(e)/l(g(e)))$ where $e$ is the edge of $g$ containing $\widetilde{e}$.
We take the length of an edge $\widetilde{e}$ not in the span of $C$ to be $l(g(\widetilde{e}))/\lambda_*$, where $\lambda<\lambda_*<\lambda_\epsilon$.
This gives the required length function.
\end{proof}
Note that a stable supertrellis which is not minimal may change the topology of the graph.

\begin{theorem}[Nielsen entropy of minimal extensions]
\label{thm:minimalextensionentropy}
Let $([f];T)$ be a well-formed trellis mapping class.
If $([\widetilde{f}];\widetilde{T})$ is a minimal extension of $([f];T)$, then $\hniel[\widetilde{f};\widetilde{T}]=\hniel[f;T]$.
\end{theorem}

\begin{proof}
There exists $n$ such that $\widetilde{f}^{-n}(\widetilde{T}^U\subset T^U$.
Let $([\widehat{f}];\widehat{T})$ be a minimal extension of $([\widetilde{f}];\widetilde{T})$ such that $\widehat{f}^{-n}(\widehat{T}^U)\subset T^U$.
Then $([\widehat{f}];\widehat{T})$ is conjugate to the trellis mapping class $([\widehat{f}];(\widehat{f}^{-n}(T^U);\widehat{f}^{-n}(T^S)))$, which is a minimal stable extension of $([f];T)$, so, by Lemma~\ref{lem:minimalstableiterate}, $\hniel[f;T]=\hniel[\widehat{f};(\widehat{f}^{-n}(T^U);\widehat{f}^{-n}(T^S))]$.
Hence 
\[ \hniel[f;T]\leq\hniel[\widetilde{f};\widetilde{T}]\leq\hniel[\widehat{f};\widehat{T}]
     =\hniel[\widehat{f};(\widehat{f}^{-n}(\widehat{T}^U);\widehat{f}^{-n}(\widehat{T}^S))]=\hniel[f;T] ,\]
yielding $\hniel[\widetilde{f};\widetilde{T}]=\hniel[f;T]$ as required.
\end{proof}

We now give an example of a trellis and a minimal stable iterate, together with the graph representatives.
\begin{example}
\label{ex:minimalextensionentropy}
The minimal extension of the trellis map shown in \figrefpart{minimalextension}{a} is shown in \figrefpart{minimalextension}{b}, together with the graph representatives.
\fig{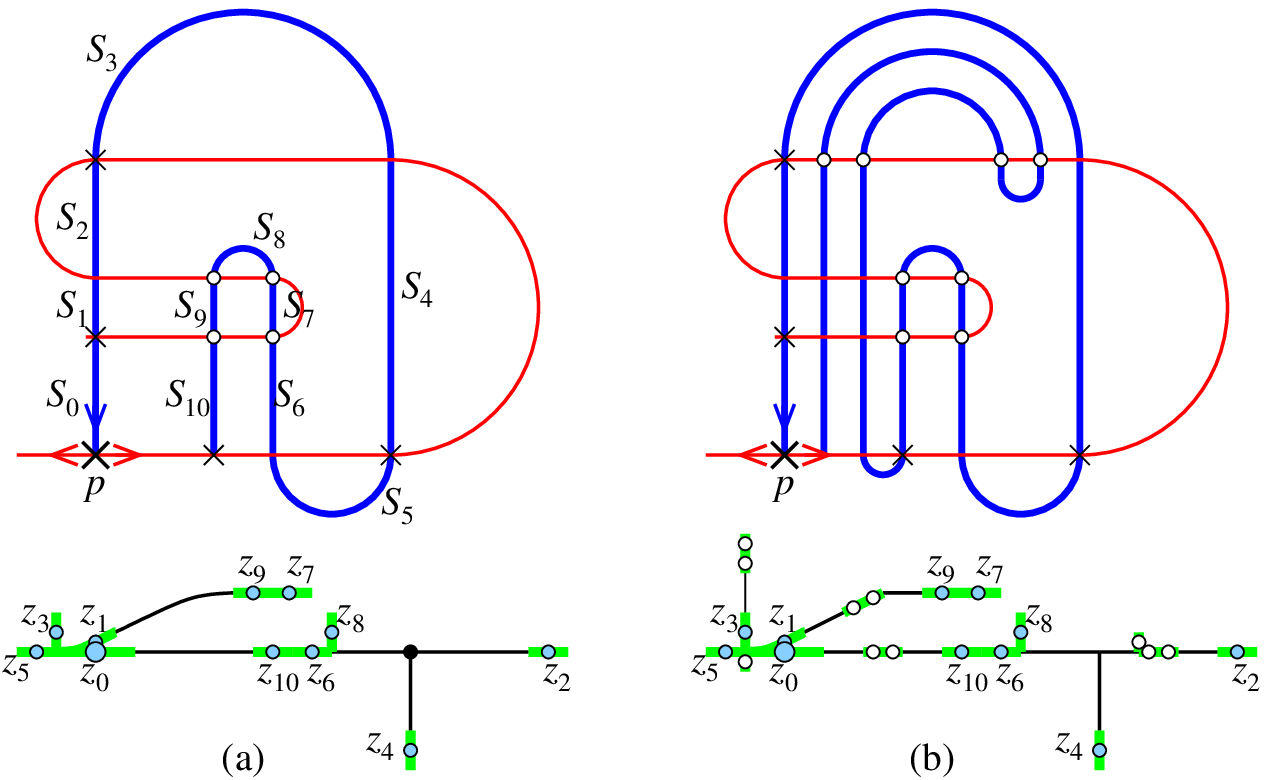}{(a) A trellis type with its graph representative, and (b) the minimal stable iterate.}

The control edges for the original trellis are shown with shaded dots, and those of the extension with white dots.
This clearly shows the embedding of the graph representative of the trellis in that of its minimal extension.
In each case, the essential graph representative is the span of $\{z_0,z_2,z_4\}$.
\end{example}

We next state a simple result which shows that we can puncture at points of an essential periodic orbit without changing the Nielsen entropy.
\begin{lemma}
\label{lem:periodicblowup}
Let $([f];T)$ be a well-formed trellis mapping class with graph representative $(g;G,W)$.
Suppose $P$ is an essential periodic orbit of $f$ which does not globally shadow $T^S$ or $\partial M$.
Let $([\widetilde{f}];\widetilde{T})$ be the trellis mapping class obtained by puncturing $M$ at the points of $P$.
Then $\hniel[\widetilde{f};\widetilde{T}]=\hniel[f;T]$.
\end{lemma}

\begin{proof}
Since $P$ is an essential periodic orbit of $f$, it is shadowed by a periodic orbit of $g$.
Replacing the periodic orbit of $g$ shadowing $P$ with a peripheral loop yields the graph representative $\widetilde{g}$ of $[\widetilde{f};\widetilde{T}]$, which clearly has the same topological entropy as that of $g$.
\end{proof}

We now turn to supertrellises formed by introducing new branches.
The next lemma applies to both stable and unstable supertrellises, since it deals only with minimal iterates of arcs and cross-cuts to $T$.
\begin{lemma}
\label{lem:invariantarcs}
Let $([f];T)$ be a trellis mapping class, and let $A$ be a set of mutually disjoint embedded curves in $\cut[U]{T}$ with endpoints, but no interior points, in $T^S$.
Additionally, if a point of $A$ in an arc $\alpha$ lies in a periodic segment $S$ of $T^S$, then $\alpha$ initially lies along an unstable segment of $T^U$.
Suppose there is a minimal iterate $B$ of $A$ such that every arc of $A$ is contained in an curve of $B$.
Then there is a minimal unstable supertrellis $([\widetilde{f}];\widetilde{T})$ of $([f];T)$ such that $\widetilde{T}^U=T^U\cup B$.

Similarly, if $A$ is a set of cross-cuts to $T^U$ which are disjoint from $T^S$, and there is a minimal iterate $B$ of $A$ such that $A\subset B$, then there is a minimal unstable supertrellis $([\widetilde{f}];\widetilde{T})$ of $([f];T)$ such that $\widetilde{T}^U=T^U\cup B$.
\end{lemma}
\begin{proof}
Let $\beta_{i+1}\in f_{\min}[\alpha_i]$ and be disjoint from $T^U$.
Without loss of generality, we can assume that $\beta_{i+1}\supset \alpha_{i+1}$ for all $i$.
Then the curves $f\circ\alpha_i$ and $\beta_{i+1}$ are isotopic by an isotopy mapping the endpoints into $f(T^S)$.
Further, since the endpoints of curves of $\alpha$ do not lie in segments of $T^S$ containing points of $T^P$, or are points of $T^P$, we can take $h_t$ to be fixed these segments. 
By the isotopy extension theorem~\ref{thm:curveisotopy} this isotopy can be extended to an ambient isotopy $h_t$ which is fixed on $T^U$ and maps $f(T^S)$ into $f(T^S)$.
Let $\widetilde{f}=h_1\circ f$, so $T^U$ and $T^S$ is are still part of the unstable and stable manifolds of $T^P$ for $\widetilde{f}$
Then $\widetilde{f}(T^U)=h_1(f(T^U))\supset h_1(T^U)=T^U$, $\widetilde{f}(T^S)=h_1(f(T^S))=f(T^S)\subset T^S$ and $\widetilde{f}(\alpha_i)=h_1(f(\alpha_i))=\beta_{i+1}\supset \alpha_{i+1}$.
Therefore $([\widetilde{f}];\widetilde{T})$, where $\widetilde{T}=(T^U\cup\bigcup\beta_i,T^S)$, is a well-formed trellis mapping class, and is a minimal supertrellis of $([f];T)$ since the curves $\beta$ are minimal iterates of the curves $\alpha$.
The statement for cross-cuts to $T^U$ follows by reversing time.
\end{proof}

We can use this lemma to show the existence of a minimal unstable supertrellis based on the graph map, and that the Nielsen entropy of this supertrellis is equal to the Nielsen entropy of the original trellis type.
\begin{lemma}
\label{lem:unstablesupertrellis}
Let $([f];T)$ be a well-formed trellis mapping class with graph representative $(g;G,W)$.
Let $\gamma_i$ for $0\leq i<n$  be a set of curves in $G$ with endpoints only in such that $g\circ\gamma_i$ contains $\gamma_{i+1\,\mod\,n}$ as a sub-curve, and each $\gamma_i$ lifts to an embedded curve $\alpha_i$ in $\cut[U]{T}$ such that the $\alpha_i$ are disjoint.
Suppose further that if an endpoint of some $\gamma_i$ lies in a periodic point $w\in W$, then $\gamma_i$ lies on the same side of $W$ as a trivial branch of $T$.
Then the collection $A=\{\alpha_i\}$ satisfies the conditions of Lemma~\ref{lem:invariantarcs}, and the minimal supertrellis $([\widetilde{f}];\widetilde{T})$ so obtained has $\hniel[\widetilde{f};\widetilde{T}]=\hniel[f;T]$.
\end{lemma}
\begin{proof}
The curves $f\circ\alpha_i$ are homotopic to $\epsilon_{i+1}=g\circ\gamma_i$, so are homotopic to minimal iterates $\beta_{i+1}$ of the $\alpha_i$.
By a further homotopy we can ensure $\alpha_i$ is contained in $\beta_i$, and so satisfies the hypotheses of Lemma~\ref{lem:invariantarcs}.
The graph representative $(\widetilde{g};\widetilde{G},\widetilde{W})$ is given by cutting along the curves $\epsilon_i$, so there is a projection $\pi:(\widetilde{G},\widetilde{W})\exto(G,W)$ such that $\pi\circ\widetilde{g}=g\circ\pi$.
Since $\widetilde{g}$ is optimal, this ensures $\htop(\widetilde{g})=\htop(g)$.
\end{proof}

We have already shown that the Nielsen entropy of a minimal stable supertrellis is the same as that of the original trellis type.
The following result shows how to construct minimal stable supertrellises.
\begin{lemma}
\label{lem:stablesupertrellis}
Let $([f];T)$ be a well-formed trellis mapping class with graph representative $(g;G,W)$.
Suppose $P$ is an essential periodic orbit of $f$ which does not globally shadow $T^S$ or $\partial M$.
Then $f$ is isotopic relative to $P$ to a diffeomorphism $\widetilde{f}$ with a trellis $\widetilde{T}$ such that $([\widetilde{f}];\widetilde{T})$ is a well-formed supertrellis of $([f];T)$ such that $P$ is contained in an attractor of $([\widetilde{f}];\widetilde{T})$ and $\hniel[\widetilde{f};\widetilde{T}]=\hniel[f;T]$.
\end{lemma}

\begin{proof}
Take cross-cuts $A$ in the statement of Lemma~\ref{lem:invariantarcs} to have essential intersections with the edges of $g$ containing $P$, giving trellises and controlled graphs shown in \figref{stablesupergraph}.
\fig{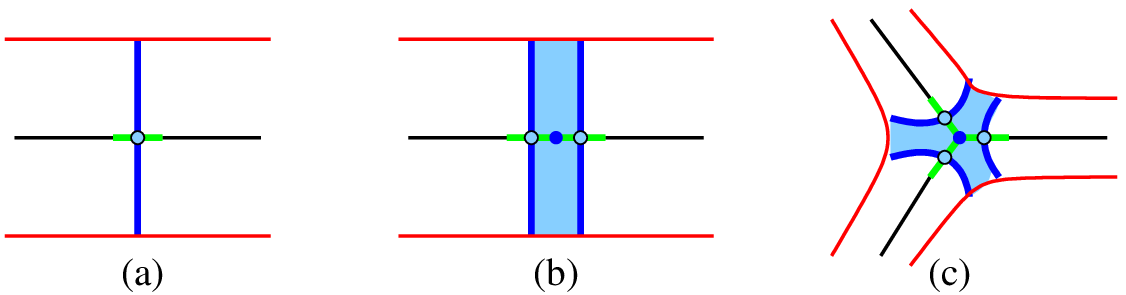}{Trellises and graphs formed by introducing new stable curves edges. The in (a) formed by introducing a single stable curves does not have any attractors. The trellises in (b) and (c) formed by introducing one stable curve crossing each incident edge at a vertex do yield new attracting regions.}
Then the supertrellis $([\widetilde{f}];\widetilde{T})$ given by Lemma~\ref{lem:invariantarcs} is the required trellis mapping class.
\end{proof}

\fig{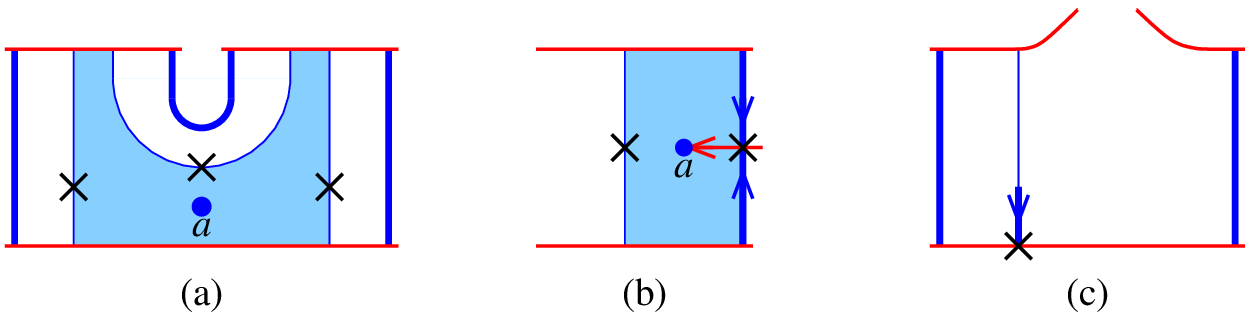}{Stable supertrellises. In (a) we introduce an attractor, in (b) we bound an attractor, and in (c) we extend a trivial branch.}
We use these results to construct new branches at essential periodic orbits and trivial branches of a trellis.
These cases are shown in \figref{stablesupertrellis}.
These branches are may used to create new attractors and repellors, or subdivide a region into an attractor/repellor and a collection of rectangles, as in \figrefpart{stablesupertrellis}{a}, to bound the end of a trivial branch in an attractor or repellor, as in (b), or to extend a trivial branch to a nontrivial branch, as in (c).
These cases correspond to the three cases of Theorem~\ref{thm:minimalsupertrellis}
\begin{theorem}
\label{thm:minimalsupertrellis}
Let $([f];T)$ be a well-formed trellis mapping class.
Then:
\begin{enumerate}
\renewcommand{\theenumi}{(\alph{enumi})}
\item\label{item:periodicblowup}
 If $P$ is an essential periodic orbit of $([f];T)$ which does not shadow $T^P$, there is are minimal supertrellises of $([f];T)$ for which $P$ is contained in a region which is attracting, repelling, or both.
\item\label{item:periodictrivialend} If $P$ is a periodic orbit at the ends of trivial branches of $T^P$, then there is a minimal supertrellis $([\widetilde{f}];\widetilde{T})$ of $([f];T)$ for which $P$ is contained in a non-chaotic region of $([\widetilde{f}];\widetilde{T})$.
\item\label{item:trivialbranch} If $T^{\us}(p,b)$ is a trivial branch of $T$ which lies in a chaotic region of $T$, then there is a minimal supertrellis $([\widetilde{f}];\widetilde{T})$ of $([f];T)$ for which $\widetilde{T}^{\us}(p,b)$ is contained in a nontrivial branch of $T^P$.
\end{enumerate}
\end{theorem}
All these statements follow from Lemma~\ref{lem:unstablesupertrellis} and  Lemma~\ref{lem:stablesupertrellis} by constructing the required curves or finding the right control edges.


\subsection{Regular domains}

Any diffeomorphism is conjugate to a linear hyperbolic diffeomorphism in a neighbourhood of a point of a periodic saddle point.
We therefore expect a trellis map $([f];T)$ to behave in a fairly predictable way in a neighbourhood of a point of $T^P$.
In particular, we know how a rectangle with sides parallel to the local stable and unstable foliations behaves.
Unfortunately, for a given trellis mapping class, the neighbourhood of $T^P$ on which we have hyperbolic behaviour may be arbitrarily small.
To deal with this problem, we introduce the concept of a \emph{regular domain}, which is a rectangular domain which behaves similarly to a sufficiently small rectangular neighbourhood of a periodic saddle point.

\begin{definition}[Regular domain]
Let $([f];T)$ be a trellis mapping class.
A rectangular domain $D$ of $T$ is a \emph{regular domain} for a period-$n$ quadrant $Q$ of $T$ if $D$ has sides $T^U[p,q^u]\subset f^{-n}(T^U)$, $T^S[p,q^s]\subset f^n(T^S)$, $T^U[q^s,r]$ and $T^S[q^u,r]$, such that $f^n(T^U[p,q^u])\cap D=T^U[p,q^u]$, $f^{-n}(T^S[p,q^s])\cap D=T^S[p,q^s]$, and $f^{-n}(T^U[q^s,r])\cap D=f^n(T^S[q^u,r])\cap D=\emptyset$.
\end{definition}
The sides $T^U[p,q^u]$ and $T^S[p,q^s]$ are called \emph{adjacent sides} of $D$, and the sides $T^U[q^s,r]$ and $T^S[q^u,r]$ are \emph{opposite sides}.
Note that the definition is given purely in terms of the topology of the trellis and the mapping of its vertices; no extensions are needed.

\fig{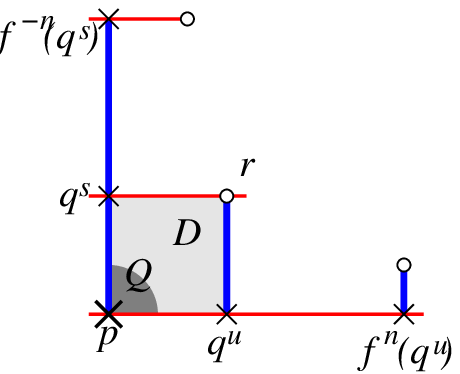}{The domain $D$ is a regular domain for $Q$.}
A regular domain is shown in \figref{regularquadrant}.
We shall always denote the vertices of a regular domain by $p$, $q^u$, $q^s$ and $r$ as shown.
We now give a number of elementary properties of a regular domain.
The first two lemmas are trivial.
\begin{lemma}
\label{lem:regulardomainimage}
Let $([\widehat{f}];\widehat{T})$ be an extension of $([f];T)$.
Then if $D$ is a regular domain for a quadrant $Q$ in $([f];T)$, then $D$ is also a regular domain for $Q$ in $([\widetilde{f}];\widetilde{T})$.
\end{lemma}
\begin{lemma}
Let $D$ be a regular domain of $Q$, and suppose that both $f^{-(n+1)}(T^S[p;q^s])$ and $f^{-1}(T^S[q^u,r])$ are subsets of $T^S$.
Then $f^{-1}(D)$ is a regular domain for $f^{-1}(Q)$.
\end{lemma}

\begin{lemma}
\label{lem:regularsubdomain}
Let $D$ be a regular domain for $Q$, and $T^S[\tilde{q}^u,\tilde{r}]\subset D$ be a curve with endpoints $\tilde{q}^u\in T^U(p,q^u)$ and $\tilde{r}\in T^U(q^s,r)$.
Then the rectangular domain $\widetilde{D}$ with vertices at $\{p,\tilde{q}^u,q^s,\tilde{r}\}$ is a regular domain for $Q$.
Similarly, if $T^U[\tilde{q}^s,\tilde{r}]\subset D$ is a curve with endpoints $\tilde{q}^s\in T^S(p,q^s)$ and $\tilde{r}\in T^S(q^u,r)$, then the rectangular domain vertices at $\{p,q^u,\tilde{q}^s,\tilde{r}\}$ is a regular domain for $Q$.
\end{lemma}

\begin{proof}
The only nontrivial step is to show $f^n(T^S[\tilde{q}^u,\tilde{r}^u])\cap\widetilde{D}=\emptyset$.
First, note that $f^n(\tilde{q}^u)\in T^U(\tilde{q}^u,f^n(q^u))$, so $f^n(\tilde{q}^u)\not\in\widetilde{D}$.
Further, $f^n(T^S[\tilde{q}^u,\tilde{r}^u])\cap T^U[p,\tilde{q}^u]=f^n(T^S[\tilde{q}^u,\tilde{r}^u]\cap f^{-n}(T^U[p,\tilde{q}^u]))=f^n(\emptyset)=\emptyset$
 and $f^n(T^S[\tilde{q}^u,\tilde{r}^u])\cap T^U[p,\tilde{q}^u]=f^n(T^S[\tilde{q}^u,\tilde{r}^u]\cap f^{-n}(T^U[p,\tilde{q}^u]))\subset f^n(D\cap f^{-n}(T^U[p,\tilde{q}^u))=\emptyset$, so $f^n(T^S[\tilde{q}^u,\tilde{r}^u])$ does not intersect the unstable boundary of $\widetilde{D}$.
Additionally, $f^n(T^S[\tilde{q}^u,\tilde{r}])\cap T^S[p,q^s]=f^n(T^s[\tilde{q}^u,\tilde{r}]\cap f^{-n}(T^S[p,q^s]))\subset f^n(T^S[\tilde{q}^u,\tilde{r}]\cap T^S[p,q^s])=\emptyset$.
It remains to show that $f^n(T^S[\tilde{q}^u,\tilde{r}])$ does not intersect itself.
We have already seen that $f^n(T^S[\tilde{q}^u,\tilde{r}])$ does not contain $\tilde{q}^u$ or $\tilde{r}$.
Further, $f^n(\{\tilde{q}^u\})\cap T^S[\tilde{q}^u,\tilde{r}]\subset f^n(T^U(p,q^u))\cap T^S[\tilde{q}^u,\tilde{r}]=\{\tilde{q}\}$, but clearly $f^n(\tilde{q}^u)\neq\tilde{q}^u)$, so $f^n(\tilde{q}^u)\not\in T^S[\tilde{q}^u,\tilde{r}]$.
Since $T^S[\tilde{q}^u,\tilde{r}]$ is an interval, this is enough to show that $f^n(T^S[\tilde{q}^u,\tilde{r}])\cap T^S[\tilde{q}^u,\tilde{r}]$.
Hence $f^n(T^S[\tilde{q}^u,\tilde{r}])\cap \partial\widetilde{D}=\emptyset$, so $f^n(T^S[\tilde{q}^u,\tilde{r}])\cap \widetilde{D}=\emptyset$.
A similar analysis proves the statement for a curve $T^U[\tilde{q}^s,\tilde{r}]$.
\end{proof}

\fig{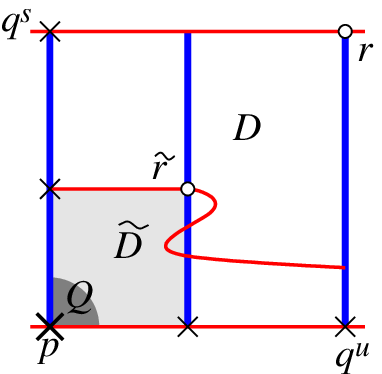}{$\widetilde{D}$ is a regular domain if $D$ is a regular domain.}
Two applications of Lemma~\ref{lem:regularsubdomain} show that the region $\widetilde{D}\subset D$ shown in \figref{regularsubdomain} is a regular domain for $Q$ if $D$ is a regular domain.

One of the most important properties of a regular domain is that iterates of the sides opposite $Q$ must cross each other.
\begin{lemma}
\label{lem:regulardomainiterate}
Suppose $D$ is a regular domain for a period-$n$ quadrant $Q$, and suppose $f^{-n}(T^S[q^u,r])\subset T^S$.
Then $f^{-n}(T^S[q^u,r])$ contains a subinterval of the form $T^S[\tilde{q}^u,\tilde{r}]\subset D$ such that $\tilde{q}^u=f^{-n}(q^u)\in T^U(p,q^u)$ and $\tilde{r}\in T^S(q^s,r)$.
\end{lemma}

\begin{proof}
Since $f^n(Q)=Q$, $f^n$ is orientation-preserving, so the orientation of the intersection at $f^{-n}(q^u)$ is the same as that at $q^u$.
Further, $f^{-n}(T^S[q^u,r])\cap T^U[p,q^u]=f^{-n}(T^S[q^u,r]\cap f^n(T^U[p,q^u]))=f^{-n}(T^S[q^u,r]\cap T^U[p,q^u])=\{f^{-n}(q^u)\}$, so $f^{-n}(T^S(q^s,r])$ does not cross $T^U[p,q^u]$.
Finally, $\{f^{-n}(r)\}\cap D\subset f^{-n}(T^U[q^s,r])\cap D=\emptyset$, so $f^{-n}(r)\not\in D$.
Therefore, $f^{-n}(T^S[q^u,r])$ intersects $T^U[q^s,r]$, and we let $\tilde{r}$ be this first intersection.
\end{proof}
Notice that Lemma~\ref{lem:regulardomainiterate} is enough to show that if $D$ is a regular domain for a trellis mapping class $([f];T)$, and $([\widetilde{f}];\widetilde{T})$ is an extension of $([f];T)$ with $\widehat{f}^{-n}(T^S[q^u,r])\subset \widehat{T}^S$, then the biasymptotic orbit through $\tilde{r}$ is forced by $([f];T)$.

\fig{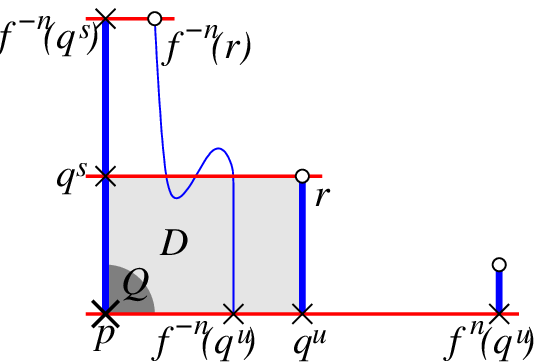}{The backward iterate of $T^S[q^u,r]$ intersects $T^U[q^s,r]$ but does not cross $T^U[p,q^u]$.}
In \figref{regularsubdomain} we show a regular domain $D$ together with $f^{-n}(T^S[q^u,r])$.
Notice that $f^{-n}(T^S[q^u,r])$ may cross $T^U[q^s,r]$ several times, but must cross at least once since $f^{-n}(r)\not\in D$.

The most important case of a regular domain for a quadrant $Q$ is when $D$ is a region, in which case we have a \emph{regular region} for $Q$.
It is also useful to consider the less restrictive case for which the interior of $D$ does not intersect $T^U$, or in other words, if $D\cap T^U\subset\partial D$.
In this case, there is a unique homotopy class of exact curves $\alpha_D:(I,J)\exto(\cut[U]{T},T^S)$ such that $\alpha_D(I)\subset D$, $\alpha_D(0)\in T^U[p,q^s]$, $\alpha_D(1)\in T^U[q^u,r]$, and $\alpha_D$ has minimal intersections with $T^S$.
It is an important but trivial observation that if $D\cap T^U=\partial D^U$ and if $\widehat{T}$ is an extension of $T$ with $\widehat{T}^U=T^U$, then $D\cap\widehat{T}^U=\partial D^U$.
The following lemma shows that any curve in such a regular domain with initial point in $T^S[p,q_s]$ iterates to a curve which crosses $T^S[q^u,r]$.
\begin{lemma}
Let $D$ be a regular domain for a quadrant $Q$ such that $D\cap T^U=\partial D^U$ and $D$ contains no smaller regular domains, and $D\cap T^P=\{p\}$.
Then for any exact curve $\alpha:(I,\{0,1\})\fto(\cut[U]{T},T^S)$ with $\alpha(I)\subset D$, $\alpha(0)\in T^S[p,q^s]$ and $\alpha(1)\not\in T^S[p,q^s]$, there exists $k$ such that $f_{\min}^{nk}[\alpha]$ crosses $T^S[q^u,r]$.
\end{lemma}

\begin{proof}
Let $S$ be the segment containing the final endpoint of $\alpha$.
Suppose $f^{in}(S)\subset D$ for all $i$.
Then $f^{in}(S)\cap T^S[p,q^s]=\emptyset$ for all $i$, and since there are only finitely many stable segments in $D$, one must map into itself under $f^{jn}$ for some $j$, giving a periodic point in $D$ distinct from $p$, a contradiction.
Therefore $f^{kn}(S)\not\subset D$ for some least $k$
For this $k$, we must have $f_{\min}^{kn}[\alpha]\cap T^S[q^u,r]\neq\emptyset$, since crossing $T^S[p,q^s]$ would imply a point of $f^{kn}\circ\alpha$ in $T^S(q^s,f^{-n}(q^s))$, a contradiction.
\end{proof}

Recall the meaning of a curve $\widetilde{\alpha}$ tightening onto $\alpha$ given in Definition~\ref{defn:curvetightening}.
We show that if $D$ and $\widetilde{D}$ are regular domains for a quadrant, then some iterate of $\alpha_D$ tightens onto $\alpha_{\widetilde{D}}$.
The importance of this result is that if a curve $\widehat{\alpha}$ tightens onto $\alpha$, and both curves are tight curves embedded in a graph, then $\alpha(I)\subset\widehat{\alpha}(I)$.

\begin{lemma}
\label{lem:subdomaintighten}
Let $D$ be a regular domain for a quadrant $Q$ such that $D\cap T^U=\partial D^U$, and let $\widetilde{D}\subset D$ be a regular domain.
Then there exists $k$ such that $f_{\min}^{kn}[\alpha_{D}]$ tightens onto $[\alpha_{\widetilde{D}}]$.
\end{lemma}

\begin{proof}
Let $k$ be the least integer such that $f^{kn}(\tilde{q}^u)\not\in T^U[p,q^u]$.
Then $f^{kn}(\tilde{q}^u)\not\in D$, so the endpoint of $f_{\min}^{kn}[\alpha_{\widetilde{D}}]$ does not lie in $D$.
Further, the first intersection of $f_{\min}^{kn}[\alpha_{\widetilde{D}}]$ with $\partial D^S$ must be with $T^S[q^u,r]$, so $f_{\min}^{kn}[\alpha_{\widetilde{D}}]$ tightens onto $[\alpha_D]$, which is the curve from $T^S[p,q^s]$ to $T^S[q^u,r]$ in $D$ with minimal intersections with $T^S$.
\end{proof}

\fig{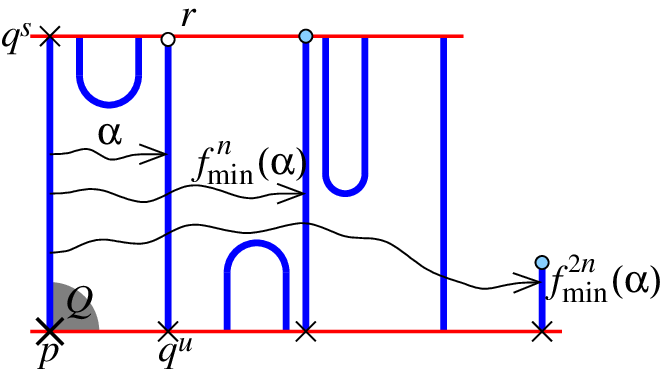}{A curve $\alpha_D$ in a regular domain $D$ and some of its minimal iterates.}
An example of the minimal iterates of a curve $\alpha=\alpha_D$ for a regular domain $D$ is shown in \figref{alphacover}.

\begin{lemma}
\label{lem:subdomainiteratetighten}
Let $D$ be a regular domain for a quadrant $Q$, and $\widetilde{D}\subset D$ a regular domain for $Q$ such that $\widetilde{D}\cap T^U=\partial\widetilde{D}^U$.
Let $\alpha$ be a curve in $D$ such that $\alpha(0)\in T^S[p,q^s]$ and $\alpha(1)\in T^S[q^u,r]$.
Then there exists $k$ such that $f_{\min}^{kn}[\alpha]$ tightens onto $[\alpha_{\widetilde{D}}]$.
\end{lemma}

\begin{proof}
For any $i$, $f_{\min}^{in}[\alpha]$ must cross $T^S[q^u,r]$, and do so before it crosses $T^S[p,q^s]$.
Choose $k$ such that $f_{\min}^{kn}(\alpha(0))\in T^S[p,\tilde{q}^s]$.
Then $f_{\min}^{kn}[\alpha]$ tightens onto $[\alpha_{\widetilde{D}}]$.
\end{proof}

If $Q$ is a quadrant with a regular domain $D$ such that $D\cap T^U=D^U$, we let $D(Q)$ be the smallest regular domain for $Q$, and write $[\alpha_Q]$ for $[\alpha_{D(Q)}]$.
Since for any regular domain $D$, the $n$th iterate of $\alpha_D$ tightens onto itself, there is some recurrent behaviour for the graph representative.
If such a curve than has some iterate which tightens onto $\alpha_D$ for a different quadrant, the graph representative maps edges from one quadrant to another.
This provides the definition of an alpha-chain.
\begin{definition}[Alpha chain]
Let $([f];T)$ be a well-formed trellis mapping class and $Q_U$ and $Q_S$ be quadrants of $T$ which are contained in regular domains which do not intersect $T^U$ in their interiors.
We say there is an \emph{alpha chain} from $Q_U$ to $Q_S$ if there exists $n$ such that $f_{\min}^n[\alpha_{Q_U}]$ tightens onto $[\alpha_{Q_S}]$.
\end{definition}

\fig{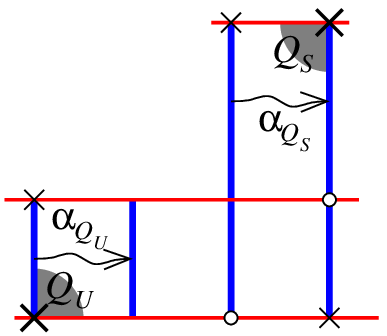}{A trellis with an alpha chain from $Q_U$ to $Q_S$.}
An example of an alpha chain from a quadrant $Q_U$ to $Q_S$ is shown in \figref{alphachain}.
If there is an alpha chain from any segment to any other, the trellis is said to be \emph{transitive}.
\begin{definition}[Transitive trellis]
We say that $T$ is \emph{transitive} if every quadrant $Q$ is contained in a regular domain regular $D_Q$ such that $D_Q\cap T^U=\partial D_Q$, and for every pair of quadrants $Q_U$ and $Q_S$, there is an alpha chain from $Q_U$ to $Q_S$.
\end{definition}

The following lemma gives a sufficient condition for the existence of an alpha-chain from$Q_U$ to $Q_S$.
\begin{lemma}
Suppose $D_U$ and $D_S$ are regular domains for $Q_U$ and $Q_S$ respectively such that $D_U\cap D^S$ is a rectangular domain with unstable edges contained in $\partial D_U^U$ and stable edges in $\partial D_S^S$.
Then there is an alpha chain from $Q_U$ to $Q_S$.
\end{lemma}

\begin{proof}
Let $\beta$ be a curve in $D^U\cap D^S$ homotopic to the initial interval of $T^U(p_u)\cap D^S$.
Then $f_{\min}^{k_un_u}[\alpha_{Q_U}]$ tightens onto $\beta$ for some $k_u$, and since $\beta$ is a curve joining the stable sides of $D_S$, we also have $f_{\min}^{k_sn_s}[\beta]$ tightens onto $\alpha_{Q_S}$ for some $k_s$.
Hence $f_{\min}^{k_un_u+k_sn_s}[\alpha_{Q_U}]$ tightens onto $[\alpha_{Q_S}]$ as required.
\end{proof}

Transitivity is preserved by taking stable extensions.
\begin{lemma}
If there is an alpha chain from $Q_U$ to $Q_S$ for $([f];T)$ and $(l\widehat{f}];\widehat{T})$ is a stable extension of $T$, then there is an alpha chain from $Q_U$ to $Q_S$ for $([\widehat{f}];\widehat{T})$.
In particular, if $([f];T)$ is transitive, and $(l\widehat{f}];\widehat{T})$ is a stable extension of $T$, then $([\widehat{f}];\widehat{T})$ is transitive.
\end{lemma}

\begin{proof}
By Lemma~\ref{lem:subdomaintighten}, $\widehat{f}_{\min}^{k_un_u}[\widehat{\alpha}_{Q_U}]$ tightens onto $[\alpha_{Q_U}]$ for some $k_u$, and by Lemma~\ref{lem:subdomainiteratetighten}, $\widehat{f}_{\min}^{k_sn_s}[\alpha_{Q_S}]$ tightens onto $[\alpha_{Q_S}]$ for some $k_s$.
Further, $f_{\min}^{n}[\alpha_{Q_U}]$ tightens onto $[\alpha_{Q_S}]$ for some $n$ since there is an alpha-chain from $Q_S$ to $Q_S$ in $T$.
Therefore, there is an alpha chain from $Q_U$ to $Q_S$ in $\widetilde{T}$.
\end{proof}

We now show that a transitive trellis has a transitive graph representative.
\begin{lemma}
\label{lem:transitivegraph}
Let $([f];T)$ be an transitive trellis mapping class.
Then the graph representative $(g;G,W)$ of $([f];T)$ has a single transitive component with positive topological entropy.
\end{lemma}

\begin{proof}
Let $\overline{G}=\bigcup_{n=0}^{\infty}g^n(G)$ be the essential subgraph, and let $e$ be any expanding edge of $G$.
By irreducibility, there is a quadrant $Q$ such that $g^n(e)$ tightens to $\alpha_{Q}$ for some $n$.

We now consider preimages of edges.
Again, let $e$ be an expanding edge of $G$.
There exists a tight curve $\beta_0$ in $\widetilde{G}$ with endpoints in $W$ such that $\beta_0$ is homotopic to a subinterval of a branch $T^U(p,b)$.
We then find a curve $\beta_1$ in $\widetilde{G}$ homotopic to a subinterval of a branch at $f^{-1}(p)$.
Proceeding recursively gives an edge-path $\beta_n=\alpha_Q$ where $\beta_i$ is a sub-path of $g(\beta_{i+1})$ for $0\leq i<n$.
Therefore, there exists a quadrant $Q$ and an integer $n$ such that $g^n(\alpha_Q)\supset e$. 

Now if $Q_U$ and $Q_S$ are any two quadrants, there exists $n$ such that $g^n(\alpha_{Q_U})$ contains a sub-curve $\alpha_{Q_S}$.
Since $(f;T)$ is irreducible, we can find $N$ such that for any two quadrants and any $n\geq N$, the iterate $g^n(\alpha_{Q_U})$ contains a sub-path $\alpha_{Q_S}$.
Therefore, any such curve $\alpha$ generates the same graph component under iteration.

Combining these results shows that there exists $N$ such that if $e_1$ and $e_2$ are any two edges, and $n\geq N$, then $e_2\subset g^n(e_1)$.
\end{proof}

Finally, we show that the hyperbolicity near $T^P$ is enough to create intersections from which we can deduce regularity and transitivity.
Note that for this result we are concerned with $f$-extensions of a trellis map $(f;T)$ rather than an isotopy class.
\begin{lemma}
\label{lem:regularextension}
Let $(f;T)$ be a trellis map and $Q$ a quadrant of $T$.
Let $q\in T^U(Q)$ be the endpoint of an unstable segment $S$ on the $Q$-side of $T^S(Q)$, and let $q^s\in T^S(Q)$ be the endpoint of a stable segment $U$ on the $Q$-side of $T^S(Q)$.
Then there exists $k$ such that $f^{-kn}(S)$ intersects $U$ at a point $r$ such that $\{p,f^{-kn}(q),q^s,r\}$ are the vertices of a regular domain for $Q$.
\end{lemma}

\begin{proof}
By the Lambda lemma, as $i\tendsto\infty$, $f^{-in}(S)$ limits on $W^S(Q)$ in the $C^1$ topology.
Take a neighbourhood $K$ of $q^s$ such that $f^{-n}(K)\cap T^S[p,q^s]=\emptyset$, and choose $k$ such that $f^{-kn}(S)$ intersects $U$ in at a point $r$ in $K$ such that $T^S[q^s,r]\subset K$, and the domain $D$ with vertices at $\{p,f^{-kn}(q),q^s,r\}$ is a rectangle which does not intersect $f^{-n}(K)$, and such that $f^n(T^U(p,f^{-kn}(q^u)))\cap D=T^U[p,f^{-kn}(q^u)]$.
Let $\widehat{T}=(T^U,f^{-nk}(T^S))$.
Then $f^{-n}(T^U[q^s,r])\cap D\subset f^{-n}(K)\cap D=\emptyset$, so $D$ is a regular domain of $(f;\widehat{T})$, as shown in \figref{regularisequadrant}.
\fig{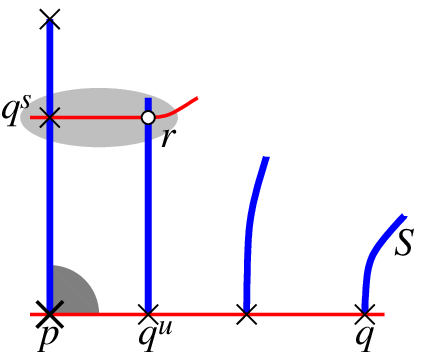}{Backward iterates of $S$ eventually form the opposite stable side of a regular domain.}
\end{proof}


\subsection{Existence of entropy minimisers}

To prove the existence of a diffeomorphism in a trellis mapping class whose topological entropy is the Nielsen entropy of the class, we first show how to construct such a diffeomorphism for a particularly simple class of trellises.
\begin{lemma}
\label{lem:simpleminimalmodel}
Let $([f];T)$ be a trellis mapping class.
Suppose every chaotic region $R$ of $([f];T)$ is a rectangle, and that every such rectangle is contained in a rectangular domain $D$ with $D^U\subset f^{-1}(T^U)$ such that every region in $D$ and $f(D)$ is a rectangle.
Then there exists a diffeomorphism $\widehat{f}\in([f];T)$ such that $\htop(\widehat{f})=\htop([f];T)$.
\end{lemma}

\begin{proof}
For every domain $D$ containing a chaotic region $R$, foliate $D$ and $f(D)$ by an unstable foliation $\cal{F}^U$ parallel to $T^U$ and a transverse stable foliation $\cal{F}^S$ parallel to $T^S$.
Isotope $f$ to obtain a diffeomorphism $\widetilde{f}$ which preserves the stable and unstable foliations,
 and for which all points of non-chaotic regions are in the basin of a stable or unstable periodic orbit.
Let $(G,W)$ be a graph representative of $([f];T)$ for which $G$ is transverse to $T^S$, and $\pi:(\cut[U]{T},T^S)\exto(G,W)$ be a deformation-retract which collapses each leaf of $\cal{F}^S$ onto a point of $g$.
Then $\pi\circ \widetilde{f}=g\circ\pi$ on every chaotic region $R$, so $\htop(\widetilde{f})=\htop(g)=\hniel[f;T]$.
\end{proof}

We can use this to prove the general case.
\begin{theorem}[Existence of entropy minimisers]
\label{thm:minimalmodel}
Let $([f];T)$ be a trellis mapping class.
Suppose there is a diffeomorphism $\widehat{f}\in([f];T)$ such that every extension of $T$ by $\widehat{f}$ is minimal.
Then there is a uniformly-hyperbolic diffeomorphism $\widetilde{f}\in([f];T)$ such that every extension of $T$ by $\widetilde{f}$ is minimal, and $\htop(\widetilde{f})=\hniel[f;T]$.
\end{theorem}

\begin{proof}
By Theorem~\ref{thm:minimalextensionentropy}, any trellis mapping class of the form $[\widehat{f};(T^U,\widehat{f}^{-n}(T^S)]$ has the same Nielsen entropy at $([f];T)$.
By Lemma~\ref{lem:regularextension} we can therefore take a $\widehat{f}$-extension $\widehat{T}_1$ of $T$ such that every quadrant of $\widehat{T}_1$ lies in a rectangular region, 
 and by irreducibility, we can ensure that every region of  $\widehat{T}_1$ is a topological disc or annulus.
If $r_0$ is a point of $\widehat{T}_1^S$ at the end of a trivial branch and lies in a chaotic region, by introducing new unstable curves as in Theorem~\ref{thm:minimalsupertrellis}, we can take a minimal supertrellis $([\widetilde{f}_2];\widetilde{T}_2)$ such that the orbit of $r_0$ lies in a repelling or chaotic region.
Similarly, if $a_0$ is a point of $\widetilde{T}_2^S$ at the end of a trivial branch of $\widehat{T}_1^U$, by introducing new stable curves as in Theorem~\ref{thm:minimalsupertrellis}, we can take a minimal supertrellis $([\widetilde{f}_3];\widetilde{T}_3)$ such that the orbit of $a_0$ lies in a attracting or chaotic region.

Now suppose there is a chaotic region $R$ of $([\widetilde{f}_3];\widetilde{T}_3)$ which is not a rectangle.
Then the graph representative $(\widetilde{g}_3];\widetilde{G}_3,\widetilde{W}_3)$ of $([\widetilde{f}_3];\widetilde{T}_3)$ has a peripheral loop or a valence-$n$ vertex in $R$ which corresponds to a boundary component or essential periodic orbit which does not shadow $T^S$.
Introducing new stable curves for all such $R$ as in Theorem~\ref{thm:minimalsupertrellis} gives a minimal supertrellis $([\widetilde{f}_4];\widetilde{T}_4)$ for which $R$ is a domain containing an attractor and some chaotic rectangles.

Every chaotic region $R$ of $([\widetilde{f}_4];\widetilde{T}_4)$ is now a rectangle. 
Taking a minimal backwards stable iterate $([\widetilde{f}_5];\widetilde{T}_5)$ of $([\widetilde{f}_4];\widetilde{T}_4)$ gives a trellis mapping class satisfying the conditions of Lemma~\ref{lem:simpleminimalmodel}.
Hence there is a diffeomorphism $\widetilde{f}\in([\widetilde{f}_5];\widetilde{T}_5)$ for which $\htop(\widetilde{f})=\hniel[\widetilde{f}_5;\widetilde{T}_5]$, which equals $\hniel[f;T]$ since every extension and supertrellis is minimal.
Further, it is clear that every $\widetilde{f}$-extension of $\widetilde{T}_5$ is minimal, so every $\widetilde{f}$-extension of $T$ is minimal.
\end{proof}


\subsection{Approximate entropy minimising diffeomorphisms}
\label{sec:minimisingsequence}

Throughout this section we take $([f];T)$ to be a well-formed irreducible trellis mapping class, $\lambda=\exp(\hniel[f;T])$ and $\lambda_\epsilon=\exp(\hniel[f;T]+\epsilon)$.
We aim to show that we can always find a diffeomorphism $\widehat{f}\in([f];T)$ which has entropy less than $\lambda_\epsilon$.
We first consider how to perform non-minimal extensions without increasing the Nielsen entropy above $\hniel[f;T]+\epsilon$.
The technical details are given in Lemma~\ref{lem:stableiteratelength} and Lemma~\ref{lem:entropyincreasingextension}, which deal directly with the graph representatives.
We then prove Lemma~\ref{lem:transitiveextensionincrease}, which allows us to consider non-minimal extensions without reference to the graph representative.

Our aim is to to construct a trellis mapping class which satisfies the conditions of Theorem~\ref{thm:minimalmodel}.
To do this, we may need to introduce new periodic points to the trellis to create attracting and repelling regions.
We then iterate curves bounding an attractor or repellor into regular domains, and finally take non-minimal iterates to move bigon boundaries into attractors and repellors.

\fig{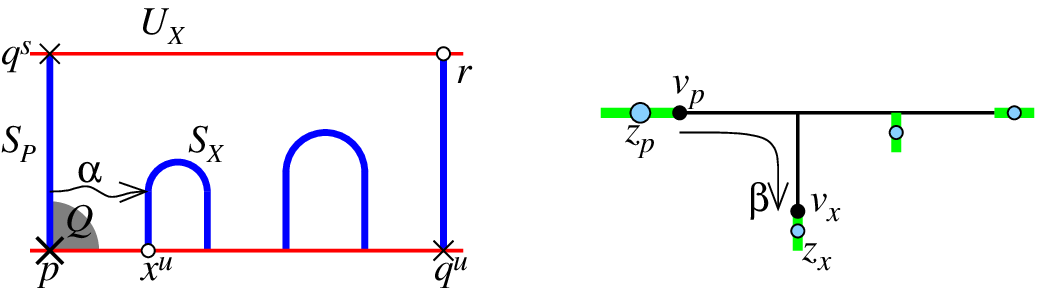}{A trellis and its graph representative in the neighbourhood of a periodic control edge $z_p$.}
To show that we can take non-minimal iterates without increasing the Nielsen entropy, we consider a length function on a graph near a periodic control edge, as shown in \figref{quadrantgraph}.
We let $x^u$ be the first intersection of $T^U(Q)$ with $T^S$, and $S_X(Q)$ be the stable segment on the $Q$-side of $T^U(Q)$ and $x^u$, as shown.
Let $z_p$ be the control edge crossing $S_P(Q)$, and $z_x$ the control edge crossing $S_X(Q)$.
We let $\beta(Q)$ be the curve from the periodic control vertex $v_p$ to the periodic control vertex $v_x$.
Note that $\alpha$ is homotopic to the inclusion of the path $z_p\beta z_x$ in $M$.

Note that for any graph representative of a trellis mapping class, we can choose a length function $l$ such that $l(z)=l(p)=l_p$ for all control edges $z$ and all peripheral edges $p$, and such that $l(g(e))<\lambda_{\epsilon}l(e)$ for all edges $e$.

\begin{lemma}
\label{lem:stableiteratelength}
Let $([f_i];T_i)$ be the minimal extension of $([f];T)$ with $([f_0];T_0)=([f];T)$ and $T_i=(T^U,f_i^{-1}(T_{n-1}))$.
Let $(g_i;G_i,W_i)$ be the graph representative of $([f_i];T_i)$ with natural embeddings $(G_{i-1};W_{i-1})\embed(G_{i},W_{i})$, and let $l_i$ be a length function on $(G_i;W_i)$ preserved by this embedding.
Let $Q$ be a quadrant of $([f];T)$, let $z_p$ be the control edge crossing $S_P(Q)$, and let $z_{i}$ be the control edge crossing the first stable segment crossing $T^U(Q)$.
Let $\beta_i$ be the edge-path in $G_i$ with initial vertex at $v_p$ and final vertex at $v_{i}$ in $z_{i}$.
Then $l_i(\beta_i)\tendsto 0$ as $i\tendsto\infty$, and $l(z_i)\tendsto0$ as $i\tendsto\infty$.
\end{lemma}

\begin{proof}
The result follows since $g_i(\beta_i)=\beta_{i-1}$, and $\beta_i$ does not contain any periodic control edges.
\end{proof}

\begin{lemma}
\label{lem:entropyincreasingextension}
Let $Q_0$ is a quadrant with $f(Q_0)=Q_1$ and $2l(\beta(Q_1))+2l(z_x(Q_1))<(\lambda_{\epsilon}-1)l_p$
Then there is a well-formed trellis mapping class $([\widehat{f}];\widehat{T})$ such that $\hniel([\widehat{f}];\widehat{T})<\lambda_\epsilon$ and $U_X(Q)$ intersects $\widehat{f}^{-1}(S_X(f(Q)))$.
\end{lemma}

\begin{proof}
To obtain the graph $\widehat{G}$, subdivide the control edge $z_p(Q_0)$ into two control edges $\widehat{z}_p$ and $\widehat{z}_q$, and two expanding edges, one edge $\widehat{e}_p$ joining $\widehat{p}$ to $\widehat{q}$, and one edge $\widehat{e}_q$ joining $\widehat{q}$ to $v_p(Q_0)$ as shown in \figref{quadrantbifurcation}.
\fig{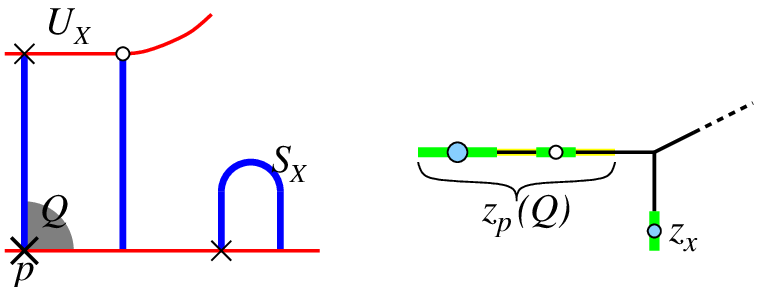}{A non-minimal extension.}

We then recursively subdivide control edges mapping to $z_p$ into a control edge eventually mapping to $\widehat{z}_p$ and an expanding edge.
A compatible graph map is obtained by taking $\widehat{e}_p$ to $\beta(Q_1)$, $\widehat{z}_q$ to $z_x(Q_1)$, and $\widehat{e}_q$ to $\bar{z}_x(Q_1)\bar{\beta(Q_1)}$.
Since this graph map takes $z_p(Q_0)$, with length $l_p$, to $z_p(Q_1)\beta(Q_1)z_x(Q_1)\bar{z}_x(Q_1)\bar{\beta}(Q_1)$ with length less that $\lambda_{\epsilon}l_p$, but is unchanged on all other edges, the growth rate is still less than $\lambda_{\epsilon}$.
Hence the graph representative of $([\widehat{f}];\widehat{T})$, which can be obtained by Algorithm~\ref{alg:main}, has growth rate less than $\lambda_{\epsilon}$.
\end{proof}

We now use these technical results to show how to create intersections bounding regular domains.
\begin{lemma}
\label{lem:transitiveextensionincrease}
Let $([f];T)$ be a trellis mapping class and $Q$ a quadrant of $([f];T)$.
Let $q^u\in T^U(Q)$ be the endpoint of an unstable segment $S_0$ on the same side of $T^S(Q)$ as $Q$, and let $q^s\in T^S(Q)$ be the endpoint of a stable segment $U$ on the same side of $T^S(Q)$ as $Q$.
Further, for any $\epsilon>0$ there is a minimal stable extension $([\widehat{f}];\widehat{T})$ and an integer $k$ such that $\widehat{f}^{-kn}(S)\subset \widehat{T}^S$ and intersects $U$ at a point $r$ such that $\{p,f^{-kn}(q^u),q^s,r\}$ are the vertices of a regular domain for $Q$.
\end{lemma}

\begin{proof}
Let $S_i$ be the segment of $f_{in}^{-k}(S_0)$ containing $f_i^{-in}(q^u)$, let $z_i$ be the control edge crossing $S_i$.
Let $z_p$ be the control edge crossing $S_P(Q)$, and let $\beta_i$ be the edge path from $z_p$ to $z_i$ which contains $z_i$ as its final edge, but does not contain $z_p$.
Then the hypotheses of the length functions for the graph representative imply that $l(\beta_i)\tendsto0$ as $k\tendsto\infty$.
Fix $k$ such that $1+2l(\beta_k)<\lambda_\epsilon$.
\end{proof}

Before proving the main theorem we show how to ensure that every quadrant $Q$ is contained in a regular domain whose opposite sides bound attractors or repellors.
\begin{lemma}
Let $([f];T)$ be a transitive trellis mapping class, and $P$ be an essential periodic orbit of $([f];T)$ which does not shadow $T^P$.
Then there is a minimal supertrellis $([\widetilde{f}];\widetilde{T})$ of $([f];T)$ obtained by blowing up at $P$ such that every quadrant $Q$ is contained in a regular rectangular region with the opposite sides not in $T$.
\end{lemma}

\begin{proof}
Since $P$ is an essential periodic orbit, every curve $\alpha(Q)$ eventually maps across $P$, so crosses a stable segment of $P$.
Hence iterating this segment backwards gives a curve bounding a regular domain, and further backward iterates give a regular region.
A similar argument holds for forward iterates of unstable segments, since in this case we can simply reverse time.
\end{proof}

\begin{theorem}[Existence of approximate entropy minimisers]
\label{thm:minimisingsequence}
Let $([f];T)$ be a well-formed trellis mapping class.
Then for every $\epsilon>0$, there exists a diffeomorphism $\widehat{f}\in([f];T)$ such that $\htop(\widehat{f})<\hniel[f;T]+\epsilon$.
\end{theorem}

\begin{proof}
We repeatedly construct supertrellises $([f_i];T_i)$ with $\hniel[f_i;T_i]<\hniel[f;T]+\epsilon=\log\lambda_{\epsilon}$.

Take $([f_0];T_0)=([f],T)$.
By irreducibility, there is a minimal extension $([f_1];T_1)$ of $([f_0];T_0)$ such that for every pair of nontrivial branches $T_1^U(p_u,b_u)$ and $T_1^S(p_s,b_s)$,
 there are points $p_u=p_0,p_1,\ldots,p_n=p_s$ such that $T_1^U(p_u,b_u)\cap T_1^S(p_1)\neq\emptyset$, $T_1^U(p_i)\cap T_1^S(p_{i+1})\neq\emptyset$ for $1\leq i<n-1$
 and $T_1^U(p_{n-1})\cap T_1^S(p_s,b_s)\neq\emptyset$.
By Lemma~\ref{lem:transitiveextensionincrease}, there is a non-minimal extension $([f_2];T_2)$ of $([f_1];T_1)$ such that every quadrant of $T_2$ is contained in a regular rectangular region.

We next construct a non-minimal extension $([f_3];T_3)$ of $([f_2];T_2)$ such that every nontrivial branch of $T_3^U$ intersects every nontrivial branch of $T_3^S$.
If $T^U(Q_U)$ intersects $T^S(Q)$ and $T^U(Q)$ intersects $T^S(Q_S)$, we can take minimal iterates of $T^U(Q_U)$ and $T^S(Q_S)$ until $T^U(Q_U)\supset U_X(Q)$ and $T^S(Q_S)\supset S_X(Q)$.
Then an application  of Lemma~\ref{lem:transitiveextensionincrease} shows that there is a non-minimal intersection such that $T^U(Q_U)\cap T^S(Q_S)\neq\emptyset$.

Now construct a non-minimal extension $([f_4];T_4)$ of $([f_3];T_3)$ such that there is an alpha-chain from any quadrant $Q_U$ of $T_4$ to any other quadrant $Q_S$.
Since $T_3^S(Q_S)$ intersects $T_3^U(Q_U)$, by taking minimal backward iterates of $T_3^S(Q_S)$ we can ensure that $S_X(Q_U)\subset T^S(Q_S)$.
Then by Lemma~\ref{lem:transitiveextensionincrease}, there is a non-minimal extension such that$T^S(Q_S)$ crosses a regular region $R(Q_U)$.
Similarly, we can ensure $T^U(Q^U)$ crosses the regular region $R(Q_S)$.
A further application of Lemma~\ref{lem:transitiveextensionincrease} to $S_X(Q_U)$ gives a non-minimal extension with an alpha-chain from $Q_U$ to $Q_S$.
Repeating this construction for all regions gives the required trellis $T_4$.

Let $([f_5];T_5)$ be a minimal supertrellis of $([f_4];T_4)$ such that all trivial branches of $T_4$ are contained in a non-chaotic region, as given by Theorem~\ref{thm:minimalsupertrellis}.
If $([f_5];T_5)$ has no attractors or no repellors, take a further minimal supertrellis $([f_6];T_6)$, with at least on attractor and one repellor, as given by Theorem~\ref{thm:minimalsupertrellis} and Theorem~\ref{thm:minimalsupertrellis}.
Since $([f_6];T_6)$ is transitive, we can take a minimal extension $([f_7];T_7)$ by iterating the stable and unstable curves bounding stable and unstable regions such that all quadrants of $T_7$ are contained in a regular region with the opposite stable side bounding an attracting domain, and the opposite unstable side bounding a repelling domain.

\fig{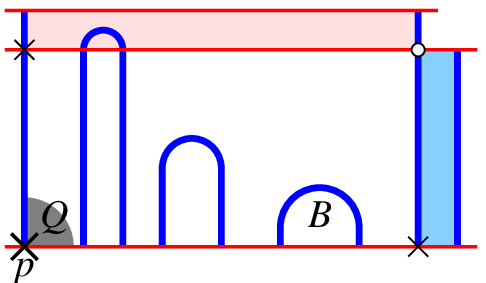}{Backward iterates of $B$ give an inner bigon in a repelling domain.}
Now, for any inner bigon $B$ such that $B^S$ is nonwandering, we take backward iterates of $B^S$ to obtain a stable segment $S$ in a regular domain of a quadrant $Q$.
Applying Lemma~\ref{lem:transitiveextensionincrease} gives a non-minimal extension such that $f^{-n}(B)$ contains a single inner bigon, and this inner bigon lies in a repelling domain, as shown in \figref{nonwanderingbigon}.
Similarly, by taking forward iterates of $B^U$ we can ensure $f^n(B)$ contains a single inner bigon in an attracting domain.
Applying this procedure to all bigons gives a trellis mapping class $([f_8];T_8)$.

Finally, we take a minimal supertrellis $([f_9];T_9)$ satisfying the conditions of Lemma~\ref{lem:simpleminimalmodel} by introducing new stable curves as in Theorem~\ref{thm:minimalsupertrellis}
Then there is a uniformly-hyperbolic diffeomorphism $\widehat{f}\in([f_9];T_9)$, and hence in $([f];T)$ such that $\htop(\widehat{f})=\hniel[f_9;T_9]<\hniel[f;T]+\epsilon$ as required.
\end{proof}

We now prove the existence of pseudo-Anosov map which approximate the Nielsen entropy.
The condition on attractors and repellors is to ensure that the trellis mapping class contains a pseudo-Anosov map, since pseudo-Anosov maps have no attractors or repellors.
The strategy is to create the periodic orbits which will give the one-prong singularities of the pseudo-Anosov map.
Many of the steps of the proof mimic those the proof of Theorem~\ref{thm:minimisingsequence}.
\begin{theorem}[Existence of pseudo-Anosov representatives]
\label{thm:pseudoanosovmodel}
Let $([f];T)$ be a trellis mapping class with no attractors or repellors.
Then for any $\epsilon>0$ there exists a pseudo-Anosov diffeomorphism $\widehat{f}\in([f];T)$ such that $\htop(\widehat{f})<\htop([f];T)+\epsilon$.
\end{theorem}

\begin{proof}
If $T$ has an adjacent pair of trivial branches, take an extension $([f_1];T_1)$ such that these two branches intersect in a single transverse homoclinic point.
This does not affect the Nielsen entropy.
If $T$ has any other trivial branches, take a minimal extension $([f_2];T_2)$ such that these branches have an intersection point, which is possible by Theorem~\ref{thm:minimalsupertrellis} since $T$ has no attracting or repelling regions.
As in the proof of Theorem~\ref{thm:minimisingsequence}, take an extension $([f_3];T_3)$ which is transitive.
Take a further extension $([f_4];T_4)$ such that every branch intersects every other with both orientations, which is possible by transitivity since we can take minimal iterates such that every unstable branch has an extension with and iterate of every stable branch which does not map to an intersection of $T_3$, and an isotopy in the neighbourhood of this intersection yields three intersections, one with the opposite orientation.

\fig{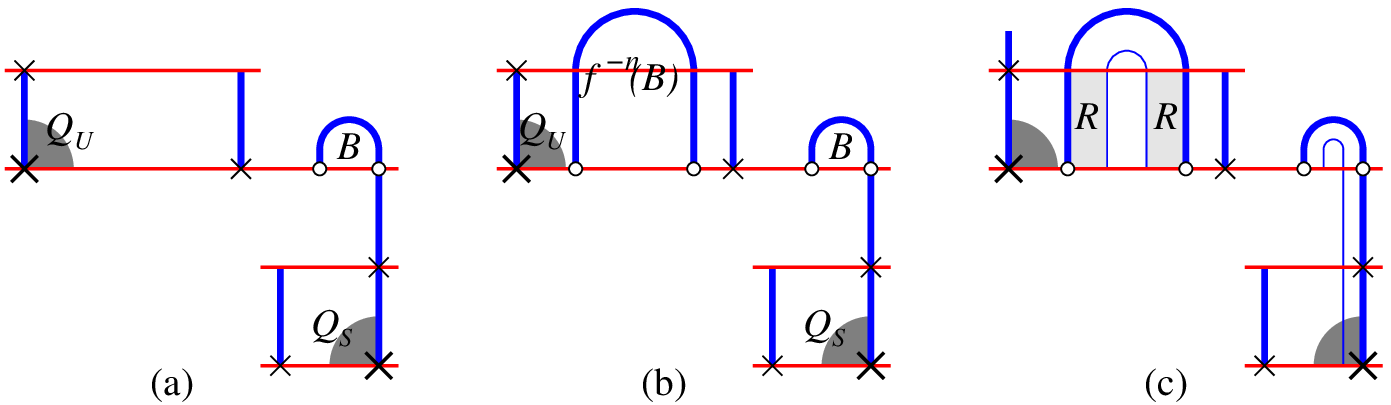}{The domain with vertices $q_0$ and $q_1$ must have periodic orbits in the rectangles $R_0$ and $R_1$.}
Let $B$ be an inner bigon of $T_4$, as shown in \figrefpart{chaoticbigon}{a}.
By Lemma~\ref{lem:transitiveextensionincrease}, we can find a stable extension such that $f{-n}(B^S)$ crosses a regular domain $R(Q_U)$ for some $n$.
Further, by removing intersections if necessary, we can ensure that $f{-n}(B^S)$ crosses $D$ twice and gives a new inner bigon $\widehat{B}$, as shown in \figrefpart{chaoticbigon}{b}.
A further application of Lemma~\ref{lem:transitiveextensionincrease} shows that we can find a stable extension $[\widehat{f}_3;\widehat{T}_3]$ such that $\widehat{B}^U$ crosses some regular domain, as shown in \figrefpart{chaoticbigon}{c}.
Since the regions $R$ are mapped over by $\alpha(Q_U)$ and map over $\alpha(Q_S)$, they must contain a periodic orbit, since $T$ is transitive.
Applying this construction for every inner bigon of $([f_4];T_4)$ gives a non-minimal extension $([f_5];T_5)$ such that every bigon of $T_4$ is a domain of $T_5$ containing an essential periodic orbit of $([f_5];T_5)$.

Let $([f_6];T_6)$ be the trellis mapping class obtained by puncturing at a periodic orbit in every inner bigon of $T^4$ to give a surface $M_6$.
Then $\hniel[f_6;T_6]=\hniel[f_5;T_5]$.
Since $T_4$ is a subtrellis of $T_6$, we can take a trellis mapping class $([f_7];T_7)=([f_6];T_4)$ in the surface $M_6$.
Since every inner bigon of $T_4$ contains component of $\partial M_6$, the trellis mapping class $([f_7];T_7)$ has no inner bigons.
The graph representative $(g_7;G_7,W_7)$ is locally injective except at cusps, so is efficient and hence is a train-track map for a pseudo-Anosov diffeomorphism $\widetilde{f}\in([f_7];T_7)$.
Then $\htop(\widetilde{f})=\hniel[f_7;T_7]\leq\hniel[f_6;T_6]<\hniel[f;T]+\epsilon$ as required.
\end{proof}


\subsection{Non-existence of entropy minimisers}
\label{sec:noentropyminimiser}

Suppose that an irreducible trellis mapping class has a uniformly hyperbolic diffeomorphism realising the Nielsen entropy.
The following lemma shows that any isotopy removing intersections results in a trellis mapping class with strictly smaller Nielsen entropy.
\begin{lemma}
\label{lem:entropydecrease}
Let $([f];T)$ be a well-formed irreducible trellis mapping class, and suppose $f$ be a uniformly-hyperbolic diffeomorphism such that $\htop(f)=\hniel[f;T]$.
Then if $([\overline{f}];\overline{T})$ is a trellis mapping class which does not force $([f];T)$, then $\hniel[\overline{f};\overline{T}]<\hniel[f;T]$.
\end{lemma}
\begin{proof}
\fig{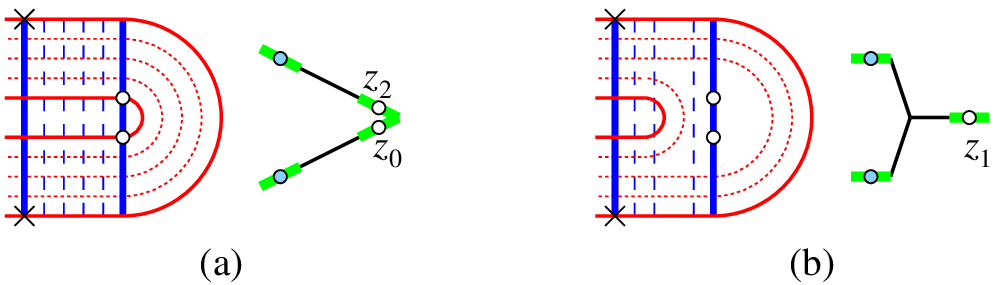}{Change in trellis and graph at a homoclinic bifurcation.}
By taking an $f$-extension $T_1$ of $T$, we can ensure that every inner bigon of $([f];T_1)$ is contained in a larger domain with the topology of \figrefpart{bifurcationdecrease}{a}, and that these domains are separate.
Further, we can ensure that every quadrant is contained in a rectangular region.
Now consider the effect of a homotopy to a trellis mapping class $([f_2];T_2)$ removing one pair of orbits on the same inner bigons.
The resulting trellis and graph locally have the topology of the right of \figrefpart{bifurcationdecrease}{b}.
The edges $a_0$ and $a_2$ must be expanding, since they forward iterate to a closed segment, and since every quadrant is contained in a rectangular region, $\partial g^n(a_0)=\partial g^n(a_2)=b$ for some edge $b$.
Therefore $a_0$ and $a_2$ are folded together to obtain the graph representative for $([f_2];T_2)$.
Further, there must be an edge mapping $\ldots\bar{a}_0a_2\ldots$.
Hence the entropy of the graph representative of $([f_2];T_2)$ is less than that of $[f;T]$, so $\hniel[f_2;T_2]<\hniel[f;T]$.
Since we can then prune $([f_2];T_2)$ to obtain $([\overline{f}];\overline{T})$, we have $\hniel[\overline{f};\overline{T}]\leq\hniel[f_2;T_2]$, hence $\hniel[\overline{f};\overline{T}]<\hniel[f;T]$ as required.
\end{proof}

We use this to prove that the sufficient condition given for the existence of entropy minimisers in Theorem~\ref{thm:minimalmodel} is necessary.
\begin{theorem}[Non-existence of entropy minimisers]
Let $([f];T)$ be an irreducible trellis mapping class.
Suppose that for every diffeomorphism $\widehat{f}\in([f];T)$ there is a $\widehat{f}$-extension of $T$ which is not minimal.
Then there does not exist a diffeomorphism in $([f];T)$ whose topological entropy equals $\hniel[f;T]$.
\end{theorem}

\begin{proof}
Suppose there is an $f$-extension $T_1$ of $T$ which is not minimal.
Let $T_2$ be a transitive $f$-extension of $T_1$, and $T_3$ be a further extension such that every non-wandering segment of $T$ crosses a regular domain of $T_3$.
Then we have entropy bound $\htop(f)\geq\hniel[f;T_3]$.

Take $([f_4];T_4)$ to be a minimal supertrellis of $([f];T_3)$ for which the opposite sides of every regular region bound a stable or unstable region.
Then every nonwandering segment of $T$ enters a stable or unstable region of $T_4$.
Prune $([f_4];T_4)$ to obtain a trellis mapping class $([f_5];T_5)$ which forces $([f];T_1)$ satisfying the conditions of Theorem~\ref{thm:minimalmodel}.
Then $\hniel[f_3;T_3]=\hniel[f_4;T_4]\geq\hniel[f_5;T_5]$, but by Lemma~\ref{lem:entropydecrease}, $\hniel([f_5];T_5)>\hniel[f;T]$.
Combining these inequalities we have $\htop(f)\geq\hniel[f_5;T_5]>\hniel[f;T]$ as required.
\end{proof}

The results of this section show that many fundamental properties of an irreducible trellis type $[f;T]$ depends on whether there is a diffeomorphism $\widetilde{f}\in[f;T]$ for which $\htop(\widetilde{f})=\hniel[f;T]$.
If such a diffeomorphism exists, then the entropy of the trellis type is carried in a uniformly hyperbolic diffeomorphism, but is fragile in the sense that any pruning will reduce the Nielsen entropy, and any diffeomorphism $\widehat{f}$ for which some extension is non-minimal must have strictly greater topological entropy.
If no such diffeomorphism exists, every diffeomorphism in the trellis type has an extension which is non-minimal, but pruning this extension gives a trellis type for which the entropy is still greater than the Nielsen entropy of $[f;T]$.


\begin{appendix}

\section{Intersections and Isotopies of Curves in Surfaces}
\label{sec:isotopytheory}

In this section we prove some technical results on isotopies of curves in surface.
(A \emph{curve} is a function $\alpha:I\fto M$; a \emph{path} is the image of such a function.)
We consider a surface $M$ with a set $\Gamma=\{\gamma_1\ldots\gamma_k\}$ of finitely-many cross-cuts.
(A \emph{cross-cut} is a simple curve or path with both endpoints but no interior points in $\partial M$.) 
All homotopies and isotopies of cross-cuts will be assumed to fix endpoints.
Homotopies and isotopies of curves with endpoints on cross-cuts will be assumed to keep this property.

The main results are Theorem~\ref{thm:minimalnumber},
 which gives conditions under which two sets of curves can be isotoped to reduce the number of interestions with each other,
 Theorem~\ref{thm:twoisotopy}, which shows how to isotope curves in a controlled way, and Theorem~\ref{thm:minimalcurve},
 which shows that there is essentially a unique configuration of curves with minimal intersections in any isotopy class.

We begin by quoting some standard results of (surface) topology. 
The first is due to Epstein \cite{Epstein66}, the can be found in \cite{Hirsch76}, and the third is an easy consequence of the classification of orientable surfaces.
All these results also hold in the differentiable category.
\begin{theorem}
\label{thm:epstein}
Let $\alpha_0$ and $\alpha_1$ be two simple curves in a surface $M$ with the same endpoints. 
Then $\alpha_0$ and $\alpha_1$ are homotopic if and only if they are isotopic (relative to endpoints).
\end{theorem}
\begin{theorem}
\label{thm:curveisotopy} \label{thm:hirsch}
Suppose $\alpha_0$ and $\alpha_1$ are two simple curves in a surface $M$ which are isotopic via an isotopy $\alpha_t$.
Then there is an isotopy $h_t:M\fto M$ with $h_0=\id$ such that $h_t(\alpha_0)=\alpha_s$.
\end{theorem}
\begin{theorem}
\label{thm:crosscut}
Let $M$ be a surface and $\Gamma$ a cross-cut in $M$.
Then if $\alpha_0$ and $\alpha_1$ are curves in $M\setminus\Gamma$ which are homotopic in $M$, then $\alpha_0$ and $\alpha_1$ are homotopic in $M\setminus\Gamma$.
An alternative statement is that the inclusion-induced map $\pi_1(M\setminus\Gamma)\fto \pi_1(M)$ is injective.
\end{theorem}

We also have the following differentiable version of the Alexander trick:
\begin{theorem}
\label{thm:polygondiffeo} \label{thm:polygon}
Suppose $M$ is a polygon, and $f:M\fto M$ is a smooth map which preserves the vertices and edges of $M$.
Then $f$ is isotopic to the identity.
\end{theorem}


\subsection{Minimal intersections}
\label{sec:minimalintersections}

A particularly important concept will be that of two sets of cross-cuts having \emph{minimal intersections} with each other.
We give a definition based on a local definition, and then show that this implies that the two sets have minimise the total number of intersections in the isotopy class.
\begin{definition}[Minimal intersections]
\label{defn:minimalintersections}
A set of mutually disjoint cross-cuts $\Alpha$ has \emph{minimal intersections} with $\Gamma$ if $\Alpha$ and $\Gamma$ are transverse and there is no disc in $M$ bounded by an arc in $\Alpha$ and an arc in $\Gamma$.
\end{definition}

\fig{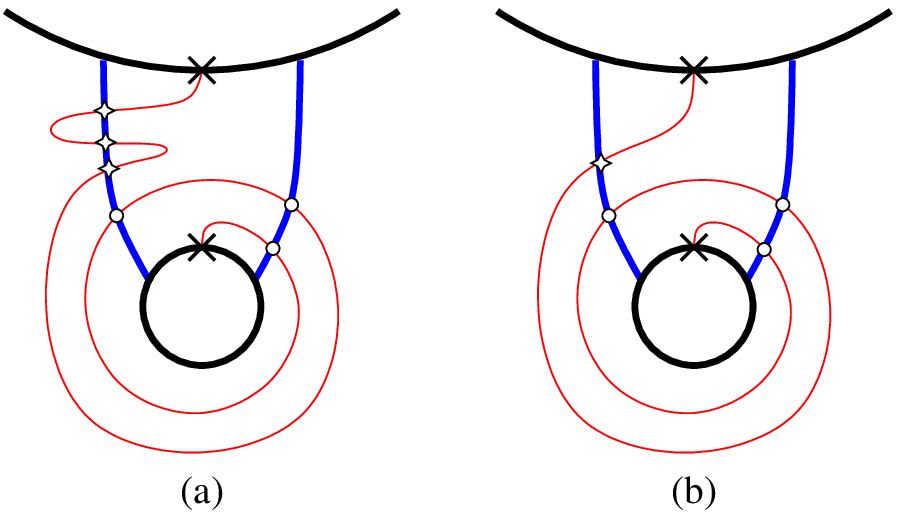}{Intersections of a curve with two cross-cuts. (a) Not minimal intersections. (b) An isotopic curve with minimal intersections.}

The following elementary result relates this local definition to the global property
 of having fewest intersections in the isotopy class.
\begin{theorem}
\label{thm:minimalnumber}
Let $M$ be a surface, $\Gamma$ a set of mutually disjoint cross-cuts in $M$, and $\alpha$ a cross-cut.
Then there is an isotopy $h_t$ such that $h_0=\id$ and $h_1\circ\alpha$ has minimal intersections with $\Gamma$.
Further, $\card{(h_1\circ\alpha)\cap\Gamma} \leq \card{\alpha\cap\Gamma}$,
 with equality if and only if $\alpha$ has minimal intersections with $\Gamma$.
\end{theorem}

The proof is an easy consequence of the following lemma:
\begin{lemma}
\label{lem:minimalintersections}
Suppose $\alpha$ is a curve in $M$, and $D$ a disk bounded by an arc in $\Gamma$ and an arc of $\alpha$, and that the intersection of $\alpha$ and $\gamma$ is (topologically) transverse.
Let $U$ be a simply-connected neighbourhood of $D$ such that $\alpha\cap\Gamma\cap U\subset D$.
Then there is an isotopy $h_t:M\fto M$ supported in $U$ such that $h_1\circ\alpha$ and $\alpha$ agree outside of $U$, but $h_1\circ\alpha$ has no intersections with $\Gamma$ in $U$.
\end{lemma}
\begin{proof}
Let $\alpha[t_1,t_2]$ be the arc of $\alpha$ bounding $D$, and $t_0$, $t_3$ be parameters such that $\alpha[t_0,t_1)$ and $\alpha(t_2,t_3]$ lie in $U\setminus\Gamma$.
Then there is a smooth curve $\alpha$ which agrees with $\alpha$ except in $(t_0,t_2)$, and $\alpha\cap\Gamma\cap U=\emptyset$.
There is a homeomorphism $h_1$ of $M$ supported in $U$ such that $h_1\circ\alpha=\alpha$, and since $U$ is simply-connected, $h_1$ is isotopic to $\id$ by an isotopy supported in $U$.
\end{proof}

\begin{proof}[Theorem~\ref{thm:minimalnumber}]
If the intersections of $\alpha$ with $\Gamma$ are not isolated,
 we take a small perturbation $h_{t_0}$ which is isotopic to the identity
 and for which $h_{t_0}\circ\alpha$ has isolated intersections with $\Gamma$.
In particular, $\alpha_{t_0}$ has finitely many intersections with $\Gamma$ while $\alpha$ has infinitely many.
Now, if $\alpha_{t_0}$ does not have minimal intersections with $\Gamma$,
 there is a disc $D$ bounded by an arc in $\Gamma$ and an arc of $\alpha_{t_0}$,
 so by Lemmas~\ref{lem:minimalintersections}, we can isotope further by $h_t$ to obtain $\alpha_{t_1}=h_{t_1}\circ\alpha$
 with two fewer intersections with $h_{t_0}\circ\alpha$.
(If one of the bounding points of $D$ is an endpoint of $\alpha$, we can still isotope to remove one intersection.)
We continue isotoping until we reach a curve $\alpha_1$ with minimal intersections with $\Gamma$.
The procedure must terminate in a finite number of steps, since $\alpha_{t_0}$ has finitely many intersections.
\end{proof}


\subsection{Ambient isotopies of curves}
\label{sec:isotopy}

\begin{theorem}
\label{thm:twoisotopy}
Let $M$ be a surface and $\Gamma$ a set of cross-cuts in $M$. 
Then if $h$ is isotopic to $\id$ and $h(\Gamma)=\Gamma$
 there is an isotopy $h_t$ such that $h_0=\id$, $h_1=h$ and $h_t(\Gamma)=\Gamma$.
\end{theorem}

\begin{lemma}
\label{lem:disjointisotopy}
Suppose $\alpha_0$ and $\alpha_1$ are homotopic curves which do not intersect $\Gamma$.
Then there is an isotopy $h_t$ such that $h_0=\id$, $h_1(\alpha_0)=\alpha_1$ and $h_t\restrict{\Gamma}=\id$.
\end{lemma}

\begin{proof}
By Theorem~\ref{thm:crosscut}, 
Then $\hat{\alpha}_t=\pi\circ\tilde{\alpha}$ is a homotopy between $\alpha_0$ and $\alpha_1$ in $M\setminus\Gamma$, and hence in a subsurface $K$ of $M$ such that $K\cap\Gamma=\emptyset$.
Then by Theorem~\ref{thm:epstein}, there is an isotopy $\alpha_t$ between $\alpha_0$ and $\alpha_1$ supported in $K$.
By Theorem~\ref{thm:hirsch}, there is an isotopy $h_t$ supported in $K$ such that $h_1(\alpha_0)=\alpha_1$.
\end{proof}

This result can be extended to collections of curves.
\begin{lemma}
\label{lem:collectionisotopy}
Let $\Alpha_0=\{\alpha_0^1\ldots \alpha_0^n\}$ and $\Alpha_1=\{\alpha_1^1\ldots \alpha_1^n\}$ be two sets of disjoint curves in a surface $M$ such that $\alpha_0^i$ is homotopic to $\alpha_1^i$ for $i=1,\ldots,n$.
Then there is an isotopy $h_t$ such that $h_0=\id$ and $h_1(\alpha_0^i)=\alpha_1^i$.
\end{lemma}

\begin{proof}
Assuming the result for $n-1$ curves, we can find an isotopy such that $h_t(\alpha_0^i)=\alpha_1^i$ for $i\leq n-1$, and $h_1(\alpha_0^n)$ and $\alpha_1^n$ are both disjoint from $\alpha_1^1\ldots\alpha_1^{n-1}$ and isotopic in $M$.
The result follows by taking $\Gamma=\{\alpha_1^1,\ldots,\alpha_1^{n-1}\}$ in Lemma~\ref{lem:disjointisotopy}.
\end{proof}

We can now prove Theorem~\ref{thm:twoisotopy}.
\begin{proof}
Let $\Alpha$ be a set of cross-cuts which is disjoint from $\Gamma$ such that $M\setminus(\Alpha\cup\Gamma)$ is simply-connected.
By Lemma~\ref{lem:collectionisotopy}, there is a isotopy $\tilde{h}_t$ such that $\tilde{h}_0=\id$, $\tilde{h}_1(\Alpha)=h(\Alpha)$ and $\tilde{h}_t(\Gamma)=\Gamma$.
Then $\tilde{h}_1^{-1}\circ h(\Gamma)=\Gamma$ and $\tilde{h}_1^{-1}\circ h(\Alpha)=\Alpha$, so $\tilde{h}_1^{-1}\circ h$ is a homeomorhism/diffeomorphism of a polygon, and clearly must preserve the sides of the polygon. 
Hence by Theorem~\ref{thm:polygondiffeo}, there is an isotopy $\hat{h}_t$ such that $\hat{h}_0=\id$, $\hat{h}_1=\tilde{h}_1^{-1}\circ h$ and $\hat{h}_t(\Gamma)=\Gamma$.
Let $h_t=\tilde{h}_t\circ \hat{h}_t$.
Then $h_0=\id$, $h_1=h$ and $h_t(\Gamma)=\tilde{h}+t(\hat{h}_t(\Gamma))=\tilde{h}_t(\Gamma)=\Gamma$ as required.
\end{proof}


\subsection{Nielsen intersection theory}
\label{sec:nielsenintersection}

\begin{theorem}
\label{thm:minimalcurve}
Let $M$ be a surface, $\Gamma$ a set of cross-cuts in $M$. 
Then if $\Alpha_0$ and $\Alpha_1$ are homotopic sets of cross-cuts with minimal intersections with $\Gamma$,
 there is a isotopy $h_t$ such that $h_0=\id$, $h_1(\Alpha_0)=\Alpha_1$ and $h_t(\Gamma)=\Gamma$.
In particular, $h_1$ is a homeomorphism $(M,\Alpha_0,\Gamma)\fto(M,\Alpha_1,\Gamma)$.
\end{theorem}

To prove Theorem~\ref{thm:minimalcurve}, we need to develop a form of Nielsen theory for intersections of curves on surfaces.
We first define an \emph{intersection class}, and then show that each class has a well-defined \emph{intersection number}.
\begin{definition}[Intersection class]
\label{defn:itersectionclass}
Let $\Gamma$ be a set of cross-cuts in a surface $M$, $\alpha$ a curve in $M$, and $x_0$ and $x_1$ points in $\alpha\cap\Gamma$, with $x_i=\alpha(s_i)$ for $i=0,1$.
We say $x_0$ and $x_1$ lie in the same \emph{intersection class} if $\alpha\restrict{[s_0,s_1]}$ is homotopic to a curve in $\Gamma$.
\end{definition}

If $\alpha_t$ is an isotopy, and there are functions $\sigma,\tau:I\fto I$ such that $\alpha_{\tau(t)}(\sigma(t))\in\Gamma$ for all $t$, then we say $\alpha_{\tau(0)}(\sigma(0))$ and $\alpha_{\tau(1)}(\sigma(1))$ are \emph{connected} by the isotopy; essentially they lie in equivalent intersection classes.

\begin{lemma}
Suppose $\alpha_t$ is an isotopy of curves, and there are functions $\sigma$ and $\tau$ such that
 $\alpha_{\tau(0)}(\sigma(0))$ and $\alpha_{\tau(1)}(\sigma(1))$ are connected by the isotopy $\alpha_t$, and $\tau(0)=\tau(1)=\tau$.
Then $\alpha_{\tau}(\sigma(0))$ and $\alpha_\tau(\sigma(1))$ lie in the same intersection class.
\end{lemma}

\begin{proof}
Let 
\[ \hat{\alpha}_t(s)=\alpha_{(1-t)\tau+t\tau(s)}((1-t)((1-s)\sigma(0)+\sigma(1))+t\sigma(s)) .\]
Then $\hat{\alpha}_0(s)=\alpha_\tau((1-s)\sigma(0)+\sigma(1))$, $\hat{\alpha}_1(s)=\alpha_{\tau(s)}(\sigma(s))$ and $\hat{\alpha}_t(i)=\alpha_{\tau}(\sigma(i))$ for $i=0,1$.
\end{proof}

\begin{lemma}
\label{lem:intersectrelation}
Suppose $\alpha_t$ is a homotopy, and there are functions $\sigma_i,\tau_i:I\fto I$ for $i=0,1$ such that $\alpha_{\tau_i(t)}(\sigma_i(t))\in\Gamma$ for $i=1,2$.
Then the curves $\alpha_i\restrict{[\sigma_i(0),\sigma_i(1)]}$ are homotopic via curves with endpoints in $\Gamma$.
\end{lemma}

\begin{proof}
A homotopy is given by 
 \[ \hat{\alpha}_t(s)=\alpha_{(1-s)\tau_0(t)+s\tau_1(t)}((1-s)\sigma_0(t)+\sigma_1(t)) ,\] since we can check that
$\hat{\alpha}_i(s)=\alpha_i((1-s)\sigma_0(i)+s\sigma_1(i))=\alpha_i\restrict{[\sigma_i(0),\sigma_i(1)]}(s)$ and $\hat{\alpha}_t(i)=\alpha_{\tau_i(t)}(\sigma_i(t))\in\Gamma$.
\end{proof}

It is clear that for a given curve, the intersection classes are open in the intersection set.
\begin{definition}[Intersection number]
The \emph{intersection number} of an intersection class is the topological intersection number of that class.
An intersection class is \emph{essential} if its intersection number is non-zero.
\end{definition}

We now show that essential intersection classes preserve their relative ordering under homotopies of curves.
\begin{lemma}
\label{lem:intersectionorder}
Let $\alpha_0$ and $\alpha_1$ be isotopic curves in $M$ with endpoints in $\Gamma$ and minimal intersections with $\Gamma$.
Let $t_0^k$ and $t_1^k$ be ordered parameter values of the intersecions of $\alpha_0$ and $\alpha_1$ with $\Gamma$ respectively.
Then for each curve $\gamma_i$, there are integers $k_1\ldots k_{n_i}$ such that the intersections of $\alpha_0$ and $\alpha_1$ with $\gamma_i$ are $t_0^{k_1}\ldots t_0^{k_{n_i}}$ and $t_0^{k_1}\ldots t_0^{k_{n_i}}$ respectively.
\end{lemma}

\begin{proof}
Since $\alpha_0$ and $\alpha_1$ have minimal intersections with $\Gamma$,
 any two intersections of $\alpha_0$ or $\alpha_1$ with $\Gamma$ are transverse and lie in different intersection classes.
Let $\alpha_t$ be an isotopy between $\alpha_0$ and $\alpha_1$, let $X=\{s,t:\alpha_t(s)\in\Gamma\}$,
 and let $X_t=\{s:\alpha_t(s)\in\Gamma\}$.
Since points in the same component of $X$ with the same $t$-value must correspond to points in the same intersection class of $\alpha_t$, for each intersection point of $\alpha_0$ there is a connected subset of $X$ which contains a point of $X_t$ for any $t$;
 in particular, it contains a point of $X_1$, and this point must be unique
 since $\alpha_1$ also has minimal intersections with $\Gamma$.
This gives the required bijection from $\alpha_0\cap\Gamma$ to $\alpha_1\cap\Gamma$.
\end{proof}

\begin{lemma}
\label{lem:neighbourhoodmap}
Let $\alpha_0$ and $\alpha_1$ be isotopic curves in $M$ with endpoints in $\Gamma$ and minimal intersections with $\Gamma$.
Parameterise $\alpha_0$ and $\alpha_1$ so that the intersections occur at the same parameter values $t_1,t_2,\ldots t_n$.
Then there are tubular neighbourhoods $U,V$ of $\Gamma$ with $\cl{U}\subset V$, and a homeomorphism supported on $V$, such that $h(\alpha_0)$ and $\alpha_1$ agree on $U$, and consist of curves, each crossing $\Gamma$ once transversely.
Further, on each interval $[t_{i},t_{i+1}]$ for $i<n$ and for the intervals $[0,t_1]$ and $[t_n,1]$, the curves $h\circ\alpha_0$ and $\alpha_1$ are homotopic.
\end{lemma}

\begin{proof}
By Lemma~\ref{lem:intersectionorder}, we can find a $U$, $V$ and $h$ such that $h(\alpha_0)$ and $\alpha_1$ agree on $U$, so need only show that $h(\alpha_0)$ and $\alpha_1$ are consist of mutually homotopic curves outside $U$.
By Lemma~\ref{lem:intersectrelation}, the curves are homotopic on each interval.
\end{proof}

Using these results, we can prove Theorem~\ref{thm:minimalcurve}.
\begin{proof}
By Lemma~\ref{lem:intersectionorder} we can find an isotopy $h_t$ such that $h_1(\alpha_0)$ and $\alpha_1$ agree on a neighbourhood $U$ of $\Gamma$.
Without loss of generality, therefore, we can assume that $\alpha_0$ and $\alpha_1$ agree on $U$.
Now let $V$ be a tubular neighbourhood of $\Gamma$ such that $\cl{V}\subset U$, and $K$ be the surface $M\setminus V$.
Then $\alpha_0\cap K$ and $\alpha_1\cap K$, give families of mutually disjoint homotopic curves, and the result follows from Lemma~\ref{lem:collectionisotopy}.
\end{proof}

\end{appendix}



\end{document}
